\documentclass[11pt]{amsproc}
\usepackage{amssymb,amsmath,amsthm,anysize}
\usepackage{longtable}
\usepackage{booktabs}
\usepackage{multirow,bigstrut,bigdelim,color}
\usepackage{lscape}

\usepackage{rotating}
\newtheorem{Thm}{Theorem}[section]
\newtheorem{Cor}[Thm]{Corollary}

\newtheorem{Pro}[Thm]{Proposition}
\newtheorem{La}[Thm]{Lemma}

\newenvironment{Prf}{\noindent\textbf{Proof.}}{\hfill $\Box$ \medskip}

\newcommand{\cV}{\mathcal{V}}
\newcommand{\cC}{\mathcal{C}}

\newcommand{\bF}{\mathbb{F}}

\newcommand{\cZ}{\mathcal{Z}}
\newcommand{\cN}{\mathcal{N}}

\newcommand{\GL}{\mathrm{GL}}

\newcommand{\PGL}{\mathrm{PGL}}

\newcommand{\GO}{\mathrm{GO}}

\newcommand{\PG}{\mathrm{PG}}

\renewcommand\le{\leqslant}

\title{Combinatorial invariants for nets of conics in $\PG(2,q)$}

\author{Michel Lavrauw, Tomasz Popiel, John Sheekey}
\address{Michel Lavrauw, Sabanc{\i} University, Istanbul, Turkey; Email: \textsc{mlavrauw@sabanciuniv.edu}}
\address{Tomasz Popiel; Email: \textsc{tomasz.popiel@uwa.edu.au}}
\address{John Sheekey, University College Dublin, Ireland; Email: \textsc{john.sheekey@ucd.ie}}
\date{\today}

\begin{document}

\begin{abstract}
The problem of classifying linear systems of conics in projective planes dates back at least to Jordan, who classified pencils (one-dimensional systems) of conics over $\mathbb{C}$ and $\mathbb{R}$ in 1906--1907. 
The analogous problem for finite fields $\mathbb{F}_q$ with $q$ odd was solved by Dickson in 1908.  In 1914, Wilson attempted to classify nets (two-dimensional systems) of conics over finite fields of odd characteristic, but his classification was incomplete and contained some inaccuracies. 
In a recent article, we completed Wilson's classification of nets of {\em rank one}, namely those containing a repeated line. 
The aim of the present paper is to introduce and calculate certain combinatorial invariants of these nets, which we expect will be of use in various applications. 
Our approach is geometric in the sense that we view a net of rank one as a plane in $\PG(5,q)$ that meets the quadric Veronesean in at least one point; two such nets are then equivalent if and only if the corresponding planes belong to the same orbit under the induced action of $\PGL(3,q)$ viewed as a subgroup of $\PGL(6,q)$. 
We have previously determined the orbits of {\em lines} in $\PG(5,q)$ under this action, which correspond to the aforementioned pencils of conics in $\PG(2,q)$. 
The main contribution of this paper is to determine the {\em line-orbit distribution} of a plane $\pi$ corresponding to a net of rank one, namely, the number of lines in $\pi$ belonging to each line orbit. 
It turns out that this list of invariants completely determines the orbit of $\pi$, and we will use this fact in forthcoming work to develop an efficient algorithm for calculating the orbit of a given net of rank one. 
As a more immediate application, we also determine the stabilisers of nets of rank one in $\PGL(3,q)$, and hence the orbit sizes. 
\end{abstract}

\maketitle

\section{Introduction}

The forms of degree $d$ on an $n$-dimensional projective space $\PG(V)$ comprise a vector space $W$ of dimension ${n+d}\choose{d}$. 
Subspaces of the projective space $\PG(W)$ are called {\em linear systems of hypersurfaces of degree $d$}. 
One-dimensional linear systems are called a {\em pencils}, and two-dimensional linear systems are called {\em nets}. 
In this paper we are concerned with linear systems of conics, namely the case $d=n=2$. 
For an elementary exposition of the more general case $d=2$ (that is, pencils of quadrics), we refer the reader to \cite[Chapters 9 and 11]{Casas-Alvero2010}.

The classification of linear systems of conics in $\PG(2,\mathbb{F})$ consists of determining the orbits of the subspaces of $\PG(5,\mathbb{F})$ under the induced action of the projectivity group $\PGL(3,\mathbb{F})$. 
Pencils of conics over the fields $\mathbb{F}=\mathbb{C}$ and $\mathbb{R}$ were classified by Jordan~\cite{Jordan1906, Jordan1907} in 1906--1907, and nets of conics over these fields were treated several decades later by Wall~\cite{Wall1977}. 
Here we are concerned with linear systems of conics over {\em finite} fields. 
Pencils of conics over $\mathbb{F}_q$ with $q$ odd were classified by Dickson~\cite{Dickson1908} in 1908, and an incomplete classification for $q$ even was obtained by Campbell~\cite{Campbell1927} in 1927. 
In 1914, Wilson~\cite{Wilson1914} obtained an incomplete classification of nets of conics for $q$ odd. 
In a recent paper~\cite{LaPoSh2020}, we completed Wilson's classification of nets of conics of {\em rank one}, namely those containing a repeated line. 
It turns out that there are $15$ orbits of such nets for every odd $q$. 
Representatives are given in Table~\ref{table:main}, the notation of which is explained in Section~\ref{s:prelim}.

Wilson's approach to this classification problem was purely algebraic, based on explicit coordinate transformations. 
Our approach was geometric, based on the observation that a net of conics in $\PG(2,q)$ corresponds to a plane in $\PG(5,q)$, and more specifically that a net of rank one corresponds to a plane that meets the quadric Veronesean in at least one point. 
The details of this correspondence, and of the aforementioned action of $\PGL(3,q)$ on subspaces of $\PG(5,q)$, are explained in Section~\ref{s:prelim}. 
In an earlier paper \cite{LaPo2020}, we had also investigated the classification of {\em pencils} of conics from this geometric viewpoint, determining the corresponding orbits of {\em lines} in $\PG(5,q)$. 
Although the classification itself was necessarily in agreement with that of Dickson~\cite{Dickson1908}, our geometric approach had certain advantages, and in particular allowed us to calculate various auxiliary data about the line orbits, including the {\em point-orbit distribution} of a line $\ell$, namely the number of points on $\ell$ belonging to each of the four $\PGL(3,q)$-orbits of points in $\PG(5,q)$. 
Indeed, our original motivation was not to reproduce Dickson's classification, but was instead related to the classification in \cite{LaSh2015} of tensors in $\mathbb{F}_q^2 \otimes \mathbb{F}_q^3 \otimes \mathbb{F}_q^3$. 
The point-orbit distributions obtained in \cite{LaPo2020} were useful for several arguments in \cite{LaPoSh2020}, and have since been used to develop an algorithm \cite{AlLa2019} for determining the tensor rank of an element of $\mathbb{F}_q^2 \otimes \mathbb{F}_q^3 \otimes \mathbb{F}_q^3$. 
More generally, of course, the problem of calculating the rank of a tensor has many wide-ranging applications, including, famously, to determining the computational complexity of matrix multiplication.

The aim of the present paper is to calculate analogous data for nets of conics of rank one in $\PG(2,q)$, viewed as planes in $\PG(5,q)$ intersecting the quadric Veronesean. 
The point-orbit distributions of the planes corresponding to nets of rank one were determined in \cite{LaPoSh2020}; they are recorded here in Table~\ref{table:line-rank-dist}, the notation of which is explained in Section~\ref{s:prelim}. 
It is natural to ask also about the {\em line-orbit distribution} of a plane $\pi$, that is, the number of lines in $\pi$ belonging to each of the orbits determined in \cite{LaPo2020}. 
These data are determined in Section~\ref{s:dist}, which comprises the bulk of the paper and constitutes the proof of the following theorem. 
In order to state this theorem in terms of nets rather than planes, we abuse terminology and take the ``line-orbit distribution'' of a net to mean the line-orbit distribution of the corresponding plane.

\begin{Thm} \label{thm:dists}
The line-orbit distributions of nets of conics of rank one in $\PG(2,q)$, $q$ odd, are given in Table~\ref{fullTable}.
\end{Thm}

It turns out that the line-orbit distribution of a given net of rank one determines the orbit of the net. 
In fact, for this purpose it suffices to know which line orbits intersect the corresponding plane (that is, we do not need to know the exact number of lines of each type). 
This information is summarised in Table~\ref{TableLODsummary}, which, upon inspection, implies the following corollary. 

\begin{Cor} \label{cor:complete}
If two nets of conics of rank one in $\PG(2,q)$, $q$ odd, have the same line-orbit distribution, then they belong to the same orbit under the natural action of $\PGL(3,q)$.
\end{Cor}

This suggests a natural application, which we will address in a forthcoming paper, namely, to develop an efficient algorithm for determining the orbit of a given net of rank one. 
We expect that Theorem~\ref{thm:dists} will have also have other applications. 
Here we use it to aid in calculating the stabiliser subgroups of nets of rank one in $\PGL(3,q)$. 
We obtain the following result, which we state here in terms of nets but prove in Section~\ref{s:stabilisers} in terms of planes. 

\begin{Thm} \label{thm:stabs}
The stabilisers in $\PGL(3,q)$ and orbit sizes of nets of rank one in $\PG(2,q)$, $q$ odd, are given in Table~\ref{TableStabilisers}.
\end{Thm}

Note that the stabilisers of pencils (equivalently, lines in $\PG(5,q)$) were determined in \cite{LaPo2020}. 
They are recorded here in Table~\ref{table:line-rank-dist} for reference.

For further discussion of the history of the problem of classifying linear systems of conics, we refer the reader to \cite[Section 1.2]{LaPo2020}, where we also outline the technical difficulties that arise when working over a finite field as opposed to an algebraically closed field.

\begin{sidewaystable}
\small
\begin{tabular}{lrccccccccccccccc}
\toprule
Orbit & $q\;(3)$ & $o_5$ & $o_6$ & $o_{8,1}$ & $o_{8,2}$ & $o_9$ & $o_{10}$ & $o_{12}$ & $o_{13,1}$ & $o_{13,2}$ & $o_{14,1}$ & $o_{14,2}$ & $o_{15,1}$ & $o_{15,2}$ & $o_{16}$ & $o_{17}$ \\ 
\midrule
$\Sigma_1$ & & $\tfrac{q(q+1)}{2}$ & $q+1$ &&&& $\tfrac{q(q-1)}{2}$ &&&&&&&&& \\
$\Sigma_2$ & & $3$ && $\tfrac{3(q-1)}{2}$ & $\tfrac{3(q-1)}{2}$ &&&&&& $\tfrac{(q-1)^2}{4}$ & $\tfrac{3(q-1)^2}{4}$ &&&& \\
$\Sigma_3$ & & $1$ & $1$ & $q$ && $q-1$ &&& $\tfrac{q(q-1)}{2}$ & $\tfrac{q(q-1)}{2}$ &&&&&& \\
$\Sigma_4$ & & $1$ && $2q$ &&&& $1$ & $\tfrac{q(q-3)}{2}$ & $\tfrac{q(q-1)}{2}$ &&&&& $q-1$ & \\
$\Sigma_5$ & & $1$ && $q-1$ & $q-1$ & $2$ &&& $\tfrac{q-1}{2}$ & $\tfrac{q-1}{2}$ & $\tfrac{(q-1)(q-3)}{8}$ & $\tfrac{(q-1)(3q-5)}{8}$ & $\tfrac{(q-1)^2}{4}$ & $\tfrac{(q+1)(q-1)}{4}$ && \\
$\Sigma_6$ & &&& $\tfrac{q+1}{2}$ & $\tfrac{q+1}{2}$ && $1$ &&&&&& $\tfrac{(q+1)(q-1)}{2}$ & $\tfrac{(q+1)(q-1)}{2}$ && \\
$\Sigma_7$ & && $q+1$ &&&&& $q^2$ &&&&&&&& \\
$\Sigma_8$ & && $1$ & $q$ &&&& $1$ & $q(q-1)$ &&&&&& $q-1$ & \\
$\Sigma_9$ & && $1$ & $q$ &&&&& $q$ && $\tfrac{q(q-1)}{2}$ && $\tfrac{q(q-1)}{2}$ && \\
$\Sigma_{10}$ & && $1$ && $q$ &&&&& $q$ && $\tfrac{q(q-1)}{2}$ & $\tfrac{q(q-1)}{2}$ && \\
$\Sigma_{11}$ & $0$ &&& $q$ && $1$ &&& $q$ && $\tfrac{q(q-3)}{6}$ && $\tfrac{q(q-1)}{2}$ && & $\tfrac{q^2}{3}$ \\
 & $\not \equiv 0$ &&& $q$ && $1$ &&& $q-1$ && $\tfrac{(q-1)(q-2)}{6}$ && $\tfrac{q(q-1)}{2}$ && $1$ & $\tfrac{(q+1)(q-1)}{3}$ \\
$\Sigma_{12}$ & $1$ &&& $\tfrac{q-1}{2}$ & $\tfrac{q-1}{2}$ & $2$ &&& $\tfrac{q-7}{2}$ & $\tfrac{q-1}{2}$ & $\tfrac{(q-1)(q-7)}{24}+1$ & $\tfrac{(q-1)(q-3)}{8}$ & $\tfrac{(q-1)^2}{4}$ & $\tfrac{(q+1)(q-1)}{4}$ & $3$ & $\tfrac{(q-1)(q+2)}{3}$ \\
 & $\not \equiv 1$ &&& $\tfrac{q-1}{2}$ & $\tfrac{q-1}{2}$ & $2$ &&& $\tfrac{q-3}{2}$ & $\tfrac{q-1}{2}$ & $\tfrac{(q-3)(q-5)}{24}$ & $\tfrac{(q-1)(q-3)}{8}$ & $\tfrac{(q-1)^2}{4}$ & $\tfrac{(q+1)(q-1)}{4}$ & $1$ & $\tfrac{q(q+1)}{3}$ \\
$\Sigma_{13}$ & $-1$ &&& $\tfrac{q+1}{2}$ & $\tfrac{q+1}{2}$ &  &&& $\tfrac{q-5}{2}$ & $\tfrac{q+1}{2}$ & $\tfrac{(q+1)(q-5)}{24}+1$ & $\tfrac{(q+1)(q-1)}{8}$ & $\tfrac{(q+1)(q-3)}{4}$ & $\tfrac{(q+1)(q-1)}{4}$ & $3$ & $\tfrac{(q+1)(q-2)}{3}$ \\
 & $\not \equiv -1$ &&& $\tfrac{q+1}{2}$ & $\tfrac{q+1}{2}$ &  &&& $\tfrac{q-1}{2}$ & $\tfrac{q+1}{2}$ & $\tfrac{(q-1)(q-3)}{24}$ & $\tfrac{(q+1)(q-1)}{8}$ & $\tfrac{(q+1)(q-3)}{4}$ & $\tfrac{(q+1)(q-1)}{4}$ & $1$ & $\tfrac{q(q-1)}{3}$ \\
$\Sigma_{14}$ & $1$ &&& $\tfrac{q-1}{2}$ & $\tfrac{q-1}{2}$ & $2$ &&& $\tfrac{q-1}{2}$ & $\tfrac{q-1}{2}$ & $\tfrac{(q-1)(q-7)}{24}$ & $\tfrac{(q-1)(q-3)}{8}$ & $\tfrac{(q-1)^2}{4}$ & $\tfrac{(q+1)(q-1)}{4}$ & & $\tfrac{(q-1)^2}{3}+q$ \\
 & $-1$ &&& $\tfrac{q+1}{2}$ & $\tfrac{q+1}{2}$ & &&& $\tfrac{q+1}{2}$ & $\tfrac{q+1}{2}$ & $\tfrac{(q+1)(q-5)}{24}$ & $\tfrac{(q+1)(q-1)}{8}$ & $\tfrac{(q-1)(q-3)}{4}$ & $\tfrac{(q+1)(q-1)}{4}$ & & $\tfrac{(q-1)^2}{3}-q$ \\
$\Sigma_{14}'$ & $0$ &&& $q$ && $1$ &&&&& $\tfrac{q(q-1)}{6}$ && $\tfrac{q(q-1)}{2}$ && $q$ & $\tfrac{q(q-1)}{3}$ \\
$\Sigma_{15}$ & && $1$ &&& $q$ &&&&&&&&& $q^2$ & \\
\bottomrule
\end{tabular}
\caption{Line-orbit distributions of planes in $\PG(5,q)$, $q$ odd, that intersect the quadric Veronesean. Where applicable, the second column indicates the congruence class(es) of $q$ modulo $3$.}
\label{fullTable}
\end{sidewaystable}

\begin{table}[!t]
\begin{tabular}{lccccccccccccccc}
\toprule
Orbit & $o_5$ & $o_6$ & $o_{8,1}$ & $o_{8,2}$ & $o_9$ & $o_{10}$ & $o_{12}$ & $o_{13,1}$ & $o_{13,2}$ & $o_{14,1}$ & $o_{14,2}$ & $o_{15,1}$ & $o_{15,2}$ & $o_{16}$ & $o_{17}$ \\ 
\midrule
$\Sigma_1$ & $\times$ & $\times$ &&&& $\times$ &&&&&&&&& \\
$\Sigma_2$ & $\times$ && $\times$ & $\times$ &&&&&& $\times$ & $\times$ &&&& \\
$\Sigma_3$ & $\times$ & $\times$ & $\times$ && $\times$ &&& $\times$ & $\times$ &&&&&& \\
$\Sigma_4$ & $\times$ && $\times$ &&&& $\times$ & $\times$ & $\times$ &&&&& $\times$ & \\
$\Sigma_5$ & $\times$ && $\times$ & $\times$ & $\times$ &&& $\times$ & $\times$ & $\times$ & $\times$ & $\times$ & $\times$ && \\
$\Sigma_6$ &&& $\times$ & $\times$ && $\times$ &&&&&& $\times$ & $\times$ && \\
$\Sigma_7$ && $\times$ &&&&& $\times$ &&&&&&&& \\
$\Sigma_8$ && $\times$ & $\times$ &&&& $\times$ & $\times$ &&&&&& $\times$ & \\
$\Sigma_9$ && $\times$ & $\times$ &&&&& $\times$ && $\times$ && $\times$ && \\
$\Sigma_{10}$ && $\times$ && $\times$ &&&&& $\times$ && $\times$ & $\times$ && \\
$\Sigma_{11}$ &&& $\times$ && $\times$ &&& $\times$ && $\times$ && $\times$ && $\bullet$ & $\times$ \\
$\Sigma_{12}$ &&& $\times$ & $\times$ & $\times$ &&& $\times$ & $\times$ & $\times$ & $\times$ & $\times$ & $\times$ & $\times$ & $\times$ \\
$\Sigma_{13}$ &&& $\times$ & $\times$ &  &&& $\times$ & $\times$ & $\times$ & $\times$ & $\times$ & $\times$ & $\times$ & $\times$ \\
$\Sigma_{14}$ &&& $\times$ & $\times$ & $\bullet$ &&& $\times$ & $\times$ & $\times$ & $\times$ & $\times$ & $\times$ & & $\times$ \\
$\Sigma_{14}'$ &&& $\times$ && $\times$ &&&&& $\times$ && $\times$ && $\times$ & $\times$ \\
$\Sigma_{15}$ && $\times$ &&& $\times$ &&&&&&&&& $\times$ & \\
\bottomrule
\end{tabular}
\caption{Summary of line-orbit distributions of planes in $\PG(5,q)$, $q$ odd, that intersect the quadric Veronesean. 
Here $\times$ means that a plane of the indicated type contains at least one line of the indicated type, for all $q$, and $\bullet$ means that the existence of a line of the indicated type depends on $q$. 
Specifically: a plane of type $\Sigma_{11}$ contains a line of type $o_{16}$ if and only if $q \not \equiv 0 \pmod 3$, and a plane of type $\Sigma_{14}$ contains a line of type $o_9$ if and only if $q \equiv 1 \pmod 3$. 
}
\label{TableLODsummary}
\end{table}

\begin{table}[!t]
\begin{tabular}{crcc}
\toprule
Orbit & $q \pmod 3$ & Stabiliser & Orbit size \\
\midrule
\addlinespace[2pt]
$\Sigma_1$ & & $E_q^2:\textnormal{GL}(2,q)$ & $q^2+q+1$ \\
\addlinespace[2pt]
$\Sigma_2$ & & $C_{q-1}^2:S_3$ & $\tfrac{1}{6}q^3(q^2+q+1)(q+1)$ \\
\addlinespace[2pt]
$\Sigma_3$ & & $E_q:C_{q-1}^2$ & $q^2(q^2+q+1)(q+1)$ \\
\addlinespace[2pt]
$\Sigma_4$ & & $(E_q:C_{q-1}):C_2$ & $\tfrac{1}{2}q^2(q^3-1)(q+1)$ \\
\addlinespace[2pt]
$\Sigma_5$ & & $C_{q-1} : C_2$ & $\tfrac{1}{2}q^3(q^3-1)(q+1)$ \\
\addlinespace[2pt]
$\Sigma_6$ & & $\GO^-(2,q)$ & $\tfrac{1}{2}q^3(q^3-1)$ \\
\addlinespace[2pt]
$\Sigma_7$ & & $E_q^2:\textnormal{GL}(2,q)$ & $q^2+q+1$ \\
\addlinespace[2pt]
$\Sigma_8$ & & $E_q:C_{q-1}^2$ & $q^2(q^2+q+1)(q+1)$ \\
\addlinespace[2pt]
$\Sigma_9$ & & $(E_q:C_{q-1}):C_2$ & $\tfrac{1}{2}q^2(q^3-1)(q+1)$ \\
\addlinespace[2pt]
$\Sigma_{10}$ & & $(E_q:C_{q-1}):C_2$ & $\tfrac{1}{2}q^2(q^3-1)(q+1)$ \\
\addlinespace[2pt]
$\Sigma_{11}$ & $\not \equiv 0$ & $C_{q-1}$ & $q^3(q^3-1)(q+1)$ \\
\addlinespace[2pt]
& $0$ & $E_q$ & $q^2(q^3-1)(q^2-1)$ \\
\addlinespace[2pt]
$\Sigma_{12}$ & $1$ & $S_3$ & $\tfrac{1}{6}q^3(q^3-1)(q^2-1)$  \\
\addlinespace[2pt]
& $\not \equiv 1$ & $C_2$ & $\tfrac{1}{2}q^3(q^3-1)(q^2-1)$ \\
\addlinespace[2pt]
$\Sigma_{13}$ & $-1$ & $S_3$ & $\tfrac{1}{6}q^3(q^3-1)(q^2-1)$ \\
\addlinespace[2pt]
& $\not \equiv -1$ & $C_2$ & $\tfrac{1}{2}q^3(q^3-1)(q^2-1)$ \\
\addlinespace[2pt]
$\Sigma_{14}$ & $\not \equiv 0$ & $C_3$ & $\tfrac{1}{3} q^3(q^3-1)(q^2-1)$ \\
\addlinespace[2pt]
$\Sigma_{14}'$ & $0$ & $E_q : C_{q-1}$ & $q^2(q^3-1)(q+1)$ \\
\addlinespace[2pt]
$\Sigma_{15}$ & & $E_q^{1+2}:C_{q-1}$ & $(q^3-1)(q+1)$ \\
\bottomrule
\end{tabular}
\caption{Stabilisers and orbit sizes of planes in $\PG(5,q)$, $q$ odd, that intersect the quadric Veronesean, under the action of $\PGL(3,q)$ defined in Section~\ref{s:prelim}. 
Where applicable, the second column indicates the congruence class(es) of $q$ modulo $3$. 
The group-theoretic notation used in the second column is explained in Section~\ref{s:stabilisers}.}
\label{TableStabilisers}
\end{table}

\section{Preliminaries} \label{s:prelim}

Here we summarise some necessary background material, most of which is drawn from our earlier paper \cite{LaPoSh2020}, to which we refer the reader for further details.

We denote the finite field of order $q$ by $\bF_q$, assuming throughout the paper that the prime power $q$ is odd. 
A {\em form} on a vector space $V$ (or on $\PG(V)$) is a homogeneous polynomial in the polynomial ring over the coefficient field of $V$ in $\dim(V)$ variables. The zero locus in $\PG(V)$ of a form $f$ on $\PG(V)$ is denoted by $\cZ(f)$.
A ternary quadratic form $f$ on $\bF_q^3$ defines a {\em conic} $\cC=\cZ(f)$ in $\PG(2,q)$. 
Two distinct conics $\cZ(f_1)$, $\cZ(f_2)$ define a {\em pencil of conics} 
\[\{\cZ(af_1+bf_2):a,b\in \bF_q, (a,b)\neq (0,0)\}.\]
A {\em net of conics} $\cN$ is defined by three conics $\cC_i=\cZ(f_i)$, $i\in\{1,2,3\}$, not contained in a pencil:
\[
\cN=\{\cZ(af_1+bf_2+cf_3):a,b,c\in \bF_q, (a,b,c)\neq (0,0,0)\}.
\]
Given a net, one can consider
\begin{equation} \label{eqn:qform}
xf_1+yf_2+zf_3=a_{00}(x,y,z)X_0^2+a_{01}(x,y,z)X_0X_1+\dots + a_{22}(x,y,z)X_2^2
\end{equation}
as a quadratic form whose coefficients are linear forms in $x,y,z$. 
Given $a,b,c\in \bF_q$, not all zero, we obtain a conic 
\[
\cN(a,b,c)=\cZ(af_1+bf_2+cf_3).
\]
Since $q$ is odd, we can consider the matrix $A_\cN$ of the bilinear form associated to the quadratic form (\ref{eqn:qform}). 
The {\em discriminant} 
\[
\Delta_\cN=\det (A_\cN) 
\]
of the net $\cN$ defines a cubic curve $\cZ(\Delta_\cN)$ in $\PG(2,q)$. 
The conic $\cN(a,b,c)$ is singular if and only if $(a,b,c)$ lies on the cubic $\cZ(\Delta_\cN)$: see \cite[Lemma~2.1]{LaPoSh2020}.
A net has {\em rank one} if it contains a repeated line, {\em rank two} if it contains no repeated lines but contains a conic which is not absolutely irreducible, and {\em rank three} if every conic in the net is absolutely irreducible.

We represent points $y=(y_0,y_1,y_2,y_3,y_4,y_5,y_6)$ of $\PG(5,q)$ by symmetric matrices
\begin{equation} \label{eq:symmetricMatrix}
M_y= \left[ 
\begin{matrix} y_0 & y_1 & y_2 \\ y_1 & y_3 & y_4\\ y_2 & y_4 & y_5\\\end{matrix}
\right].
\end{equation}
The Veronese surface $\cV(\bF_q)$ in $\PG(5,q)$ is defined by setting the $2\times 2$ minors of the above matrix equal to zero, and we have the corresponding Veronese map from $\PG(2,q)$ to $\cV(\bF_q) \subset \PG(5,q)$: 
\[
\nu:(x_0,x_1,x_2)\mapsto (x_0^2,x_0x_1,x_0x_2,x_1^2,x_1x_2,x_2^2). 
\]
The {\em rank} of a point in $\PG(5,q)$ is the rank of the matrix $M_y$. 
The points of rank $1$ are those in $\cV(\bF_q)$; the points of rank $2$ are those in its secant variety. 
The image of a line $\ell$ in $\PG(2,q)$ under $\nu$ is a conic in $\cV(\bF_q)$. 
A plane in $\PG(5,q)$ intersecting $\cV(\bF_q)$ in a conic is called a {\it conic plane}, and each conic plane is equal to $\langle \nu(\ell)\rangle$ for some line $\ell$ in $\PG(2,q)$. 
Each two points $x,y \in \cV(\bF_q)$ lie on a unique conic $\cC(x,y) \subset \cV(\bF_q)$ given by
\[
\cC(x,y)=\nu(\langle \nu^{-1}(x),\nu^{-1}(y)\rangle).
\]
Each rank-$2$ point $z$ lies in a unique conic plane $\langle \cC_z \rangle$. 
If $z$ is on the secant $\langle x,y\rangle$ then $\cC_z = \cC(x,y)$. 

We have a mapping $\delta^*$ from set of conics in $\PG(2,q)$ to the set of points in $\PG(5,q)$, taking the conic $\cC=\cZ(f)$ with $f(X_0,X_1,X_2)=\sum_{i\le j}a_{ij}X_iX_j$ to the point $(a_{00},a_{01},a_{02},a_{11},a_{12},a_{22})$, which may in turn be viewed as a $3 \times 3$ matrix as per \eqref{eq:symmetricMatrix}. 
Under this mapping, a pencil (respectively, net) of conics in $\PG(2,q)$ becomes a line (respectively, plane) in $\PG(5,q)$. 
We consider the action of the group $\PGL(3,q)$ on the points of $\PG(5,q)$ defined as follows: 
if $\varphi_A \in \PGL(3,q)$ is induced by $A\in \GL(3,q)$ then we define $\alpha(\varphi_A) \in \PGL(6,q)$ by
\[
\alpha(\varphi_A):y\mapsto z, \quad \text{where} \quad M_z=AM_yA^T,
\] 
and write
\[
K:=\alpha(\PGL(3,q))\le \PGL(6,q).
\]
The action of $K$ on points of $\PG(5,q)$ induces an action on subspaces, and we are able to make the following crucial observation:

\begin{Pro}
The classification of pencils (respectively, nets) of conics in $\PG(2,q)$, $q$ odd, up to coordinate transformations is equivalent to the classification of $K$-orbits of lines (respectively, planes) in $\PG(5,q)$. 
\end{Pro}

Moreover, since we assume $q$ to be odd, the $K$-orbits of points (respectively, lines) in $\PG(5,q)$ are in one-to-one correspondence with the $K$-orbits of hyperplanes (respectively, solids) in $\PG(5,q)$. 
Let us now describe the orbits of $K$ on points, lines and planes in $\PG(5,q)$. 
The points of rank $1$ in $\PG(5,q)$ form one $K$-orbit, and the points of rank $3$ form a second $K$-orbit. 
There are two $K$-orbits of points of rank $2$, namely the {\em exterior} and {\em interior} rank-$2$ points, where a rank-$2$ point $z$ is said to be {\em exterior} if it lies on a tangent to the conic $\mathcal{C}_z$, and {\em interior} otherwise. 
The following lemma is useful for determining whether a rank-$2$ point is exterior or interior.
Here $M_{ij}(A)$ denotes the matrix obtained from a matrix $A$ by removing the $i$th row and the $j$th column, $|\cdot|$ denotes the determinant, and the matrix $M_y$ is defined as in \eqref{eq:symmetricMatrix}.

\begin{La}\label{lem:extcriteria}
A point $y\in \PG(5,q)$ is an exterior rank-$2$ point if and only if $|M_y|=0$ and 
$-|M_{11}(M_y)|$, $-|M_{22}(M_y)|$, $-|M_{33}(M_y)|$
are all squares with at least one being non-zero.
\end{La}

\begin{Prf}
See \cite[Lemma~3.7]{LaPoSh2020}.
\end{Prf}

According to Dickson's classification \cite{Dickson1908}, the lines in $\PG(5,q)$ form $15$ orbits under the action of $K$. 
Representatives are given in Table~\ref{table:lines}, as per our earlier papers \cite{LaPo2020,LaPoSh2020}. 
Note that for certain line orbits, namely $o_{i,j}$ with $i \in \{8,13,14\}$, we have used different representatives than those originally given in \cite[Table~2]{LaPo2020}; the reason for this is explained in \cite[Remarks 3.6 and 4.2]{LaPoSh2020}, and the representatives used here are the same as those used in \cite{LaPoSh2020}.
As shown in \cite{LaPoSh2020}, the planes of rank one, namely those intersecting $\mathcal{V}(\mathbb{F}_q)$ (and therefore containing a point of rank $1$ in the sense defined above), also form $15$ $K$-orbits; representatives are given in Table~\ref{table:main}. 

The {\em point-orbit distribution} of a subspace $W$ of $\PG(5,q)$ is the ordered list $[n_1,n_{2\text{e}},n_{2\text{i}},n_3]$, where $n_1$ is the number of rank $1$ points in $W$, $n_{2\text{e}}$ is the number of exterior rank-$2$ points in $W$, $n_{2\text{i}}$ is the number of interior rank-$2$ points in $W$, and $n_3$ is the number of rank $3$ points in $W$. 
The point-orbit distributions of lines of $\PG(5,q)$ are given in Table~\ref{table:line-rank-dist}. 
The point-orbit distributions of planes of rank one are given in Table~\ref{table:pt-orbit-dist}. 
We occasionally refer also to the {\em rank distribution} of the subspace $W$, which is the list $[n_1,n_{2\text{e}}+n_{2\text{i}},n_3]$. 

Note that the $K$-orbit of a line $\ell$ in $\PG(5,q)$ is determined by its point-orbit distribution, unless the point-orbit distribution is $[0,1,0,q]$, in which case $\ell$ can be either of type $o_{15,1}$ or of type $o_{16}$. 
(The rank distribution is a strictly coarser invariant: for each $i \in \{8,13,14,15\}$, the orbits $o_{i,1}$ and $o_{i,2}$ have equal rank distributions but non-equal point-orbit distributions.) 
The following lemma gives a geometric criterion (as opposed to a combinatorial one) for distinguishing between these two orbits. 
Here $S_{n,n}(\mathbb{F}_q)$ is the {\em Segre variety} in $\PG(n^2+2n,q)$, that is, the image of the map taking $(\langle x \rangle, \langle y \rangle) \in \PG(n-1,q) \times \PG(n-1,q)$ to $\langle x \otimes y \rangle$. 

\begin{La} \label{lemma_o15_o16}
Let $\ell$ be a line of type $o_{15,1}$ or $o_{16}$ in $\langle \mathcal{V}(\mathbb{F}_q) \rangle \cong \PG(5,q)$. 
Let $w$ be the unique point of rank $2$ on $\ell$, and let $W=W(w)$ be the unique solid in $\langle S_{3,3}(\mathbb{F}_q) \rangle \cong \PG(8,q)$ that contains $w$ and intersects the Segre variety $S_{3,3}(\mathbb{F}_q)$ in a subvariety $Q=Q(w)$ equivalent to a Segre variety $S_{2,2}(\mathbb{F}_q)$. 
Write $U = \langle W,\ell \rangle$. 
Then $\ell$ has type $o_{15,1}$ if there exists a point of rank $1$ in $U \setminus Q$, and $\ell$ has type $o_{16}$ otherwise.
\end{La}

\begin{Prf}
This is \cite[Lemma~1]{AlLa2019}.
\end{Prf}

We often use the next lemma when applying Lemma~\ref{lemma_o15_o16}.

\begin{La} \label{lemma_o15_o16_application}
Let $c_1,c_2,c_3$ be fixed elements of $\mathbb{F}_q$, not all zero, and consider the matrix
\[
M = \left[ \begin{matrix} 
d_{11} & d_{12} & c_1 \\ 
d_{12} & d_{22} & c_2 \\ 
c_1 & c_2 & c_3 \end{matrix} \right],
\]
where $d_{11}, d_{12}, d_{21}, d_{22}$ are variables in $\mathbb{F}_q$. 
If $c_3=0$ then $M$ has rank at least $2$ for all choices of the $d_{ij}$.
If $c_3 \neq 0$ then $M$ has rank $1$ if each $d_{ij} = c_ic_jc_3^{-1}$. 
\end{La}

\begin{Prf}
If $c_3 \neq 0$ and each $d_{ij} = c_ic_jc_3^{-1}$ then $M$ is row equivalent to the matrix
\[
\left[ \begin{matrix} 
c_1^2 & c_1c_2 & c_1c_3 \\ 
c_1c_2 & c_2^2 & c_2c_3 \\ 
c_1 & c_2 & c_3 \end{matrix} \right]. 
\]
The first and second rows of this matrix are (possibly zero) scalar multiples of the third row (which is non-zero by assumption), so $M$ has rank $1$ as claimed. 
This verifies the second assertion. 
To verify the first assertion, suppose that $c_3=0$. 
Then either $c_1 \neq 0$ or $c_2 \neq 0$, by assumption. 
If $c_1 \neq 0$ then the first and third rows of $M$ are linearly independent, so $M$ has rank at least $2$. 
Similarly, if $c_2 \neq 0$ then the second and third rows of $M$ are linearly independent.
\end{Prf}


\begin{table}[!t]
\begin{tabular}{lll}
\toprule
Orbit & Representative & Conditions \\
\midrule
$o_5$ & $\left[ \begin{matrix} \alpha & \cdot & \cdot \\ \cdot & \beta & \cdot \\ \cdot & \cdot & \cdot \end{matrix} \right]$ & \\
\addlinespace[2pt]
$o_6$ & $\left[ \begin{matrix} \alpha & \beta & \cdot \\ \beta & \cdot & \cdot \\ \cdot & \cdot & \cdot \end{matrix} \right]$ & \\
\addlinespace[2pt]
$o_{8,1}$ & $\left[ \begin{matrix} \alpha & \cdot & \cdot \\ \cdot & \beta & \cdot \\ \cdot & \cdot & -\beta \end{matrix} \right]$& \\
\addlinespace[2pt]
$o_{8,2}$ & $\left[ \begin{matrix} \alpha & \cdot & \cdot \\ \cdot & \beta & \cdot \\ \cdot & \cdot & -\varepsilon\beta \end{matrix} \right]$ 
& $\varepsilon \in \boxtimes$ \\
\addlinespace[2pt]
$o_9$ & $\left[ \begin{matrix} \alpha & \cdot & \beta \\ \cdot & \beta & \cdot \\ \beta & \cdot & \cdot \end{matrix} \right]$ & \\
\addlinespace[2pt]
$o_{10}$ & $\left[ \begin{matrix} v\alpha & \beta & \cdot \\ \beta & \alpha+u\beta & \cdot \\ \cdot & \cdot & \cdot \end{matrix} \right]$ 
& ($*$) \\
\addlinespace[2pt]
$o_{12}$ & $\left[ \begin{matrix} \cdot & \alpha & \cdot \\ \alpha & \cdot & \beta \\ \cdot & \beta & \cdot \end{matrix} \right]$ & \\
\addlinespace[2pt]
		\addlinespace[26pt]
		\addlinespace[2pt]
		&& \\
		\bottomrule
		\end{tabular}
		\begin{tabular}{lll}
		\toprule
		Orbit & Representative & Conditions \\
		\midrule
$o_{13,1}$ & $\left[ \begin{matrix} \cdot & \alpha & \cdot \\ \alpha & \beta & \cdot \\ \cdot & \cdot & -\beta \end{matrix} \right]$ & \\
\addlinespace[2pt]
$o_{13,2}$ & $\left[ \begin{matrix} \cdot & \alpha & \cdot \\ \alpha & \beta & \cdot \\ \cdot & \cdot & -\varepsilon\beta \end{matrix} \right]$ 
& $\varepsilon \in \boxtimes$ \\
\addlinespace[2pt]
$o_{14,1}$ & $\left[ \begin{matrix} \alpha & \cdot & \cdot \\ \cdot & -(\alpha+\beta) & \cdot \\ \cdot & \cdot & \beta \end{matrix} \right]$ & \\
\addlinespace[2pt]
$o_{14,2}$ & $\left[ \begin{matrix} \alpha & \cdot & \cdot \\ \cdot & -\varepsilon(\alpha+\beta) & \cdot \\ \cdot & \cdot & \beta \end{matrix} \right]$ 
& $\varepsilon \in \boxtimes$ \\
\addlinespace[2pt]
$o_{15,1}$ & $\left[ \begin{matrix} v\beta & \alpha & \cdot \\ \alpha & u\alpha+\beta & \cdot \\ \cdot & \cdot & \alpha \end{matrix} \right]$
& $-v \in \Box$, ($*$) \\
\addlinespace[2pt]
$o_{15,2}$ & $\left[ \begin{matrix} v\beta & \alpha & \cdot \\ \alpha & u\alpha+\beta & \cdot \\ \cdot & \cdot & \alpha \end{matrix} \right]$
& $-v \in \boxtimes$, ($*$) \\
\addlinespace[2pt]
$o_{16}$ & $\left[ \begin{matrix} \cdot & \cdot & \alpha \\ \cdot & \alpha & \beta \\ \alpha & \beta & \cdot \end{matrix} \right]$ & \\ 
\addlinespace[2pt]
$o_{17}$ & $\left[ \begin{matrix} v^{-1}\alpha & \beta & \cdot \\ \beta & u \beta - w \alpha & \alpha \\ \cdot & \alpha & \beta \end{matrix} \right]$ 
& ($**$) \\
\addlinespace[2pt]
\bottomrule
\end{tabular}
\caption{Representatives of the 15 line orbits in $\PG(5,q)$, $q$ odd, under the action of $K\cong\PGL(3,q) \le \PGL(6,q)$ defined in Section~\ref{s:prelim}. 
Here $\cdot$ denotes $0$, $(\alpha,\beta)$ ranges over all non-zero values in $\mathbb{F}_q^2$, and $\Box$ (respectively, $\boxtimes$) is the set of squares (respectively, non-squares) in $\mathbb{F}_q$. 
Condition~($*$): $v\lambda^2+uv\lambda - 1 \neq 0$, $\forall \lambda \in \mathbb{F}_q$. 
Condition~($**$): $\lambda^3+w \lambda^2- u \lambda+ v \neq 0$, $\forall \lambda \in \mathbb{F}_q$. 
}
\label{table:lines}
\end{table}

\begin{table}[!t]
\begin{tabular}{llc}
\toprule
Orbit & Point-orbit distribution & Stabiliser \\
\midrule
$o_5$ & $[2, \frac{q-1}{2}, \frac{q-1}{2}, 0]$ & $(E_q^2 : C_{q-1}^2) : C_2$ \\
\addlinespace[2pt]
$o_6$ & $[1, q, 0, 0]$ & $E_q^{1+2} : C_{q-1}^2$ \\
\addlinespace[2pt]
$o_{8,1}$ & $[1, 1, 0, q-1]$ & $\GO^+(2,q)$ \\
\addlinespace[2pt]
$o_{8,2}$ & $[1, 0, 1, q-1]$ & $\GO^-(2,q)$ \\
\addlinespace[2pt]
$o_9$ & $[1, 0, 0, q]$ & $E_q^2 : C_{q-1}$ \\
\addlinespace[2pt]
$o_{10}$ & $[0, \frac{q+1}{2}, \frac{q+1}{2}, 0]$ & $E_q^2 : \GO^-(2,q)$ \\
\addlinespace[2pt]
$o_{12}$ & $[0, q+1, 0, 0]$ & $\GL(2,q)$ \\
\addlinespace[2pt]
$o_{13,1}$ & $[0, 2, 0, q-1]$ & $C_{q-1} \times C_2$ \\
\addlinespace[2pt]
$o_{13,2}$ & $[0, 1, 1, q-1]$ & $C_{q-1} \times C_2$ \\
\addlinespace[2pt]
$o_{14,1}$ & $[0, 3, 0, q-2]$ & $C_2^2 : S_3$ \\
\addlinespace[2pt]
$o_{14,2}$ & $[0, 1, 2, q-2]$ & $C_2^2 : C_2$ \\
\addlinespace[2pt]
$o_{15,1}$ & $[0, 1, 0, q]$ & $C_2^2$ \\
\addlinespace[2pt]
$o_{15,2}$ & $[0, 0, 1, q]$ & $C_2^2$ \\
\addlinespace[2pt]
$o_{16}$ & $[0, 1, 0, q]$ & $E_q : C_{q-1}$ \\
\addlinespace[2pt]
$o_{17}$ & $[0, 0, 0, q+1]$ & $C_3$ \\
\bottomrule
\end{tabular}
\caption{Point-orbit distributions and stabilisers (under the action of $\PGL(3,q)$ defined in Section~\ref{s:prelim}) for line orbits in $\PG(5,q)$, $q$ odd. 
The group-theoretic notation used in the third column is explained in Section~\ref{s:stabilisers}. 
(Note that there is a typo in \cite[Table~3]{LaPo2020}, where the stabilisers were originally recorded.  
The stabiliser of a line of type $o_{10}$ is listed there as $E_q^2 : \textrm{O}^-(2,q)$, where $\textrm{O}^-(2,q)$ is understood to be the isometry group of a non-degenerate quadratic form of $-$ type on a $2$-dimensional vector space over $\mathbb{F}_q$, which has order $2(q+1)$. 
The $\textrm{O}^-(2,q)$ should be replaced by $\GO^-(2,q)$, namely the similarity group of such a form, which has order $2(q^2-1)$.
The corresponding orbit size is listed correctly in \cite[Table~4]{LaPo2020}.)
}
\label{table:line-rank-dist}
\end{table}

\begin{table}[!t]
\begin{tabular}{lll}
\toprule
Orbit & Representative & Conditions \\
\midrule
$\Sigma_1$ & $\left[ \begin{matrix} \alpha & \gamma & \cdot \\ \gamma & \beta & \cdot \\ \cdot & \cdot & \cdot \end{matrix} \right]$ & \\
\addlinespace[2pt]
$\Sigma_2$ & $\left[ \begin{matrix} \alpha & \cdot & \cdot \\ \cdot & \beta & \cdot \\ \cdot & \cdot & \gamma \end{matrix} \right]$ & \\
\addlinespace[2pt]
$\Sigma_3$ & $\left[ \begin{matrix} \alpha & \cdot & \gamma \\ \cdot & \beta & \cdot \\ \gamma & \cdot & \cdot \end{matrix} \right]$& \\
\addlinespace[2pt]
$\Sigma_4$ & $\left[ \begin{matrix} \alpha & \cdot & \gamma \\ \cdot & \beta & \gamma \\ \gamma & \gamma & \cdot \end{matrix} \right]$ & \\
\addlinespace[2pt]
$\Sigma_5$ & $\left[ \begin{matrix} \alpha & \cdot & \gamma \\ \cdot & \beta & \gamma \\ \gamma & \gamma & \gamma \end{matrix} \right]$ & \\
\addlinespace[2pt]
$\Sigma_6$ & $\left[ \begin{matrix} \alpha & \beta & \cdot \\ \beta & \varepsilon\alpha & \cdot \\ \cdot & \cdot & \gamma \end{matrix} \right]$ 
& $\varepsilon \in \boxtimes$ \\
\addlinespace[2pt]
$\Sigma_7$ & $\left[ \begin{matrix} \alpha & \beta & \gamma \\ \beta & \cdot & \cdot \\ \gamma & \cdot & \cdot \end{matrix} \right]$ & \\
\addlinespace[2pt]
$\Sigma_8$ & $\left[ \begin{matrix} \alpha & \beta & \cdot \\ \beta & \cdot & \gamma \\ \cdot & \gamma & \cdot \end{matrix} \right]$ & \\
\addlinespace[2pt]
\bottomrule
\end{tabular}
\begin{tabular}{lll}
\toprule
Orbit & Representative & Conditions \\
\midrule
$\Sigma_9$ & $\left[ \begin{matrix} \alpha & \beta & \cdot \\ \beta & \gamma & \cdot \\ \cdot & \cdot & -\gamma \end{matrix} \right]$ & \\
\addlinespace[2pt]
$\Sigma_{10}$ & $\left[ \begin{matrix} \alpha & \beta & \cdot \\ \beta & \gamma & \cdot \\ \cdot & \cdot & -\varepsilon\gamma \end{matrix} \right]$ & $\varepsilon \in \boxtimes$ \\
\addlinespace[2pt]
$\Sigma_{11}$ & $\left[ \begin{matrix} \cdot & \beta & \gamma \\ \beta & \alpha & \alpha \\ \gamma & \alpha & \alpha+\gamma \end{matrix} \right]$ & \\
\addlinespace[2pt]
$\Sigma_{12}$ & $\left[ \begin{matrix} \alpha & \beta & \cdot \\ \beta & \gamma & \beta \\ \cdot & \beta & \gamma \end{matrix} \right]$ & \\
\addlinespace[2pt]
$\Sigma_{13}$ & $\left[ \begin{matrix} \alpha & \beta & \cdot \\ \beta & \gamma & \beta \\ \cdot & \beta & \varepsilon\gamma \end{matrix} \right]$
& $\varepsilon \in \boxtimes$ \\
\addlinespace[2pt]
$\Sigma_{14}$ & $\left[ \begin{matrix} \alpha& \beta & \cdot \\ \beta & c\gamma & \beta-\gamma \\ \cdot & \beta-\gamma & \gamma \end{matrix} \right]$ 
& $q \not \equiv 0 \pmod 3$, ($\dagger$)
\\ 
\addlinespace[2pt]
$\Sigma_{14}'$ & $\left[ \begin{matrix} \alpha+\gamma & \gamma & \gamma \\ \gamma & \beta+\gamma & \gamma \\ \gamma & \gamma & -\beta \end{matrix} \right]$ 
& $q \equiv 0 \pmod 3$ \\
\addlinespace[2pt]
$\Sigma_{15}$ & $\left[ \begin{matrix} \alpha & \beta & \gamma \\ \beta & \gamma & \cdot \\ \gamma & \cdot & \cdot \end{matrix} \right]$ & \\
\addlinespace[2pt]
\bottomrule
\end{tabular}
\caption{Representatives of the 15 orbits of planes in $\PG(5,q)$, $q$ odd, meeting the quadric Veronesean in at least one point, under the action of $\PGL(3,q) \le \PGL(6,q)$ defined in Section~\ref{s:prelim}. 
Here $\cdot$ denotes $0$, $(\alpha,\beta,\gamma)$ ranges over all non-zero values in $\mathbb{F}_q^3$, and $\boxtimes$ is the set of non-squares in $\mathbb{F}_q$. 
In orbit $\Sigma_{14}$, condition~($\dagger$) is: $c \in \bF_q \setminus \{0,1\}$, $-3c$ is a square in $\mathbb{F}_q$, and $\tfrac{\sqrt{c}+1}{\sqrt{c}-1}$ is a non-cube in $\mathbb{F}_q(\sqrt{-3})$. 
}
\label{table:main}
\end{table}

\begin{table}[!t]
\begin{tabular}{lll}
\toprule
Orbit & Point-orbit distribution & Condition \\
\midrule
$\Sigma_1$ & $[q+1, \frac{q(q+1)}{2}, \frac{q(q-1)}{2}, 0]$ & \\
\addlinespace[2pt]
$\Sigma_2$ & $[3, \frac{3(q-1)}{2}, \frac{3(q-1)}{2}, q^2-2q+1]$ & \\
\addlinespace[2pt]
$\Sigma_3$ & $[2, \frac{3q-1}{2}, \frac{q-1}{2}, q^2-q]$ & \\
\addlinespace[2pt]
$\Sigma_4$ & $[2, \frac{3q-1}{2}, \frac{q-1}{2}, q^2-q]$ & \\
\addlinespace[2pt]
$\Sigma_5$ & $[2, q-1, q-1, q^2-q+1]$ & \\
\addlinespace[2pt]
$\Sigma_6$ & $[1, \frac{q+1}{2}, \frac{q+1}{2}, q^2-1]$ & \\
\addlinespace[2pt]
$\Sigma_7$ & $[1, q^2+q, 0, 0]$ & \\
\addlinespace[2pt]
$\Sigma_8$ & $[1, 2q, 0, q^2-q]$ & \\
\addlinespace[2pt]
$\Sigma_9$ & $[1, 2q, 0, q^2-q]$ & \\
\addlinespace[2pt]
$\Sigma_{10}$ & $[1, q, q, q^2-q]$ & \\
\addlinespace[2pt]
$\Sigma_{11}$ & $[1, q, 0, q^2]$ & \\
\addlinespace[2pt]
$\Sigma_{12}$ & $[1, \frac{q-1}{2}, \frac{q-1}{2}, q^2+1]$ & \\
\addlinespace[2pt]
$\Sigma_{13}$ & $[1, \frac{q+1}{2}, \frac{q+1}{2}, q^2-1]$ & \\
\addlinespace[2pt]
$\Sigma_{14}$ & $[1, \frac{q\mp 1}{2}, \frac{q\mp 1}{2}, q^2\pm 1]$ & $q \equiv \pm 1 \pmod 3$ \\
\addlinespace[2pt]
$\Sigma_{14}'$ & $[1, q, 0, q^2]$ & $q \equiv 0 \pmod 3$ \\
\addlinespace[2pt]
$\Sigma_{15}$ & $[1, q, 0, q^2]$ & \\
\bottomrule
\end{tabular}
\caption{Point-orbit distributions of planes in $\PG(5,q)$, $q$ odd, meeting the quadric Veronesean in at least one point.}
\label{table:pt-orbit-dist}
\end{table}

\section{Line-orbit distributions} \label{s:dist}

We now determine the line-orbit distributions of the planes in $\PG(5,q)$ that meet the quadric Veronesean in at least one point. 

\begin{La} \label{dist1}
A plane of type $\Sigma_1$ contains 
\begin{itemize}
\item $\tfrac{q(q+1)}{2}$ lines of type $o_5$, 
\item $q+1$ lines of type $o_6$, 
\item $\tfrac{q(q-1)}{2}$ lines of type $o_{10}$.
\end{itemize}
\end{La}

\begin{Prf}
The point-orbit distribution of a plane $\pi \in \Sigma_1$ is $[q+1, \frac{q(q+1)}{2}, \frac{q(q-1)}{2}, 0]$, and the points of rank $1$ lie on a conic $\mathcal{C}$. 
(For example, if $\pi$ is the representative of $\Sigma_1$ given in Table~\ref{table:main}, then the points of rank $1$ comprise the conic $\mathcal{C}:\alpha\beta-\gamma^2=0$.) 
In particular, $\pi$ is spanned by the points of $\mathcal{C}$; that is, $\pi$ is a conic plane. 
Each line $\ell$ in $\pi$ meets $\mathcal{C}$ in either two, one or zero points; that is, $\ell$ is a secant, tangent or external line to $\mathcal{C}$.
If $\ell$ is a secant to $\mathcal{C}$, then Table~\ref{table:line-rank-dist} implies that $\ell$ has type $o_5$, because these are the only types of lines containing two points of rank $1$. 
Hence, the number of lines of type $o_5$ in $\pi$ is equal to the number of secants to $\mathcal{C}$, namely $\binom{q+1}{2} = \tfrac{q(q+1)}{2}$. 
If $\ell$ is a tangent to $\mathcal{C}$, then $\ell$ has type $o_6$, because these are the only types of lines containing exactly one point of rank $1$ and no points of rank $3$.
Hence, the number of lines of type $o_6$ in $\pi$ is equal to the number of tangents to $\mathcal{C}$, which is $q+1$. 
The remaining $\tfrac{q(q-1)}{2}$ lines in $\pi$ are external to $\mathcal{C}$, so contain only points of rank $2$. 
Table~\ref{table:line-rank-dist} implies that all such lines have type $o_{10}$ or $o_{12}$.  
A simple counting argument then shows that they all have type $o_{10}$. 
Indeed, each of the $\tfrac{q(q-1)}{2}$ {\em interior} rank-$2$ points lies on $\tfrac{q+1}{2}$ secants and $\tfrac{q+1}{2}$ external lines to $\mathcal{C}$, and these external lines certainly have type $o_{10}$ because lines of type $o_{12}$ contain only exterior rank-$2$ points. 
On the other hand, each line of type $o_{10}$ contains $\tfrac{q+1}{2}$ interior rank-$2$ points, so a double count shows that $\pi$ contains a total of $\tfrac{q(q-1)}{2}$ lines of type $o_{10}$. 
That is, every external line has type $o_{10}$.
\end{Prf}

\begin{La} \label{dist2}
A plane of type $\Sigma_2$ contains 
\begin{itemize}
\item three lines of type $o_5$, 
\item $\tfrac{3(q-1)}{2}$ lines of type $o_{8,1}$, 
\item $\tfrac{3(q-1)}{2}$ lines of type $o_{8,2}$. 
\item $\tfrac{(q-1)^2}{4}$ lines of type $o_{14,1}$, 
\item $\tfrac{3(q-1)^2}{4}$ lines of type $o_{14,2}$.
\end{itemize}
\end{La}

\begin{Prf}
Let $\pi \in \Sigma_2$, and let $x,y,z$ denote the three points of rank $1$ in $\pi$. 
The three lines $\langle x,y \rangle$, $\langle y,z \rangle$, $\langle z,x \rangle$ have type $o_5$ and point-orbit distribution $[2, \frac{q-1}{2}, \frac{q-1}{2}, 0]$. 
Since $\pi$ has point-orbit distribution $[3, \frac{3(q-1)}{2}, \frac{3(q-1)}{2}, q^2-2q+1]$, any point not on one of these three lines has rank $3$. 
Now consider the other $q-1$ lines through $x$. 
Half of these lines meet $\langle y,z \rangle$ in an interior rank-$2$ point, and the other half meet $\langle y,z \rangle$ in an exterior rank-$2$ point. 
Hence, there are $\frac{q-1}{2}$ lines of each of types $o_{8,1}$ and $o_{8,2}$ through $x$. 
Repeating this argument for $y$ and $z$ yields a total of $\tfrac{3(q-1)}{2}$ lines of each of these two types. 

Now consider an exterior rank-$2$ point $u$. 
Without loss of generality, $u$ lies on $\langle y,z \rangle$ and on a line of type $o_{8,1}$ through $x$, both of which have already been counted. 
Let $\ell$ be one of the remaining $q-1$ lines through $u$. 
Then $\ell$ meets the lines $\langle x,y \rangle$ and $\langle x,z \rangle$ in points of rank $2$, and therefore contains exactly three points of rank $2$. 
It follows from Table~\ref{table:line-rank-dist} that either (i) the rank-$2$ points on $\ell$ are all exterior, in which case $\ell$ has type $o_{14,1}$; or (ii) one (namely $u$, in our case) is exterior and the other two are interior, in which case $\ell$ has type $o_{14,2}$. 
Since half of the $q-1$ points of rank $2$ on $\langle x,y \rangle$, say, are exterior and the other half are interior, it follows that $u$ lies on $\tfrac{q-1}{2}$ lines of each of the types $o_{14,1}$ and $o_{14,2}$. 
Since each line of type $o_{14,2}$ contains exactly one exterior rank-$2$ point, the total number of lines of type $o_{14,2}$ in $\pi$ is therefore equal to $\tfrac{q-1}{2}$ times the number of exterior rank-$2$ points in $\pi$, which is $\tfrac{3(q-1)}{2}$. 
That is, $\pi$ contains $\tfrac{3(q-1)^2}{4}$ lines of type $o_{14,2}$. 
On the other hand, each line of type $o_{14,1}$ contains three exterior rank-$2$ points, so the total number of lines of type $o_{14,1}$ in $\pi$ is $\tfrac{1}{3} \cdot \tfrac{3(q-1)}{2} \cdot \tfrac{q-1}{2}$. 
We have now accounted for all of the $q^2+q+1$ lines in $\pi$. 
\end{Prf}

\begin{La} \label{dist3}
A plane of type $\Sigma_3$ contains 
\begin{itemize}
\item one line of type $o_5$, 
\item one line of type $o_6$, 
\item $q$ lines of type $o_{8,1}$,  
\item $q-1$ lines of type $o_9$, 
\item $\tfrac{q(q-1)}{2}$ lines of type $o_{13,1}$, 
\item $\tfrac{q(q-1)}{2}$ lines of type $o_{13,2}$. 
\end{itemize}
\end{La}

\begin{Prf}
Let $\pi$ denote the representative of $\Sigma_3$ given in Table~\ref{table:main}, namely
\[
\left[ \begin{matrix}
\alpha & \cdot & \gamma \\
\cdot & \beta & \cdot \\
\gamma & \cdot & \cdot
\end{matrix} \right].
\]
As explained in the caption of the table, here $\cdot$ denotes $0$, and the triple of parameters $(\alpha,\beta,\gamma)$ ranges over $\mathbb{F}_q^3 \setminus \{(0,0,0)\}$.  
Recall that $\pi$ has point-orbit distribution $[2, \frac{3q-1}{2}, \frac{q-1}{2}, q^2-q]$. 
Let $x$ and $y$ denote the rank-$1$ points given by $\beta=\gamma=0$ and $\alpha=\gamma=0$, respectively. 
The line $\langle x,y \rangle$ has type $o_5$. 
The rank-$2$ points not on this line are all on the line $\ell : \beta = 0$, which has type $o_6$ (and contains $x$). 
In particular, every line through $x$ other than $\langle x,y \rangle$ or $\ell$ has point-orbit distribution $[1,0,0,q]$ and hence type $o_9$, yielding a total of $q-1$ such lines in $\pi$. 
Every line through $y$ other than $\langle x,y \rangle$ meets $\ell$ in an exterior point of rank $2$ and therefore has type $o_{8,1}$, giving a total of $q$ such lines. 
Now consider a point $w$ of rank $2$ on the line $\langle x,y \rangle$. 
Each of the $q$ other lines through $w$ meets $\ell$ in an exterior point of rank $2$, so has point-orbit distribution $[0,2,0,q-1]$ or $[0,1,1,q-1]$, and hence type $o_{13,1}$ or $o_{13,2}$, according to whether $w$ is exterior or interior. 
This gives a total of $\tfrac{q(q-1)}{2}$ lines of each of these two types. 
\end{Prf}

\begin{La} \label{Sigma4} \label{dist4}
A plane of type $\Sigma_4$ contains 
\begin{itemize}
\item one line of type $o_5$, 
\item $2q$ lines of type $o_{8,1}$, 
\item one line of type $o_{12}$,
\item $\tfrac{q(q-3)}{2}$ lines of type $o_{13,1}$, 
\item $\tfrac{q(q-1)}{2}$ lines of type $o_{13,2}$, 
\item $q-1$ lines of type $o_{16}$.
\end{itemize}
\end{La}

\begin{Prf}
Let $\pi$ denote the representative of $\Sigma_4$ given in Table~\ref{table:main}, namely
\[
\left[ \begin{matrix}
\alpha & \cdot & \gamma \\
\cdot & \beta & \gamma \\
\gamma & \gamma & \cdot
\end{matrix} \right].
\]
Recall that $\pi$ has point-orbit distribution $[2, \frac{3q-1}{2}, \frac{q-1}{2}, q^2-q]$. 
Let $x$ and $y$ denote the rank-$1$ points given by $\beta=\gamma=0$ and $\alpha=\gamma=0$, respectively. 
The line $\langle x,y \rangle$ has type $o_5$. 
The rank-$2$ points not on this line are all on the line $\ell : \beta = -\alpha$, which has type $o_{12}$, as can be verified by applying Lemma~\ref{lem:extcriteria}. 
In particular, every line through $x$ or $y$ other than $\langle x,y \rangle$ meets $\ell$ in an exterior point of rank $2$, so has point-orbit distribution $[1,1,0,q-1]$ and hence type $o_{8,1}$. 
This yields a total of $2q$ lines of type $o_{8,1}$.

The lines $\ell$ and $\langle x,y \rangle$ meet in the exterior rank-$2$ point $w:\alpha = -\beta = 1$, $\gamma = 0$.
Let $\ell'$ be any of the other $q-1$ lines through $w$. 
Then $\ell'$ has point-orbit distribution $[0,1,0,q]$ and hence type $o_{15,1}$ or $o_{16}$. 
We apply Lemma~\ref{lemma_o15_o16} to show that it has type $o_{16}$. 
Note that $\ell'$ is represented by a matrix of the form
\[
M = \left[ \begin{matrix}
a+b\alpha & \cdot & b\gamma \\
\cdot & -a+b\beta & b\gamma \\
b\gamma & b\gamma & \cdot
\end{matrix} \right], 
\]
where $\alpha,\beta,\gamma$ are fixed and $(a,b)$ ranges over $\mathbb{F}_q^2 \setminus \{(0,0)\}$. 
In the notation of Lemma~\ref{lemma_o15_o16}, the solid $W$ determined by $w$ is represented by the matrices whose third rows and third columns are zero. 
Hence, $U := \langle W,\ell' \rangle$ is represented by
\[
\left[ \begin{matrix}
d_{11} & d_{12} & b\gamma \\
d_{21} & d_{22} & b\gamma \\
b\gamma & b\gamma & \cdot
\end{matrix} \right], 
\]
where $b$ and the $d_{ij}$ range over $\mathbb{F}_q$. 
Since every point on $\ell'$ with $b \neq 0$ has rank $3$, we have $\gamma \neq 0$ and so the above matrix cannot have rank $1$ unless $b=0$. 
Hence, all rank-$1$ points in $U$ lie in $Q = W \cap S_{3,3}(\mathbb{F}_q)$, and the lemma therefore implies that $\ell'$ has type $o_{16}$. 

Finally, consider a rank-$2$ point $u \neq w$ on $\langle x,y \rangle$. 
Each of the other $q$ lines through $u$ meets $\ell$ in an exterior rank-$2$ point, so has point-orbit distribution $[0,2,0,q-1]$ or $[0,1,1,q-1]$, and hence type $o_{13,1}$ or $o_{13,2}$, according to whether $u$ is exterior or interior. 
There are $\tfrac{q-1}{2}$ choices for $u$ that are interior, and $\tfrac{q-1}{2}-1$ that are exterior (because $u \neq w$). 
Hence, $\pi$ contains in total $\tfrac{q(q-1)}{2}$ and $\tfrac{q(q-3)}{2}$ lines of type $o_{13,2}$ and $o_{13,1}$, respectively. 
\end{Prf}

\begin{La} \label{dist5}
A plane of type $\Sigma_5$ contains 
\begin{itemize}
\item one line of type $o_5$, 
\item $q-1$ lines of type $o_{8,1}$, 
\item $q-1$ lines of type $o_{8,2}$, 
\item two lines of type $o_9$, 
\item $\tfrac{q-1}{2}$ lines of type $o_{13,1}$, 
\item $\tfrac{q-1}{2}$ lines of type $o_{13,2}$, 
\item $\tfrac{(q-1)(q-3)}{8}$ lines of type $o_{14,1}$, 
\item $\tfrac{(q-1)(3q-5)}{8}$ lines of type $o_{14,2}$, 
\item $\tfrac{(q-1)^2}{4}$ lines of type $o_{15,1}$, 
\item $\tfrac{(q+1)(q-1)}{4}$ lines of type $o_{15,2}$.
\end{itemize}
\end{La}

\begin{Prf}
Let $\pi$ denote the representative of $\Sigma_5$ given in Table~\ref{table:main}, namely
\begin{equation} \label{sigma5rep}
\left[ \begin{matrix}
\alpha & \cdot & \gamma \\
\cdot & \beta & \gamma \\
\gamma & \gamma & \gamma
\end{matrix} \right].
\end{equation}
Recall that $\pi$ has point-orbit distribution $[2, q-1, q-1, q^2-q+1]$. 
Let $x$ and $y$ denote the rank-$1$ points in $\pi$ given by $\beta=\gamma=0$ and $\alpha=\gamma=0$. 
The line $\ell = \langle x,y \rangle$, given by $\gamma=0$, has type $o_5$. 
The $q-1$ rank-$2$ points in $\pi$ not on this line are all on the conic $\mathcal{C}:\alpha\beta-(\alpha+\beta)\gamma=0$; half are exterior and half are interior. 
The conic $\mathcal{C}$ meets the line $\ell$ in the points $x,y$. 
All but one of the other lines through $x$ meet $\mathcal{C}$ in a unique point of rank $2$; half of these $q-1$ lines have type $o_{8,1}$ and the other half have type $o_{8,2}$. 
The final line through $x$ is a tangent to $\mathcal{C}$ and so has type $o_9$. 
The same argument holds for the lines through $y$, so in total we obtain $q-1$ lines of type $o_{8,i}$ for each $i \in \{1,2\}$, and two lines of type $o_9$. 

Next, consider the tangent lines to $\mathcal{C}$ at its points of rank $2$. 
Let $t_u(\mathcal{C})$ denote the tangent line to $\mathcal{C}$ at a rank-$2$ point $u$. 
We claim that 
\begin{itemize}
\item[(i)] $t_u(\mathcal{C})$ meets $\ell$ in an exterior rank-$2$ point, regardless of whether $u$ is exterior or interior. 
\end{itemize}
Indeed, denote the coordinates of $u$ by $(\alpha_u,\beta_u,\gamma_u)$. 
(That is, suppose that $u$ is represented by the matrix in \eqref{sigma5rep} with a subscript `$u$' on each parameter.)
Then $\beta_u \neq \gamma_u$, because otherwise the equation of $\mathcal{C}$ would imply that $\gamma_u=0$, contradicting that assumption that $u$ does not lie on $\ell$. 
With this noted, a direct calculation shows that $t_u(\mathcal{C})$ meets $\ell$ in the point with coordinates $(-\gamma_u^2/(\beta_u-\gamma_u)^2,1,0)$, which is exterior by Lemma~\ref{lem:extcriteria}. 
The line $t_u(\mathcal{C})$ therefore has type $o_{13,1}$ or $o_{13,2}$ according to whether $u$ is exterior or interior, giving a total of $\tfrac{q-1}{2}$ lines of each type. 

Next, we verify that
\begin{itemize}
\item[(ii)] $\pi$ contains in total $\tfrac{(q-1)(q-3)}{8}$ lines of type $o_{14,1}$, and $\tfrac{(q-1)(3q-5)}{8}$ lines of type $o_{14,2}$. 
\end{itemize}

{\em Proof of} (ii). Every line of type $o_{14,1}$ passes through exactly two of the $N:=\tfrac{q-1}{2}$ exterior rank-$2$ points on $\mathcal{C}$, so there are $\binom{N}{2} = \tfrac{(q-1)(q-3)}{8}$ lines of type $o_{14,1}$ in total. 
Every line of type $o_{14,2}$ passes through either one or two of the $N$ interior rank-$2$ points on $\mathcal{C}$. 
There are $\binom{N}{2}$ lines of type $o_{14,2}$ containing two interior rank-$2$ points on $\mathcal{C}$, and $N^2$ lines of type $o_{14,2}$ containing one such point and one of the exterior rank-$2$ points on $\mathcal{C}$. 

To complete the proof of the lemma, we count the remaining lines through a rank-$2$ point $w$ on $\ell$. 
First suppose that $w$ is interior. 
By (i), $w$ lies on no tangent lines to $\mathcal{C}$, so every line through $w$ that meets $\mathcal{C}$ has type $o_{14,2}$, and has already been counted in (ii). 
There are $\tfrac{q-1}{2}$ such lines $w$. 
The remaining $(q+1)-1-\tfrac{q-1}{2} = \tfrac{q+1}{2}$ lines through $w$ do not meet $\mathcal{C}$, so have type $o_{15,2}$, giving a total of $\tfrac{(q+1)(q-1)}{4}$ such lines in $\pi$. 

Now suppose that $w$ is exterior. 
By (i), $w$ lies on the tangent lines to $\mathcal{C}$ at two rank-$2$ points. 
Any other line through $w$ that meets $\mathcal{C}$ has type $o_{14,1}$ or $o_{14,2}$. 
There are $\tfrac{q-3}{2}$ such lines through $w$, and we have already counted them in (ii). 
The remaining $(q+1)-3-\tfrac{q-3}{2}=\tfrac{q-1}{2}$ lines through $w$ have point-orbit distribution $[0,1,0,q]$ and therefore type $o_{15,1}$ or $o_{16}$. 
It remains to show that they all have type $o_{15,1}$, to give the claimed total of $\tfrac{(q-1)^2}{4}$ such lines in $\pi$. 
For this we apply Lemma~\ref{lemma_o15_o16}, as in the proof of Lemma~\ref{Sigma4}. 
The point $w$ is represented by the matrix in \eqref{sigma5rep} with $\gamma=0$, $\alpha=1$ and $\beta=-\varepsilon$ with $\varepsilon$ a non-zero square.
A line $\ell'$ through $w$ is therefore represented by 
\[
\left[ \begin{matrix}
a + b\alpha & \cdot & b\gamma \\
\cdot & -a\varepsilon + b\beta & b\gamma \\
b\gamma & b\gamma & b\gamma
\end{matrix} \right], 
\]
where $\alpha,\beta,\gamma$ are fixed and $(a,b)$ ranges over $\mathbb{F}_q^2 \setminus \{(0,0)\}$. 
In the notation of Lemma~\ref{lemma_o15_o16}, the solid $W$ determined by $w$ is represented by the matrices whose third rows and third columns are zero, so $U := \langle W,\ell' \rangle$ is represented by
\[
\left[ \begin{matrix}
d_{11} & d_{12} & b\gamma \\
d_{21} & d_{22} & b\gamma \\
b\gamma & b\gamma & b\gamma
\end{matrix} \right], 
\]
where $b$ and the $d_{ij}$ range over $\mathbb{F}_q$. 
Since every point on $\ell'$ with $b \neq 0$ has rank $3$, we have $\gamma \neq 0$. 
Hence, taking all of $b$ and the $d_{ij}$ to be the same yields a point of rank $1$ outside $Q = W \cap S_{3,3}(\mathbb{F}_q)$, and the lemma therefore implies that $\ell$ has type $o_{15,1}$. 
\end{Prf}

\begin{La} \label{dist6}
A plane of type $\Sigma_6$ contains 
\begin{itemize}
\item $\tfrac{q+1}{2}$ lines of type $o_{8,1}$, 
\item $\tfrac{q+1}{2}$ lines of type $o_{8,2}$, 
\item one line of type $o_{10}$, 
\item $\tfrac{(q+1)(q-1)}{2}$ lines of type $o_{15,1}$, 
\item $\tfrac{(q+1)(q-1)}{2}$ lines of type $o_{15,2}$. 
\end{itemize}
\end{La}

\begin{Prf}
Let $\pi$ denote the representative of $\Sigma_6$ from Table~\ref{table:main}, namely
\[
\left[ \begin{matrix} \alpha & \beta & \cdot \\ \beta & \varepsilon\alpha & \cdot \\ \cdot & \cdot & \gamma \end{matrix} \right], 
\quad \text{where } \varepsilon \text{ is a non-square in } \mathbb{F}_q.
\]
Then $\pi$ has point-orbit distribution $[1,\tfrac{q+1}{2},\tfrac{q+1}{2},q^2-1]$, and its points of rank~$2$ are all on the line $\gamma=0$, which has type $o_{10}$ (by Table~\ref{table:line-rank-dist}). 
Consider the unique point of rank~$1$ in $\pi$, namely $x : \alpha=\beta=0$. 
Each line through $x$ meets exactly one of the points of rank $2$, so $\pi$ contains $\tfrac{q+1}{2}$ lines in each of the orbits $o_{8,1}$ and $o_{8,2}$.

Now consider the lines through a rank-$2$ point $w$. 
We have already the line $\langle x,w \rangle$ (of type $o_{8,1}$ or $o_{8,2}$ according to whether $w$ is exterior or interior) and the line $\gamma=0$ (of type $o_{10}$). 
If $w$ is interior then the remaining $q-1$ lines through $u$ have point-orbit distribution $[0,0,1,q]$ and hence type $o_{15,2}$, giving a total of $\tfrac{q+1}{2} \cdot (q-1)$ such lines in $\pi$. 
If $w$ is exterior then the remaining $q-1$ lines through $w$ have point-orbit distribution $[0,1,0,q]$ and hence type $o_{15,1}$ or $o_{16}$. 
We use Lemma~\ref{lemma_o15_o16} to show that they all have type $o_{15,1}$, which gives the claimed total of $\tfrac{q+1}{2} \cdot (q-1)$ lines of type $o_{15,1}$ in $\pi$. 
Let $\ell$ be such a line. 
Since $w$ lies on the line $\gamma=0$, the solid $W$ determined by $w$ is represented by the matrices whose third rows and third columns are zero, so $U := \langle W,\ell \rangle$ is represented by
\[
\left[ \begin{matrix}
d_{11} & d_{12} & \cdot \\
d_{21} & d_{22} & \cdot \\
\cdot & \cdot & b\gamma
\end{matrix} \right], 
\]
for some fixed non-zero $\gamma$, where $b$ and the $d_{ij}$ range over $\mathbb{F}_q$. 
Taking all $d_{ij}=0$ yields a point of rank~$1$ outside $Q = W \cap S_{3,3}(\mathbb{F}_q)$, so $\ell$ has type $o_{15,1}$ as claimed.
\end{Prf}

\begin{La} \label{dist7}
A plane of type $\Sigma_7$ contains 
\begin{itemize}
\item $q+1$ lines of type $o_6$, 
\item $q^2$ lines of type $o_{12}$.
\end{itemize}
\end{La}

\begin{Prf}
A plane $\pi \in \Sigma_7$ has a unique point $x$ of rank $1$, and $q^2+q$ exterior points of rank $2$. 
The $q+1$ lines through $x$ in $\pi$ therefore have point-orbit distribution $[1,q,0,0]$ and type $o_6$. 
All other lines in $\pi$ have point-orbit distribution $[0,q+1,0,0]$, and are therefore of type $o_{12}$.
\end{Prf}

\begin{La} \label{dist8}
A plane of type $\Sigma_8$ contains 
\begin{itemize}
\item one line of type $o_6$, 
\item $q$ lines of type $o_{8,1}$, 
\item one line of type $o_{12}$, 
\item $q(q-1)$ lines of type $o_{13,1}$, 
\item $q-1$ lines of type $o_{16}$.
\end{itemize}
\end{La}

\begin{Prf}
Let $\pi$ denote the representative of $\Sigma_8$ given in Table~\ref{table:main}, namely
\[
\left[ \begin{matrix}
\alpha & \beta & \cdot \\
\beta & \cdot & \gamma \\
\cdot & \gamma & \cdot
\end{matrix} \right].
\]
The $2q$ points of rank $2$ in $\pi$ are all exterior and lie on the union of the line $\gamma=0$, which has type $o_6$, and the line $\alpha=0$, which has type $o_{12}$. 
Apart from the line $\gamma=0$, every other line through the unique point $x : \beta=\gamma=0$ of rank $1$ has point-orbit distribution $[1,1,0,q-1]$ and hence type $o_{8,1}$, giving a total of $q$ such lines. 

Now consider the point of intersection of the lines $\alpha=0$ and $\gamma=0$, call it $w$. 
The remaining $q-1$ lines through $w$ have point-orbit distribution $[0,1,0,q]$, and hence type $o_{15,1}$ or $o_{16}$. 
We use Lemma~\ref{lemma_o15_o16} to show that they all have type $o_{16}$. 
Let $\ell$ be such a line. 
The solid $W$ determined by $w$ is represented by the matrices whose third rows and third columns are zero, so $U := \langle W,\ell \rangle$ is represented by
\[
\left[ \begin{matrix}
d_{11} & d_{12} & \cdot \\
d_{21} & d_{22} & b\gamma \\
\cdot & b\gamma & \cdot
\end{matrix} \right], 
\]
for some fixed non-zero $\gamma$, where $b$ and the $d_{ij}$ range over $\mathbb{F}_q$. 
Such a matrix cannot have rank $1$ unless $b=0$, so all rank-$1$ points in $U$ lie in $Q = W \cap S_{3,3}(\mathbb{F}_q)$, and $\ell$ therefore has type $o_{16}$. 

Finally, let $\ell'$ be any of the $q(q-1)$ lines not considered thus far. 
Then $\ell'$ contains neither $x$ nor $w$, and intersects each of the lines $\alpha=0$ and $\gamma=0$ in an exterior point of rank $2$. 
The point-orbit distribution of $\ell'$ is therefore $[0,2,0,q-1]$, so $\ell'$ has type $o_{13,1}$. 
\end{Prf}

\begin{La} \label{dist9}
A plane of type $\Sigma_9$ contains 
\begin{itemize}
\item one line of type $o_6$, 
\item $q$ lines of type $o_{8,1}$, 
\item $q$ lines of type $o_{13,1}$, 
\item $\tfrac{q(q-1)}{2}$ lines of type $o_{14,1}$, 
\item $\tfrac{q(q-1)}{2}$ lines of type $o_{15,1}$. 
\end{itemize}
\end{La}

\begin{Prf}
Let $\pi$ denote the representative of $\Sigma_9$ given in Table~\ref{table:main}, namely
\[
\left[ \begin{matrix}
\alpha & \beta & \cdot \\
\beta & \gamma & \cdot \\
\cdot & \cdot & -\gamma
\end{matrix} \right].
\]
The $2q$ points of rank $2$ in $\pi$ are all exterior and lie on the union of the line $\ell : \gamma=0$ of type $o_6$ and the conic $\mathcal{C} : \beta^2-\alpha\gamma=0$, which meet (only) in the unique point $x : \beta=\gamma=0$ of rank $1$. 
Apart from $\ell$, every other line through $x$ meets $\mathcal{C}$ in a unique point of rank $2$, so has point-orbit distribution $[1,1,0,q-1]$ and hence type $o_{8,1}$, giving a total of $q$ such lines. 

If $u$ is one of the $q$ points of rank $2$ on $\mathcal{C}$ then $q-1$ of the lines through $u$ contain exactly one point on $\mathcal{C}$ other than $x$, and exactly one point on $\ell$ other than $x$. 
Such lines have point-orbit distribution $[0,3,0,q-2]$ and hence type $o_{14,1}$. 
Each contains two rank-$2$ points on $\mathcal{C}$, so there are in total $\tfrac{1}{2}\cdot q \cdot(q-1)$ such lines in $\pi$. 
The final line through $u$, other than $\langle x,u \rangle$, which we have already counted, is tangent to $\mathcal{C}$ and meets $\ell$ in a point of rank $2$, so has point-orbit distribution $[0,2,0,q-1]$ and hence type $o_{13,1}$. 
There are $q$ such lines in $\pi$.

Finally, consider a rank-$2$ point $w$ on $\ell$. 
Thus far, we have counted $\tfrac{q-1}{2}+2$ of the lines through $w$, namely $\ell$ itself, one line of type $o_{13,1}$, and $\tfrac{q-1}{2}$ lines of type $o_{14,1}$.
The remaining $\tfrac{q-1}{2}$ lines through $w$ have point-orbit distribution $[0,1,0,q]$, and hence type $o_{15,1}$ or $o_{16}$. 
Let $\ell'$ be such a line. 
To complete the proof, we use Lemma~\ref{lemma_o15_o16} to show $\ell'$ has type $o_{15,1}$, thereby yielding the claimed total of $\tfrac{q(q-1)}{2}$ such lines in $\pi$. 
Since $w$ lies on the line $\gamma=0$, the solid $W$ determined by $w$ is represented by the matrices whose third rows and third columns are zero, so $U := \langle W,\ell' \rangle$ is represented by
\[
\left[ \begin{matrix}
d_{11} & d_{12} & \cdot \\
d_{21} & d_{22} & \cdot \\
\cdot & \cdot & b\gamma
\end{matrix} \right], 
\]
for some fixed non-zero $\gamma$, where $b$ and the $d_{ij}$ range over $\mathbb{F}_q$. 
Taking all $d_{ij}=0$ yields a point of rank~$1$ outside $Q = W \cap S_{3,3}(\mathbb{F}_q)$, so $\ell'$ has type $o_{15,1}$ as claimed.
\end{Prf}

\begin{La} \label{dist10}
A plane of type $\Sigma_{10}$ contains 
\begin{itemize}
\item one line of type $o_6$, 
\item $q$ lines of type $o_{8,2}$, 
\item $q$ lines of type $o_{13,2}$, 
\item $\tfrac{q(q-1)}{2}$ lines of type $o_{14,2}$, 
\item $\tfrac{q(q-1)}{2}$ lines of type $o_{15,1}$.
\end{itemize}
\end{La}

\begin{Prf}
The proof is analogous to that of Lemma~\ref{dist9}, but we include some details for clarity. 
Let $\pi$ denote the representative of $\Sigma_{10}$ given in Table~\ref{table:main}, namely
\[
\left[ \begin{matrix}
\alpha & \beta & \cdot \\
\beta & \gamma & \cdot \\
\cdot & \cdot & -\varepsilon\gamma
\end{matrix} \right], 
\quad \text{where } \varepsilon \text{ is a non-square in } \mathbb{F}_q.
\]
The point-orbit distribution of $\pi$ is $[1,q,q,q^2-q]$. 
The points of rank $2$ in $\pi$ lie on the union of the line $\ell : \gamma=0$ of type $o_6$ and the conic $\mathcal{C} : \beta^2-\alpha\gamma=0$, which meet (only) in the unique point $x : \beta=\gamma=0$ of rank $1$. 
The rank-$2$ points on $\mathcal{C}$ are all interior, and those on $\ell$ are all exterior. 
Every line through $x$ other than $\ell$ therefore has point-orbit distribution $[1,0,1,q-1]$ and type $o_{8,2}$, giving a total of $q$ such lines. 

If $u$ is a (necessarily interior) rank-$2$ point on $\mathcal{C}$ then $q-1$ of the lines through $u$ contain exactly one point on $\mathcal{C}$ other than $x$, and exactly one point on $\ell$ other than $x$. 
Such lines have point-orbit distribution $[0,1,2,q-2]$ and hence type $o_{14,2}$. 
Each contains two choices of $u$, so $\pi$ contains $\tfrac{q(q-1)}{2}$ lines of type $o_{14,2}$ in total. 
The final line through $u$, other than $\langle x,u \rangle$, is tangent to $\mathcal{C}$ and meets $\ell$ in a (necessarily exterior) point of rank $2$, so has point-orbit distribution $[0,1,1,q-1]$ and type $o_{13,2}$. 
There are $q$ such lines in total.

Finally, consider a rank-$2$ point $w$ on $\ell$. 
We have counted $\tfrac{q-1}{2}+2$ of the lines through $w$, namely $\ell$ itself, one line of type $o_{13,2}$, and $\tfrac{q-1}{2}$ lines of type $o_{14,2}$.
The remaining $\tfrac{q-1}{2}$ lines through $w$ have point-orbit distribution $[0,1,0,q]$, and hence type $o_{15,1}$ or $o_{16}$. 
The same argument as in Lemma~\ref{dist9} shows that they all have type $o_{15,1}$, giving a total of $\tfrac{q(q-1)}{2}$ such lines in $\pi$.
\end{Prf}

It remains to determine the line-orbit distributions of the planes in the $K$-orbits $\Sigma_{11}$, $\Sigma_{12}$, $\Sigma_{13}$, $\Sigma_{14}$, $\Sigma_{14}'$ and $\Sigma_{15}$. 
We treat these remaining cases (roughly) in order of difficulty: 
\begin{itemize}
\item $\Sigma_{15}$ is treated first, as it is by far the most straightforward remaining case.
\item We then handle the case $\Sigma_{11}$, which is considerably more complicated. 
\item The case $\Sigma_{14}'$ is treated next as it is largely analogous to $\Sigma_{11}$, having the same point-orbit distribution, and, in particular, rank-$2$ points of only one type (exterior).
\item Finally, we deal with $\Sigma_{12}$, $\Sigma_{13}$ and $\Sigma_{14}$. These cases pose further complications, owing to the presence of both exterior and interior points of rank $2$. 
Amongst them, $\Sigma_{14}$ is arguably the most straightforward, at least in the sense that $q \not \equiv 0 \pmod 3$ and the cubic of points of rank $\le 2$ has no inflexion points, so we opt treat it first. 
We then treat $\Sigma_{13}$ before $\Sigma_{12}$ because the representative of $\Sigma_{12}$ in Table~\ref{table:lines} may be obtained from the representative of $\Sigma_{13}$ by setting $\varepsilon=1$, so various preliminary calculations in the argument for $\Sigma_{13}$ follow through immediately for $\Sigma_{12}$. 
\end{itemize}

\begin{La} \label{dist15}
A plane of type $\Sigma_{15}$ contains 
\begin{itemize}
\item one line of type $o_6$, 
\item $q$ lines of type $o_9$, 
\item $q^2$ lines of type $o_{16}$.
\end{itemize}
\end{La}

\begin{Prf}
Let $\pi$ denote the representative of $\Sigma_{15}$ given in Table~\ref{table:main}, namely
\[
\left[ \begin{matrix}
\alpha & \beta & \gamma \\
\beta & \gamma & \cdot \\
\gamma & \cdot & \cdot
\end{matrix} \right].
\]
The point-orbit distribution of $\pi$ is $[1,q,0,q^2]$. 
The unique point of rank $1$ is $x : \beta=\gamma=0$, and all of the rank-$2$ points lie on the line $\gamma=0$, which contains $x$ and has type $o_6$. 
The remaining $q$ lines through $x$ have point-orbit distribution $[1,0,0,q]$ and hence type $o_9$.
All of the other lines in $\pi$ have point-orbit distribution $[0,1,0,q]$, and hence are of type $o_{15,1}$ or $o_{16}$. 
We use Lemma~\ref{lemma_o15_o16} to show that they are all of type $o_{16}$. 
Let $\ell$ be such a line. 
The unique point $w$ of rank $2$ on $\ell$ lies on the line $\gamma=0$, so the solid $W$ determined by $w$ is represented by the matrices whose third rows and third columns are zero. 
Hence, $U := \langle W,\ell \rangle$ is represented by
\[
\left[ \begin{matrix}
d_{11} & d_{12} & b \\
d_{21} & d_{22} & \cdot \\
b & \cdot & \cdot
\end{matrix} \right], 
\]
where $b$ and the $d_{ij}$ range over $\mathbb{F}_q$. 
Such a matrix cannot have rank $1$ unless $b=0$, so all rank-$1$ points in $U$ lie in $Q = W \cap S_{3,3}(\mathbb{F}_q)$, and $\ell$ therefore has type $o_{16}$. 
\end{Prf}

\begin{La} \label{dist11}
The lines in a plane of type $\Sigma_{11}$ are as indicated in Table~\ref{sigma11table}. 
\begin{table}[!h]
\begin{tabular}{lcc}
\toprule
Type & $q \not \equiv 0 \pmod 3$ & $q \equiv 0 \pmod 3$ \\
\midrule
$o_{8,1}$ & $q$ & $q$ \\
$o_9$ & $1$ & $1$ \\
$o_{13,1}$ & $q-1$ & $q$ \\
$o_{14,1}$ & $\tfrac{(q-1)(q-2)}{6}$ & $\tfrac{q(q-3)}{6}$ \\
$o_{15,1}$ & $\tfrac{q(q-1)}{2}$ & $\tfrac{q(q-1)}{2}$ \\
$o_{16}$ & $1$ & $0$ \\
$o_{17}$ & $\tfrac{(q+1)(q-1)}{3}$ & $\tfrac{q^2}{3}$ \\
\bottomrule
\end{tabular}
\caption{Data for Lemma~\ref{dist11}.}
\label{sigma11table}
\end{table}
\end{La}

\begin{Prf}
If $M$ denotes the matrix representative of $\Sigma_{11}$ given in Table~\ref{table:main}, then we have 
\begin{equation} \label{Sigma11matrix}
XMX^T = 
\left[ \begin{matrix}
\alpha & \beta & \cdot \\
\beta & \cdot & \gamma-\beta \\
\cdot & \gamma-\beta & \gamma
\end{matrix} \right], 
\quad \text{where} \quad 
X = \left[ \begin{matrix}
\cdot & 1 & \cdot \\
1 & \cdot & \cdot \\
\cdot & -1 & 1
\end{matrix} \right] \in \GL(3,q).
\end{equation}
For convenience, we work with the above representative of $\Sigma_{11}$ instead of the one in Table~\ref{table:main}. 
Let us call the corresponding plane $\pi$. 
The points of rank at most $2$ in $\pi$ lie on the cubic $\mathcal{C}:\beta^2\gamma+\alpha(\beta-\gamma)^2=0$. 
There is a unique point of rank $1$, namely $x:\beta=\gamma=0$, and $q$ points of rank $2$, one for each choice of $(\beta,\gamma)$ with $\beta\neq\gamma$. 
All points of rank $2$ are exterior. 
In particular, $q$ of the lines through $x$ have point-orbit distribution $[1,1,0,q-1]$ and (therefore) type $o_{8,1}$, and the final line through $x$ has point-orbit distribution $[1,0,0,q]$ and type $o_9$. 

We now count the remaining lines in $\pi$ containing points of rank $2$. 
Note that if $(\alpha,\beta,\gamma)$ are the coordinates of a point of rank $2$, then $\beta \neq \gamma$ and $\alpha = -\beta^2\gamma/(\beta-\gamma)^2$. 
Note also that $\mathcal{C}$ has one point of inflexion if $q \not \equiv 0 \pmod 3$, and no points of inflexion if $q \equiv 0 \pmod 3$. 
In the former case, the unique point of inflexion is given by $2\beta+\gamma=0$. 
We need the following fact:
\begin{itemize}
\item[($*$)] With the exception of the inflexion point in characteristic $\neq 3$, every rank-$2$ point $v$ on $\mathcal{C}$ lies on exactly two tangent lines to $\mathcal{C}$ (one of which is the tangent line at $v$). 
\end{itemize}

{\em Proof of} ($*$). 
Consider a point $w$ of rank $2$ on $\mathcal{C}$, denoting its coordinates by $(\alpha_w,\beta_w,\gamma_w)$, where $\alpha_w = -\beta_w^2\gamma_w/(\beta_w-\gamma_w)^2$. 
In the `free' coordinates $(\alpha,\beta,\gamma)$, the tangent line to $\mathcal{C}$ at $w$ is
\begin{equation} \label{Sigma11tangent}
\alpha\cdot (\beta_w-\gamma_w)^3 - \beta \cdot 2\beta_w\gamma_w^2 + \gamma\cdot\beta_w^2(\beta_w+\gamma_w) = 0.
\end{equation}
Putting the coordinates of $v$, namely $(\alpha,\beta,\gamma)=(\alpha_v,\beta_v,\gamma_v) = (-\beta_v^2\gamma_v/(\beta_v-\gamma_v)^2,\beta_v,\gamma_v)$, into this equation yields
\[
\frac{(\beta_v\gamma_w - \beta_w\gamma_v)^2}{(\beta_v-\gamma_v)^2} \cdot 
(\gamma_v(\beta_w+\gamma_w)-2\beta_v\beta_w) = 0,
\]
so $v$ lies on the tangent line to $\mathcal{C}$ at $w$ if and only if either $\beta_v\gamma_w-\beta_w\gamma_v=0$, namely $w=v$, or  
\[
\gamma_v(\beta_w+\gamma_w)-2\beta_v\beta_w = 0.
\]
Regarding $v$ as fixed, the latter equation determines a unique rank-$2$ point $w$, which is distinct from $v$ unless $2\beta_v + \gamma_v = 0$, that is, unless $v$ is the unique point of inflexion of $\mathcal{C}$.

First suppose that $q \not \equiv 0 \pmod 3$, and let $z$ denote the unique point of inflexion of $\mathcal{C}$. 
Consider the lines through $z$, recalling that we have already counted the line $\langle z,x \rangle$. 
The tangent line to $\mathcal{C}$ through $z$ contains no other points of $\mathcal{C}$, so has point-orbit distribution $[0,1,0,q]$ and hence type $o_{15,1}$ or $o_{16}$. 
We show below that it has type $o_{16}$: see Claim~$\Sigma_{11}$ at the end of the proof. 
The remaining $q-1$ lines through $z$ contain either zero or two other points on $\mathcal{C}$. 
The $\tfrac{q-1}{2}$ lines containing two other points of $\mathcal{C}$ have point-orbit distribution $[0,3,0,q-2]$ and hence type $o_{14,1}$. 
The $\tfrac{q-1}{2}$ lines containing no other points of $\mathcal{C}$ have point-orbit distribution $[0,1,0,q]$; we show below (Claim~$\Sigma_{11}$) that they all have type $o_{15,1}$. 
Now consider a rank-$2$ point $w \neq z$. 
We have already counted the line $\langle w,x \rangle$ (of type $o_{8,1}$) and the line $\langle w,z \rangle$ (of type $o_{14,1}$). 
The tangent line to $\mathcal{C}$ through $w$ meets exactly one other (rank-$2$) point $w'$ on $\mathcal{C}$, so has point-orbit distribution $[0,2,0,q-1]$ and hence type $o_{13,1}$. 
By ($*$), $w$ also lies on the tangent line to $\mathcal{C}$ at a unique rank-$2$ point $w'' \neq w'$, which gives a second line $\langle w,w'' \rangle$ through $w$ of type $o_{13,1}$. 
There are $q-1$ choices for $w$, so we obtain a total of $q-1$ lines of type $o_{13,1}$ in $\pi$. 
At this point we have a further $q-3$ lines through $w$ to consider. 
There are $q-5$ points of rank~$2$ on $\mathcal{C}$ that do {\em not} lie on any of the lines $\langle w,z \rangle$, $\langle w,w' \rangle$, $\langle w,w'' \rangle$. 
Hence, $\tfrac{q-5}{2}$ of the remaining $q-3$ lines through $w$ have point-orbit distribution $[0,3,0,q-2]$ and type $o_{14,1}$, giving a further $\tfrac{1}{3} \cdot (q-1) \cdot \frac{q-5}{2}$ lines of type $o_{14,1}$ in addition to the $\tfrac{q-1}{2}$ lines of type $o_{14,1}$ through $z$ counted above. 
Therefore, $\pi$ contains in total $\frac{(q-1)(q-2)}{6}$ lines of type $o_{14,1}$, as claimed. 
The final $(q-3) - \tfrac{q-5}{2} = \tfrac{q-1}{2}$ lines through $w$ have point-orbit distribution $[0,1,0,q]$ and hence type $o_{15,1}$ or $o_{16}$. 
We show below that they all have type $o_{15,1}$. 
This yields $\tfrac{(q-1)^2}{2}$ lines of type $o_{15,1}$ in addition to the $\tfrac{q-1}{2}$ lines of type $o_{15,1}$ through $z$, giving the claimed total of $\tfrac{q(q-1)}{2}$ lines of type $o_{15,1}$ in $\pi$. 
All of the lines in $\pi$ not counted thus far contain only points of rank $3$, and so have type $o_{17}$.

Now suppose that $q \equiv 0 \pmod 3$. 
The types of the lines through the unique point $x$ of rank~$1$ do not change. 
However, now $\mathcal{C}$ has no points of inflexion, so {\em every} point $w$ of rank $2$ lies on two lines of type $o_{13,1}$ (namely, the tangent line to $\mathcal{C}$ at $w$ and the tangent line to $\mathcal{C}$ at some other point $w''$ of rank $2$, in the notation of the above argument). 
Therefore, $\pi$ now contains $q$ rather than $q-1$ lines of type $o_{13,1}$. 
The remaining $q-2$ lines through each $w$ now meet a total of $q-3$ remaining rank-$2$ points on $\mathcal{C}$. 
This gives a total of $\tfrac{1}{3} \cdot q \cdot \frac{q-3}{2}$ lines of type $o_{14,1}$, and $q\cdot ((q-2)-\tfrac{q-3}{2}) = \tfrac{q(q-1)}{2}$ lines with point-orbit distribution $[0,1,0,q-1]$. 
By Claim~$\Sigma_{11}$ below, all of the latter lines have type $o_{15,1}$. 
All of the other lines in $\pi$ have type $o_{17}$.

{\em Claim~$\Sigma_{11}$.} Let $\ell$ be a line in $\pi$ with point-orbit distribution $[0,1,0,q]$, and let $w$ denote its unique point of rank $2$. 
Then $\ell$ has type $o_{15,1}$ unless $q \not \equiv 0 \pmod 3$, $w$ is the unique point of inflexion of $\mathcal{C}$, and $\ell$ is the tangent line to $\mathcal{C}$ at $w$, in which case $\ell$ has type $o_{16}$.

{\em Proof of Claim~$\Sigma_{11}$.} 
Let $(\alpha_w,\beta_w,\gamma_w)$ denote the coordinates of $w$. 
Recall that $\alpha_w = -\beta_w^2\gamma_w/(\beta_w-\gamma_w)^2$, and that the tangent line to $\mathcal{C}$ at $w$ is given by \eqref{Sigma11tangent}. 
We apply Lemma~\ref{lemma_o15_o16}. 
In order to calculate the solid $W$ determined by $w$, we first use an appropriate element of $K$ to transform the matrix representative of $w$ to one whose third row and third column are zero. 
This allows us to apply Lemma~\ref{lemma_o15_o16_application} in a uniform way to determine the type of $\ell$. 
First suppose that $\beta_w\gamma_w \neq 0$. 
Then the matrix representative of $w$ obtained from \eqref{Sigma11matrix}, call it $M_w$, can be transformed via the action of $K$ as follows:
\[
Y_wM_wY_w^T = \left[ \begin{matrix}
-\frac{\beta_w^2\gamma_w}{(\beta_w-\gamma_w)^2} & \beta_w & \cdot \\
\beta_w & \cdot & \cdot \\
\cdot & \cdot & \cdot
\end{matrix} \right], 
\quad \text{where} \quad 
Y_w = \left[ \begin{matrix}
1 & \cdot & \cdot \\
\cdot & 1 & \cdot \\
\frac{(\beta_w-\gamma_w)^2}{\beta_w\gamma_w} & \cdot & \frac{\beta_w-\gamma_w}{\gamma_w}
\end{matrix} \right] \in \GL(3,q).
\]
After making this transformation, the solid $W$ of Lemma~\ref{lemma_o15_o16} is, as desired, represented by the matrices whose third rows and third columns are zero. 
Now choose a point $y$ (of rank $3$) such that $\ell = \langle w,y \rangle$, denoting the coordinates of $y$ by $(\alpha,\beta,\gamma)$. 
Then the image of $y$ under the element of $K$ corresponding to the matrix $Y_w$ is represented by the matrix 
\[
\left[ \begin{matrix}
\alpha & \beta & \frac{\alpha(\beta_w-\gamma_w)^2+\beta\beta_w\gamma_w}{\beta_w\gamma_w} \\
\beta & \cdot & \frac{(\gamma\beta_w-\beta\gamma_w)(\beta_w-\gamma_w)}{\beta_w\gamma_w} \\
\frac{\alpha(\beta_w-\gamma_w)^2+\beta\beta_w\gamma_w}{\beta_w\gamma_w} & \frac{(\gamma\beta_w-\beta\gamma_w)(\beta_w-\gamma_w)}{\beta_w\gamma_w} & \frac{\beta_w-\gamma_w}{\beta_w\gamma_w}f_w(\alpha,\beta,\gamma)
\end{matrix} \right],
\]
where $f_w(\alpha,\beta,\gamma)$ is the left-hand side of \eqref{Sigma11tangent}.
Hence, in the notation of Lemma~\ref{lemma_o15_o16}, the image of $U := \langle W,\ell \rangle$ is represented by
\[
\left[ \begin{matrix}
d_{11} & d_{12} & b\cdot\frac{\alpha(\beta_w-\gamma_w)^2+\beta\beta_w\gamma_w}{\beta_w\gamma_w} \\
d_{21} & d_{22} & b\cdot\frac{(\gamma\beta_w-\beta\gamma_w)(\beta_w-\gamma_w)}{\beta_w\gamma_w} \\
b\cdot\frac{\alpha(\beta_w-\gamma_w)^2+\beta\beta_w\gamma_w}{\beta_w\gamma_w} & b\cdot\frac{(\gamma\beta_w-\beta\gamma_w)(\beta_w-\gamma_w)}{\beta_w\gamma_w} & b\cdot\frac{\beta_w-\gamma_w}{\beta_w\gamma_w}f_w(\alpha,\beta,\gamma)
\end{matrix} \right],
\]
where $b$ and the $d_{ij}$ range over $\mathbb{F}_q$. 
Lemma~\ref{lemma_o15_o16} now tells us that $\ell$ has type $o_{15,1}$ if and only if there exist $d_{ij}$ such that the above matrix has rank $1$ when $b \neq 0$. 
Lemma~\ref{lemma_o15_o16_application} therefore implies that $\ell$ has type $o_{15,1}$ unless $f_w(\alpha,\beta,\gamma)=0$, that is, unless $y$ lies on the tangent line to $\mathcal{C}$ at $w$. 
(Note that the condition that the $c_i$ not be all zero holds because the point $y$ is assumed to have rank $3$.)
However, $\ell$ cannot be the tangent line to $w$ unless $q \not \equiv 0 \pmod 3$ and $w$ is the unique point of inflexion of $\mathcal{C}$, because (as explained above) in all other cases the tangent line to $\mathcal{C}$ at $w$ has point-orbit distribution $[0,2,0,q-1]$, in contradiction with the assumption of the claim. 
The proof is analogous for the remaining two choices for $w$, namely those with $\beta_w\gamma_w=0$. 
When $\gamma_w=0$, the third row of $Y_w$ is replaced by $(1,0,1)$; and when $\beta_w=0$, $Y_w$ is replaced by the matrix obtained from the identity matrix by swapping the first and third rows.
The $(3,3)$-entry of the matrix representing $U$ then changes to $b\alpha$ or $b(\alpha+\gamma)$ according to whether $\beta_w=0$ or $\gamma_w=0$, and the result follows as in the generic case because $\alpha=0$ and $\alpha+\gamma=0$ are the corresponding tangent lines to $\mathcal{C}$. 
This completes the proof of the claim, and of the lemma.
\end{Prf}

Recall that the $K$-orbit $\Sigma_{14}'$ arises (if and) only if $q$ is a power of $3$. 

\begin{La} \label{dist14dash}
Suppose that $q \equiv 0 \pmod 3$. A plane of type $\Sigma_{14}'$ contains 
\begin{itemize}
\item $q$ lines of type $o_{8,1}$, 
\item one line of type $o_9$, 
\item $\tfrac{q(q-1)}{6}$ lines of type $o_{14,1}$, 
\item $\tfrac{q(q-1)}{2}$ lines of type $o_{15,1}$, 
\item $q$ lines of type $o_{16}$, 
\item $\tfrac{q(q-1)}{3}$ lines of type $o_{17}$.
\end{itemize}
\end{La}

\begin{Prf}
If $M$ denotes the matrix representative of $\Sigma_{14}'$ given in Table~\ref{table:main}, then we have 
\begin{equation} \label{Sigma14-dash-matrix}
XMX^T = 
\left[ \begin{matrix}
-\beta & \cdot & \gamma-\beta \\
\cdot & \alpha+\beta & \beta \\
\gamma-\beta & \beta & \cdot
\end{matrix} \right], 
\quad \text{where} \quad 
X = \left[ \begin{matrix}
\cdot & \cdot & 1 \\
-1 & 1 & \cdot \\
\cdot & 1 & 1
\end{matrix} \right] \in \GL(3,q).
\end{equation}
Let $\pi \in \Sigma_{14}'$ denote the plane with this new representative. 
The points of rank at most $2$ in $\pi$ lie on the cubic $\mathcal{C}:\beta\gamma(\beta+\gamma)+\alpha(\beta-\gamma)^2=0$. 
There is a unique point of rank $1$, namely $x:\beta=\gamma=0$, and $q$ points of rank $2$, all exterior, one for each choice of $(\beta,\gamma)$ with $\beta\neq\gamma$. 
In particular, $q$ of the lines through $x$ have type $o_{8,1}$, and the final line through $x$ has type $o_9$. 

All of the rank-$2$ points on the cubic $\mathcal{C}$ are inflexion points. 
Let $w$ be such a point, recalling that we have already counted the line through $w$ and $x$. 
The tangent line through $w$ contains no other points of $\mathcal{C}$, so has point-orbit distribution $[0,1,0,q]$ and hence type $o_{15,1}$ or $o_{16}$. 
We show in Claim~$\Sigma_{14}'$ at the end of the proof that it has type $o_{16}$. 
This gives the total of $q$ lines of type $o_{16}$ in $\pi$. 
The remaining lines through $w$ meet a total of $q-1$ remaining points on $\mathcal{C}$. 
If such a line meets $\mathcal{C}$ in a point other than $w$, then it meets $\mathcal{C}$ in exactly three points, so has point-orbit distribution $[0,3,0,q-2]$ and (hence) type $o_{14,1}$. 
There are $\tfrac{q-1}{2}$ such lines through each $w$, and hence a total of $\tfrac{1}{3}\cdot q\cdot \tfrac{q-1}{2}$ lines of type $o_{14,1}$ in $\pi$. 
The final $\tfrac{q-1}{2}$ lines through each $w$ have point-orbit distribution $[0,1,0,q]$ and hence type $o_{15,1}$ or $o_{16}$. 
We show below that they all have type $o_{15,1}$, giving a total of $q \cdot \tfrac{q-1}{2}$ lines of type $o_{15,1}$. 
All other lines in $\pi$ have type $o_{17}$.

{\em Claim~$\Sigma_{14}'$.} Let $\ell$ be a line in $\pi$ with point-orbit distribution $[0,1,0,q]$, and let $w$ denote the unique point of rank $2$ on $\ell$. 
Then $\ell$ has type $o_{15,1}$ unless $\ell$ is the tangent line to $\mathcal{C}$ at $w$, in which case $\ell$ has type $o_{16}$.

{\em Proof of Claim~$\Sigma_{14}'$.} 
Let $(\alpha_w,\beta_w,\gamma_w)$ be the coordinates of $w$. 
Then $\beta_w \neq \gamma_w$ and $\alpha_w = -\beta_w\gamma_w(\beta_w+\gamma_w)/(\beta_w-\gamma_w)^2$. 
The tangent line to $\mathcal{C}$ at $w$ is
\begin{equation} \label{Sigma14-dash-tangent}
\alpha\cdot (\beta_w-\gamma_w)^3 - \beta \cdot \gamma_w^3 + \gamma\cdot\beta_w^3 = 0.
\end{equation}
First suppose that $\beta_w \neq 0$. 
Then the matrix representative $M_w$ of $w$ obtained from \eqref{Sigma14-dash-matrix} can be transformed as follows:
\[
Y_wM_wY_w^T = \left[ \begin{matrix}
-\beta_w & \cdot & \cdot \\
\cdot & -\frac{\beta_w^3}{(\beta_w-\gamma_w)^2} & \cdot \\
\cdot & \cdot & \cdot
\end{matrix} \right], 
\quad \text{where} \quad 
Y_w = \left[ \begin{matrix}
1 & \cdot & \cdot \\
\cdot & 1 & \cdot \\
\frac{\beta_w}{\beta_w-\gamma_w} & 1 & -\frac{\beta_w^2}{(\beta_w-\gamma_w)^2}
\end{matrix} \right] \in \GL(3,q).
\]
After making this transformation, the solid $W$ of Lemma~\ref{lemma_o15_o16} is represented by the matrices whose third rows and third columns are zero. 
If we now consider the image of a rank $3$ point on $\ell$, with coordinates $(\alpha,\beta,\gamma)$, under the element of $K$ corresponding to $Y_w$, we find that the image of $U := \langle W,\ell \rangle$ (in the notation of Lemma~\ref{lemma_o15_o16}) is represented by
\[
\left[ \begin{matrix}
d_{11} & d_{12} & b\cdot \frac{\beta_w(\beta\gamma_w-\gamma\beta_w)}{(\beta_w-\gamma_w)^2} \\
d_{21}& d_{22} & b\cdot \frac{\alpha(\beta_w-\gamma_w)^2 + \beta\gamma_w(\beta_w+\gamma_w)}{(\beta_w-\gamma_w)^2} \\
b\cdot \frac{\beta_w(\beta\gamma_w-\gamma\beta_w)}{(\beta_w-\gamma_w)^2} & b\cdot \frac{\alpha(\beta_w-\gamma_w)^2 + \beta\gamma_w(\beta_w+\gamma_w)}{(\beta_w-\gamma_w)^2} & b\cdot \frac{1}{(\beta_w-\gamma_w)^2}f_w(\alpha,\beta,\gamma)
\end{matrix} \right],
\]
where $b$ and the $d_{ij}$ range over $\mathbb{F}_q$, and $f_w(\alpha,\beta,\gamma)$ is the left-hand side of \eqref{Sigma14-dash-tangent}. 
The claim (for $\beta_w \neq 0$) therefore follows from Lemmas~\ref{lemma_o15_o16} and~\ref{lemma_o15_o16_application}. 
When $\beta_w=0$, replace $Y_w$ by the matrix obtained from the identity matrix by swapping the second and third rows. 
The $(3,3)$-entry of the matrix representing $U$ then becomes $b(\alpha+\beta)$, and the result follows because the tangent line to $w$ is the line $\alpha+\beta=0$. 
This completes the proof of the claim, and of the lemma.
\end{Prf}

Recall that the $K$-orbit $\Sigma_{14}$ arises (if and) only if $q$ is {\em not} a power of $3$. 

\begin{La} \label{dist14}
Suppose that $q \not \equiv 0 \pmod 3$. 
The lines in a plane of type $\Sigma_{14}$ are as indicated in Table~\ref{sigma14table}. 
\begin{table}[!h]
\begin{tabular}{lcc}
\toprule
Type & $q \equiv 1 \pmod{3}$ & $q \equiv -1 \pmod{3}$ \\ 
\midrule
$o_{8,1}$ & $\tfrac{q-1}{2}$ & $\tfrac{q+1}{2}$ \\ 
$o_{8,2}$ & $\tfrac{q-1}{2}$ & $\tfrac{q+1}{2}$ \\ 
$o_9$ & $2$ & $0$ \\ 
$o_{13,1}$ & $\tfrac{q-1}{2}$ & $\tfrac{q+1}{2}$ \\ 
$o_{13,2}$ & $\tfrac{q-1}{2}$ & $\tfrac{q+1}{2}$ \\ 
$o_{14,1}$ & $\tfrac{(q-1)(q-7)}{24}$ & $\tfrac{(q+1)(q-5)}{24}$ \\ 
$o_{14,2}$ & $\tfrac{(q-1)(q-3)}{8}$ & $\tfrac{(q+1)(q-1)}{8}$ \\ 
$o_{15,1}$ & $\tfrac{(q-1)^2}{4}$ & $\tfrac{(q+1)(q-3)}{4}$ \\ 
$o_{15,2}$ & $\tfrac{(q+1)(q-1)}{4}$ & $\tfrac{(q+1)(q-1)}{4}$ \\ 
$o_{17}$ & $\tfrac{(q-1)^2}{3}+q$ & $\tfrac{(q+1)^2}{3}-q$ \\
\bottomrule
\end{tabular}
\caption{Data for Lemma~\ref{dist14}.}
\label{sigma14table}
\end{table}
\end{La}

\begin{Prf}
Let $\pi$ denote the representative of $\Sigma_{14}$ given in Table~\ref{table:main}, namely
\begin{equation} \label{Sigma14matrix}
\left[ \begin{matrix}
\alpha & \beta & \cdot \\
\beta & c\gamma & \beta-\gamma \\
\cdot & \beta-\gamma & \gamma
\end{matrix} \right],
\end{equation}
where $c$ satisfies the condition ($\dagger$) given in Table~\ref{table:main}. 
That is, $c$ is some fixed element of $\bF_q \setminus \{0,1\}$ such that $-3c$ is a square in $\mathbb{F}_q$ and $\tfrac{\sqrt{c}+1}{\sqrt{c}-1}$ is a non-cube in $\mathbb{F}_q(\sqrt{-3})$. 
Note that this implies, in particular, that $c$ is a square if $q \equiv 1 \pmod 3$ and a non-square if $q \equiv -1 \pmod 3$. 
The points of rank at most $2$ in $\pi$ lie on the cubic 
\begin{equation} \label{sigma14cubicEqn}
\mathcal{C} : \alpha f_c(\beta,\gamma)-\beta^2\gamma = 0, 
\quad \text{where} \quad
f_c(\beta,\gamma) := (c-1)\gamma^2+2\beta\gamma-\beta^2.
\end{equation} 
The point-orbit distribution of $\pi$ is $[1, \frac{q\mp 1}{2}, \frac{q\mp 1}{2}, q^2\pm 1]$ according to whether $q \equiv \pm 1 \pmod 3$. 
If $q \equiv -1 \pmod 3$, the lines through the unique point $x : \beta=\gamma=0$ of rank~$1$ therefore comprise $\tfrac{q+1}{2}$ lines of each of the types $o_{8,1}$ and $o_{8,2}$. 
If $q \equiv 1 \pmod 3$, they instead comprise $\tfrac{q-1}{2}$ lines of each of these types, plus two lines of type $o_9$. 

Now consider a point $w$ of rank $2$, with coordinates $(\alpha_w,\beta_w,\gamma_w)$. Note the following:
\begin{itemize}
\item[(i)] $(\beta_w,\gamma_w) \neq (0,0)$, because $w\neq x$.
\item[(ii)] $f_c(\beta_w,\gamma_w) \neq 0$. (If not then $\beta_w\gamma_w=0$ by \eqref{sigma14cubicEqn}. If $\beta_w=0$ then $f_c(\beta_w,\gamma_w) = (c-1)\gamma_w^2$, so $\gamma_w=0$ because $c \neq 1$; similarly, $\gamma_w=0$ implies $\beta_w=0$. 
Either case contradicts (i).)
\item[(iii)] $\alpha_w = \beta_w^2\gamma_w/f_c(\beta_w,\gamma_w)$, by (ii).
\item[(iv)] $w$ is exterior if and only if $-f_c(\beta_w,\gamma_w)$ is a non-zero square in $\mathbb{F}_q$, by Lemma~\ref{lem:extcriteria}.
\end{itemize}
Moreover, as noted in the proof of \cite[Lemma~7.14]{LaPoSh2020}, the condition ($\dagger$) on $c$ is equivalent to the equation $(c-1)^2\theta^3+3(c-1)\theta+2=0$ having no solution $\theta \in \mathbb{F}_q$.
In particular, it implies that 
\begin{itemize}
\item[(v)] $\mathcal{C}$ has no points of inflexion. 
\end{itemize}
Therefore, the tangent line to $\mathcal{C}$ at $w$ contains exactly one other point of rank $2$, say $v$ with coordinates $(\alpha_v,\beta_v,\gamma_v)$ (where $\alpha_v$ is determined as in (iii)).  
We claim that 
\begin{itemize}
\item[(vi)] $v$ is exterior, regardless of whether $w$ is exterior or interior. 
\end{itemize}

{\em Proof of} (vi). In the `free' coordinates $(\alpha,\beta,\gamma)$, the tangent line to $\mathcal{C}$ at $w$ is 
\begin{equation} \label{sigma14tangentGeneric}
\alpha \cdot f_c(\beta_w,\gamma_w)^2 
- \beta \cdot 2\gamma_w^2\beta_w ((c-1)\gamma_w+\beta_w)
+ \gamma \cdot \beta_w^2 ((c-1)\gamma_w^2+\beta_w^2) = 0.
\end{equation}
Since the rank-$2$ point $v$ lies on this line, we have
\begin{equation} \label{sigma14tangent}
\beta_v\cdot 2\beta_w((c-1)\gamma_w+\beta_w) + \gamma_v\cdot (c-1)((c-1)\gamma_w^2+\beta_w^2) = 0.
\end{equation}
If $\gamma_v=0$ then $\beta_v \neq 0$ by (i), so $-f_c(\beta_v,\gamma_v) = \beta_v^2$ is a non-zero square and $v$ is therefore exterior by (iv). 
Now suppose that $\gamma_v \neq 0$. 
This implies that $\beta_w \neq 0$ and $(c-1)\gamma_w+\beta_w \neq 0$: if $\beta_w = 0$ then  \eqref{sigma14tangent} implies that $\gamma_v=0$ because $(c-1)\gamma_w \neq 0$, and if $(c-1)\gamma_w+\beta_w=0$ then $\beta_w^2=(c-1)^2\gamma_w^2$, so \eqref{sigma14tangent} reads $c(c-1)^2\gamma_w^2\gamma_v = 0$ and hence again $\gamma_v=0$. 
Therefore, if $\gamma_v\neq 0$ then \eqref{sigma14tangent} yields
\begin{equation} \label{sigma14ratio}
\frac{\beta_v}{\gamma_v} = - \frac{(c-1)((c-1)\gamma_w^2+\beta_w^2)}{2\beta_w((c-1)\gamma_w+\beta_w)},
\end{equation}
and so
\[
-f_c(\beta_v,\gamma_v) = \frac{f_c(\beta_w,\gamma_w)^2(c-1)^2}{4\beta_w^2((c-1)\gamma_w+\beta_w)^2} \cdot \gamma_v^2,
\]
which is always a non-zero square. 
Hence, (vi) holds, as claimed. 

By (vi), the tangent line to $\mathcal{C}$ at $w$ has type $o_{13,1}$ or $o_{13,2}$ according to whether $w$ is exterior or interior. 
There are $\tfrac{q\mp1}{2}$ exterior and interior rank-$2$ points in $\pi$ according to whether $q \equiv \pm 1 \pmod 3$, and hence a total of this many lines of the corresponding types. 

We now set about counting the remaining lines containing a point of rank $2$. 
Suppose first that $w_i$ is an interior rank-$2$ point. 
We have already counted the line $\langle w_i,x \rangle$ (of type $o_{8,2}$) and the tangent line to $\mathcal{C}$ at $w_i$ (of type $o_{13,2}$, containing also an exterior rank-$2$ point). 
Let $\ell$ be one of the remaining $q-1$ lines through $w$. 
By (vi), $\ell$ is not the tangent to any other point on $\mathcal{C}$. 
Hence, if $\ell$ meets $\mathcal{C}$ in a second point, then it meets $\mathcal{C}$ in exactly three points, and so contains exactly three points of rank $2$ (and no points of rank $1$). 
Since $w_i$ is interior, $\ell$ must have type $o_{14,2}$; in particular, it must contain a second interior rank-$2$ point, $w_i'$ say.
There are $\tfrac{q\mp1}{2}$ choices of $w_i$ and $\tfrac{q\mp1}{2}-1$ choices of $w_i'$, according to whether $q \equiv \pm 1 \pmod 3$, so $\pi$ contains in total $\tfrac{1}{2} \cdot \tfrac{q-1}{2} \cdot \tfrac{q-3}{2}$ lines of type $o_{14,2}$ if $q \equiv 1 \pmod 3$, and $\tfrac{1}{2} \cdot \tfrac{q+1}{2} \cdot \tfrac{q-1}{2}$ lines of type $o_{14,2}$ if $q \equiv -1 \pmod 3$. 
The remaining lines through $w_i$ have point-orbit distribution $[0,0,1,q]$ and hence type $o_{15,2}$. 
There are $(q-1)-(\tfrac{q\mp1}{2}-1) = \tfrac{q\pm1}{2}$ such lines through each $w_i$ according to whether $q \equiv \pm 1 \pmod 3$, with $\tfrac{q\mp1}{2}$ choices of $w_i$ in these respective cases. 
Hence, in either case, $\pi$ contains a total of $\tfrac{(q+1)(q-1)}{4}$ lines of type $o_{15,2}$. 

To count the lines through an {\em exterior} point of rank $2$ (apart from those already counted above), we need the following additional fact. 
\begin{itemize}
\item[(vii)] Every exterior rank-$2$ point $v$ on $\mathcal{C}$ lies on exactly three tangent lines to $\mathcal{C}$ (one of which is the tangent line at $v$). 
Moreover, if $v$ is on the tangent lines at the points $w$ and $w'$, say, then $w$ and $w'$ have the same type (that is, either both are exterior or both are interior) if and only if $q \equiv \pm 1 \pmod{12}$. 
\end{itemize}

{\em Proof of} (vii). Denote the (second and third) coordinates of $v$ by $(\beta_v,\gamma_v)$ and view \eqref{sigma14tangent} as a quadratic equation in the coordinates $(\beta_w,\gamma_w)$ of a rank-$2$ point $w$ whose tangent contains $v$:
\begin{equation} \label{Sigma14-vii}
\gamma_w^2 \cdot (c-1)^2\gamma_v + \gamma_w\beta_w \cdot 2(c-1)\beta_v + \beta_w^2 \cdot ((c-1)\gamma_v+2\beta_v) = 0.
\end{equation}
If $v$ is the point with coordinates $(-1,\frac{1-c}{2},1)$, then the coefficient of $\beta_w^2$ in \eqref{Sigma14-vii} is zero and we obtain the solutions $(\beta_w,\gamma_w) = (1,0)$ and $(\beta_w',\gamma_w')=(1,1)$, so that $w = (0,1,0)$ and $w' = (c^{-1},1,1)$ (say). 
We have $-f_c(0,1)=1$, so $w$ is exterior (for all $q$). 
On the other hand, $-f_c(1,1)=-c$, so, by (iv), $w'$ is also exterior if and only if $-c$ is a square. 
If $q \equiv 1 \pmod 3$ then $c$ is a square (because $-3c$ is a square), so $-c$ is a square if and only if $-1$ is a square, namely, if and only if $q \equiv 1 \pmod 4$, which is if and only if $q \equiv 1 \pmod{12}$ given that $q \equiv 1 \pmod 3$. 
Similarly, if $q \equiv -1 \pmod 3$ then $-c$ is a square if and only if $q \equiv -1 \pmod{12}$. 
If $v \neq (-1,\frac{1-c}{2},1)$ then the coefficient of $\beta_w^2$ in \eqref{Sigma14-vii} is non-zero, so any solution has $\gamma_w \neq 0$ and we may therefore view the left-hand side of \eqref{Sigma14-vii} as a quadratic in $\frac{\beta_w}{\gamma_w}$. 
The discriminant of this quadratic is $-f_c(\beta_v,\gamma_v) \cdot 4(c-1)^2$, which is a non-zero square by (iv), so there are again two solutions, given by
\[
\frac{\beta_w}{\gamma_w} = -(c-1) \cdot \frac{\beta_v+\sqrt{-f_c(\beta_v,\gamma_v)}}{(c-1)\gamma_v+2\beta_v} 
\quad \text{and} \quad 
\frac{\beta_w'}{\gamma_w'} = -(c-1) \cdot \frac{\beta_v-\sqrt{-f_c(\beta_v,\gamma_v)}}{(c-1)\gamma_v+2\beta_v}.
\]
A further calculation shows that
\[
(\beta_w'-\gamma_w')^2 f_c(\beta_w,\gamma_w) = -c \cdot \gamma_w^2 f_c(\beta_w',\gamma_w'). 
\]
Note at this point that $\beta_w'-\gamma_w'$ is non-zero: if it were zero then $w'$ would be the point $(c^{-1},1,1)$ and so $v$ would be the point $(-1,\frac{1-c}{2},1)$ considered previously. 
Therefore, $-f_c(\beta_w,\gamma_w)$ and $-f_c(\beta_w',\gamma_w')$ are either both squares or both non-squares if and only if $-c$ is a square. 
As noted above, this is the case if and only if $q \equiv \pm 1 \pmod{12}$. 
This completes the proof of (vii).

Now consider an exterior rank-$2$ point $w_e$. 
First suppose that $q \not \equiv \pm 1 \pmod{12}$, that is, $q \equiv 7 \text{ or } 5 \pmod{12}$. 
Recall that we have already counted the line through $w_e$ and the unique point of rank $1$ (which has type $o_{8,1}$).
By (vii), $w_e$ lies on the tangent lines to $\mathcal{C}$ at one exterior rank-$2$ point $z_e$ and one interior rank-$2$ point $z_i$. 
These tangent lines have type $o_{13,1}$ and $o_{13,2}$ respectively, and have already been counted. 
Moreover, by (vi), the tangent line to $w_e$ at $\mathcal{C}$ has type $o_{13,1}$ and has also been counted; denote the second exterior rank-$2$ point on this line by $w_e^T$.  
Any other line through $w_e$ containing an interior rank-$2$ point has type $o_{14,2}$ and has already been counted. 
There are $\tfrac{q\mp1}{2}-1$ interior rank-$2$ points other than $z_i$ according to whether $q \equiv \pm 1 \pmod 3$, and hence half this many lines of type $o_{14,2}$ through $w_e$. 
That is, there are $\tfrac{q-3}{4}$ or $\tfrac{q-1}{4}$ lines of type $o_{14,2}$ through $w_e$ according to whether $q \equiv 7 \text{ or } 5 \pmod{12}$. 
Any other line through $w_e$ containing an exterior rank-$2$ point, $w_e'$ say, has type $o_{14,1}$ and contains also a third exterior rank-$2$ point. 
For $q \equiv \pm 1 \pmod 3$, there are $\tfrac{q\mp1}{2}$ choices for $w_e$ and $\tfrac{q\mp1}{2} - 3$ choices for $w_e'$, namely all of the exterior rank-$2$ points other than $w_e$, $z_e$ and $w_e^T$. 
Hence, there are $\tfrac{1}{2} \cdot (\tfrac{q\mp1}{2} - 3)$ lines of type $o_{14,1}$ through each $w_e$ according to whether $q \equiv \pm 1 \pmod 3$, namely, $\tfrac{q-7}{4}$ and $\tfrac{q-5}{4}$ such lines in these respective cases. 
Since each line of type $o_{14,1}$ contains three exterior rank-$2$ points, $\pi$ contains in total $\tfrac{1}{3} \cdot \tfrac{q-1}{2} \cdot \tfrac{q-7}{4} = \tfrac{(q-1)(q-7)}{24}$ or $\tfrac{1}{3} \cdot \tfrac{q+1}{2} \cdot \tfrac{q-5}{4} =\tfrac{(q+1)(q-5)}{24}$ lines of type $o_{14,1}$ according to whether $q \equiv 7 \text{ or } 5 \pmod{12}$. 
The remaining lines through $w_e$ have point-orbit distribution $[0,1,0,q]$ and hence type $o_{15,1}$ or $o_{16}$. 
We show in Claim $\Sigma_{14}$ at the end of the proof that they all have type $o_{15,1}$. 
There are $(q-3) - \tfrac{q-3}{4} - \tfrac{q-7}{4} = \tfrac{q-1}{2}$ or $(q-3) - \tfrac{q-1}{4} - \tfrac{q-5}{4} = \tfrac{q-3}{2}$ such lines through each $w_e$ according to whether $q \equiv 7 \text{ or } 5 \pmod{12}$, and hence a total of $\tfrac{(q-1)^2}{4}$ or $\tfrac{(q+1)(q-3)}{4}$ such lines in $\pi$ in these respective cases. 

Now suppose that $q \equiv \pm 1 \pmod{12}$. 
Note in particular that $q\mp1$ is then divisible by $4$. 
In light of (vii), we must now consider two possibilities for the (otherwise arbitrary) exterior rank-$2$ point $w_e$, namely, whether or not $w_e$ lies on the tangent line to $\mathcal{C}$ at some interior rank-$2$ point. 
For convenience, let us say that $w_e$ is of {\em class I} if it lies on the tangent line to $\mathcal{C}$ at some interior rank-$2$ point, and of {\em class E} if it does not. 
Note that each class comprises exactly half of the $\tfrac{q\mp1}{2}$ exterior rank-$2$ points. 
In either case, recall yet again that we have already counted the line through $w_e$ and the point of rank $1$.

First consider a point $w_e$ of class~E.  
By (vii), $w_e$ lies on the tangent lines to $\mathcal{C}$ at two exterior rank-$2$ points, say $z_1$ and $z_2$. 
These tangent lines have type $o_{13,1}$ and have already been counted. 
The tangent line to $\mathcal{C}$ at $w_e$ has also already been counted; it has type $o_{13,1}$ and contains an exterior rank-$2$ point $w_e^T \not \in \{w_e,z_1,z_2\}$. 
There are $q-3$ lines through $w_e$ left to consider. 
Let $\ell$ be a line through $w_e$ containing an interior rank-$2$ point. 
Since $\ell$ is not the tangent line to $\mathcal{C}$ at $w_e$, it also contains a second interior rank-$2$ point, and has type $o_{14,2}$. 
All such lines have already been counted (because all lines through any interior rank-$2$ point have already been counted). 
There are $\tfrac{q\mp 1}{2}$ interior rank-$2$ points when $q \equiv \pm 1 \pmod{12}$, and hence half this many lines of type $o_{14,2}$ through $w_e$. 
This leaves $(q-3) - \tfrac{q-1}{4} = \tfrac{3q-11}{4}$ or $(q-3) - \tfrac{q+1}{4} = \tfrac{3q-13}{4}$ lines through $w_e$ to consider, according to whether $q \equiv 1 \text{ or } -1 \pmod{12}$. 
Now let $\ell$ be a line through $w_e$ containing an exterior rank-$2$ point $w_e' \not \in \{w_e^T,z_1,z_2\}$. 
Any such line has type $o_{14,1}$. 
There are $\tfrac{q\mp1}{2}-4$ choices for $w_e'$ (because also $w_e' \neq w_e$) and hence half this many lines of type $o_{14,1}$ through $w_e$. 
That is, there are $\tfrac{q-9}{4}$ or $\tfrac{q-7}{4}$ lines of type $o_{14,1}$ through $w_e$ according to whether $q \equiv 1 \text{ or } -1 \pmod{12}$. 
The remaining lines through $w_e$ have point-orbit distribution $[0,1,0,q]$ and hence type $o_{15,1}$ or $o_{16}$. 
Claim $\Sigma_{14}$ below implies that they are all of type $o_{15,1}$. 
There are $\tfrac{q-1}{2}$ or $\tfrac{q-3}{2}$ such lines through $w_e$ according to whether $q \equiv 1 \text{ or } -1 \pmod{12}$. 
Before counting the total number of lines of types $o_{14,1}$ and $o_{15,1}$ in $\pi$, let us analyse the points of class~I.

Suppose now that $w_e$ is of class~I. 
By (vii), $w_e$ lies on the tangent lines to $\mathcal{C}$ at two interior rank-$2$ points, say $z_1'$ and $z_2'$. 
The tangent lines have type $o_{13,2}$ and have already been counted. 
The tangent line to $\mathcal{C}$ at $w_e$ has also already been counted; it has type $o_{13,1}$ and contains an exterior rank-$2$ point $w_e^T \neq w_e$. 
There are $q-3$ lines through $w_e$ left to consider. 
Any line through $w_e$ containing one of the $\tfrac{q\mp1}{2}-2$ interior rank-$2$ points other than $z_1'$ and $z_2'$ contains two such points and has type $o_{14,2}$. 
These lines have already been counted. 
There are $\tfrac{1}{2} \cdot \tfrac{q-5}{2}$ or $\tfrac{1}{2} \cdot \tfrac{q-3}{2}$ of them through $w_e$ according to whether $q \equiv 1 \text{ or } -1 \pmod{12}$, leaving $\tfrac{3q-7}{4}$ or $\tfrac{3q-9}{4}$ lines through $w_e$ to consider in these respective cases. 
Now let $\ell$ be a line through $w_e$ containing an exterior rank-$2$ point $w_e' \neq w_e$. 
The only such line counted so far is the tangent line to $\mathcal{C}$ at $w_e$ (with $w_e'=w_e^T$), so we may assume that $\ell$ is not this line. 
Hence, $\ell$ also contains a third exterior rank-$2$ point, and has type $o_{14,1}$. 
The candidates for $w_e'$ are all of the exterior rank-$2$ points except $w_e$ and $w_e^T$, so there are $\tfrac{1}{2} ( \tfrac{q-1}{2}-2 ) = \tfrac{q-5}{4}$ or $\tfrac{1}{2} ( \tfrac{q+1}{2}-2 ) = \tfrac{q-3}{4}$ lines of type $o_{14,1}$ through $w_e$ according to whether $q \equiv 1 \text{ or } -1 \pmod{12}$.
The remaining lines through $w_e$ have point-orbit distribution $[0,1,0,q]$ and, by Claim $\Sigma_{14}$ below, type $o_{15,1}$. 
There are $\tfrac{q-1}{2}$ or $\tfrac{q-3}{2}$ such lines through $w_e$ according to whether $q \equiv 1 \text{ or } -1 \pmod{12}$.

Let us now count the {\em total} number of lines of types $o_{14,1}$ and $o_{15,1}$ in $\pi$ for $q \equiv \pm 1 \pmod{12}$. 
Consider first the lines of type $o_{15,1}$. 
As argued in the preceding two paragraphs, every exterior rank-$2$ point, whether of class~E or class~I, lies on $\tfrac{q-1}{2}$ or $\tfrac{q-3}{2}$ lines of type $o_{15,1}$ according to whether $q \equiv 1 \text{ or } -1 \pmod{12}$. 
There are $\tfrac{q\mp 1}{2}$ exterior rank-$2$ points in these respective cases, and hence a total of $\tfrac{(q-1)^2}{4}$ or $\tfrac{(q+1)(q-3)}{4}$ lines of type $o_{15,1}$. 
Now let $N$ denote the total number of lines of type $o_{14,1}$ in $\pi$. 
To calculate $N$, we count in two different ways the pairs $(w,\ell)$ with $w$ an exterior rank-$2$ point and $\ell$ a line of type $o_{14,1}$ containing $w$. 
On the one hand, there are $3N$ such pairs, because each line of type $o_{14,1}$ contains three exterior rank-$2$ points. 
On the other hand, the number of pairs $(w,\ell)$ is equal to $N_E+N_I$, where $N_E$ (respectively, $N_I$) is the number of such pairs with $w$ of class~E (respectively, class~I). 
Therefore, $N = \tfrac{1}{3}(N_E+N_I)$. 
By the calculations in the preceding two paragraphs, we have $N_E = \tfrac{q-1}{4}\tfrac{q-9}{4}$ and $N_I = \tfrac{q-1}{4}\tfrac{q-5}{4}$ when $q \equiv 1 \pmod{12}$, and $N_E = \tfrac{q+1}{4}\tfrac{q-7}{4}$ and $N_I = \tfrac{q+1}{4}\tfrac{q-3}{4}$ when $q \equiv -1 \pmod{12}$. 
This gives $N = \tfrac{(q-1)(q-7)}{24}$ and $N = \tfrac{(q+1)(q-5)}{24}$ in these respective cases. 

Every line in $\pi$ not counted thus far contains only points of rank $3$, so has type $o_{17}$. 
The number of such lines for each $q$ is recorded in the statement of the lemma. 
To complete the proof, it remains to show that all lines with point-orbit distribution $[0,1,0,q]$ have type $o_{15,1}$. 

{\em Claim $\Sigma_{14}$.} Every line in $\pi$ with point-orbit distribution $[0,1,0,q]$ is of type $o_{15,1}$.

{\em Proof of Claim $\Sigma_{14}$.} 
Let $\ell$ be such a line, and let $(\alpha_w,\beta_w,\gamma_w)$ denote the coordinates of the exterior rank-$2$ point on $\ell$. 
Recall that $\alpha_w = \beta_w^2\gamma_w/f_c(\beta_w,\gamma_w)$, where $f_c$ is defined as in \eqref{sigma14cubicEqn}, and that the tangent line to $\mathcal{C}$ at $w$ is given by \eqref{sigma14tangentGeneric}. 
We apply Lemma~\ref{lemma_o15_o16}. 
First suppose that $\beta_w(\beta_w-\gamma_w) \neq 0$. 
Then the matrix representative of $w$ obtained from \eqref{Sigma14matrix}, call it $M_w$, can be transformed via the action of $K$ as follows:
\[
Y_wM_wY_w^T = \left[ \begin{matrix}
\frac{\beta_w^2\gamma_w}{f_c(\beta_w,\gamma_w)} & \beta_w & \cdot \\
\beta_w & c\gamma_w & \cdot \\
\cdot & \cdot & \cdot
\end{matrix} \right], 
\quad \text{where} \quad 
Y_w = \left[ \begin{matrix}
1 & \cdot & \cdot \\
\cdot & 1 & \cdot \\
-\frac{f_c(\beta_w,\gamma_w)}{\beta_w} & \gamma_w & \beta_w-\gamma_w
\end{matrix} \right] \in \GL(3,q).
\]
Under this transformation, the solid $W$ of Lemma~\ref{lemma_o15_o16} is represented by the matrices whose third rows and third columns are zero. 
Now choose a point $y$ (of rank $3$) such that $\ell = \langle w,y \rangle$, denoting the coordinates of $y$ by $(\alpha,\beta,\gamma)$. 
Then the image of $y$ under the element of $K$ corresponding to the matrix $Y_w$ is represented by the matrix 
\[
\left[ \begin{matrix}
\alpha & \beta & \frac{\beta\beta_w\gamma_w - \alpha f_c(\beta_w,\gamma_w)}{\beta_w\gamma_w} \\
\beta & c\gamma & \frac{(\gamma\beta_w-\beta\gamma_w)((c-1)\gamma_w+\beta_w)}{\beta_w\gamma_w} \\
\frac{\beta\beta_w\gamma_w - \alpha f_c(\beta_w,\gamma_w)}{\beta_w\gamma_w} & \frac{(\gamma\beta_w-\beta\gamma_w)((c-1)\gamma_w+\beta_w)}{\beta_w\gamma_w} & \frac{1}{\beta_w^4}g_w(\alpha,\beta,\gamma)
\end{matrix} \right],
\]
where $g_w(\alpha,\beta,\gamma)$ is the left-hand side of \eqref{sigma14tangentGeneric}.
Hence, in the notation of Lemma~\ref{lemma_o15_o16}, the image of $U := \langle W,\ell \rangle$ is represented by 
\[
\left[ \begin{matrix}
d_{11} & d_{12} & b \cdot \frac{\beta\beta_w\gamma_w - \alpha f_c(\beta_w,\gamma_w)}{\beta_w\gamma_w} \\
d_{21} & d_{22} & b \cdot \frac{(\gamma\beta_w-\beta\gamma_w)((c-1)\gamma_w+\beta_w)}{\beta_w\gamma_w} \\
b \cdot \frac{\beta\beta_w\gamma_w - \alpha f_c(\beta_w,\gamma_w)}{\beta_w\gamma_w} & b \cdot \frac{(\gamma\beta_w-\beta\gamma_w)((c-1)\gamma_w+\beta_w)}{\beta_w\gamma_w} & b \cdot \frac{1}{\beta_w^4}g_w(\alpha,\beta,\gamma)
\end{matrix} \right],
\]
where $b$ and the $d_{ij}$ range over $\mathbb{F}_q$. 
Lemma~\ref{lemma_o15_o16} now tells us that $\ell$ has type $o_{15,1}$ if and only if there exist $d_{ij}$ such that the above matrix has rank $1$ when $b \neq 0$. 
Lemma~\ref{lemma_o15_o16_application} therefore implies that $\ell$ has type $o_{15,1}$ unless $g_w(\alpha,\beta,\gamma)=0$, that is, unless $y$ lies on the tangent line to $\mathcal{C}$ at $w$. 
The latter condition is impossible, because the tangent line to $\mathcal{C}$ at $w$ contains two points of rank $2$, whereas, by assumption, $\ell$ contains a unique point of rank $2$.
The proof is analogous if $\beta_w(\beta_w-\gamma_w)=0$. 
If $\beta_w=0$, respectively $\beta_w-\gamma_w=0$, replace $Y_w$ by 
\[
\left[ \begin{matrix}
\cdot & \cdot & 1 \\
\cdot & 1 & \cdot \\
1 & \cdot & \cdot
\end{matrix} \right], 
\quad \text{respectively} \quad 
\left[ \begin{matrix}
c^{-1} & 1 & \cdot \\
1 & c & \cdot \\
\cdot & \cdot & 1
\end{matrix} \right].
\]
The $(3,3)$-entry of the matrix representing $U$ then changes to $b\alpha$ or $b(\alpha c^2-2\beta c +\gamma c)$, respectively, and the result follows as in the generic case because $\alpha=0$ and $\alpha c^2-2\beta c +\gamma c=0$ are the corresponding tangent lines to $\mathcal{C}$. 
This completes the proof of the claim, and of the lemma.
\end{Prf}

\begin{La} \label{dist13}
The lines in a plane of type $\Sigma_{13}$ are as indicated in Table~\ref{sigma13table}. 
\begin{table}[!h]
\begin{tabular}{lcc}
\toprule
Type & $q \equiv -1 \pmod 3$ & $q \not \equiv -1 \pmod 3$ \\
\midrule
$o_{8,1}$ & $\tfrac{q+1}{2}$ & $\tfrac{q+1}{2}$ \\
$o_{8,2}$ & $\tfrac{q+1}{2}$ & $\tfrac{q+1}{2}$ \\
$o_{13,1}$ & $\tfrac{q-5}{2}$ & $\tfrac{q-1}{2}$ \\
$o_{13,2}$ & $\tfrac{q+1}{2}$ & $\tfrac{q+1}{2}$ \\
$o_{14,1}$ & $\tfrac{(q+1)(q-5)}{24}+1$ & $\tfrac{(q-1)(q-3)}{24}$ \\
$o_{14,2}$ & $\tfrac{(q+1)(q-1)}{8}$ & $\tfrac{(q+1)(q-1)}{8}$ \\
$o_{15,1}$ & $\tfrac{(q+1)(q-3)}{4}$ & $\tfrac{(q+1)(q-3)}{4}$ \\
$o_{15,2}$ & $\tfrac{(q+1)(q-1)}{4}$ & $\tfrac{(q+1)(q-1)}{4}$ \\
$o_{16}$ & $3$ & $1$ \\
$o_{17}$ & $\tfrac{(q+1)(q-2)}{3}$ & $\tfrac{q(q-1)}{3}$ \\
\bottomrule
\end{tabular}
\caption{Data for Lemma~\ref{dist13}.}
\label{sigma13table}
\end{table}
\end{La}

\begin{Prf}
Let $\pi$ denote the representative of $\Sigma_{13}$ given in Table~\ref{table:main}, namely
\begin{equation} \label{eq:sig13rep}
\left[ \begin{matrix}
\alpha & \beta & \cdot \\
\beta & \gamma & \beta \\
\cdot & \beta & \varepsilon\gamma
\end{matrix} \right], 
\quad \text{where } \varepsilon \text{ is a non-square in } \mathbb{F}_q.
\end{equation}
The point-orbit distribution of $\pi$ is $[1, \frac{q+1}{2}, \frac{q+1}{2}, q^2-1]$. 
The points of rank at most $2$ in $\pi$ lie on the cubic $\mathcal{C} : \alpha (\varepsilon\gamma^2-\beta^2) - \varepsilon\beta^2\gamma = 0$.  
Hence, there are $\tfrac{q+1}{2}$ lines of each of the types $o_{8,1}$ and $o_{8,2}$ through the unique point $x : \beta=\gamma=0$ of rank~$1$. 

The inflexion points of $\mathcal{C}$ are precisely the rank-$2$ points satisfying $\gamma(3\beta^2+\varepsilon\gamma^2)=0$. 
The point $y = (0,1,0)$ is an inflexion point for all $q$. 
If $q$ is a power of $3$, it is the only inflexion point, because $3\beta^2=0$. 
If $q \equiv 1 \pmod 3$ then $-3$ is a square in $\mathbb{F}_q$, so $-3\varepsilon$ is a non-square and the hence equation $3\beta^2+\varepsilon\gamma^2=0$ has no (non-zero) solutions. 
Therefore, $y$ is again the only inflexion point. 
If $q \equiv -1 \pmod 3$ then $-3$ is a non-square and there are two further inflexion points, namely $y_\pm = (-\tfrac{3\varepsilon}{4}, \pm \sqrt{-3\varepsilon}, 3)$.

Now consider a point $w$ of rank $2$, with coordinates $(\alpha_w,\beta_w,\gamma_w)$. 
Note the following:
\begin{itemize}
\item[(i)] $(\beta_w,\gamma_w) \neq (0,0)$, because $w\neq x$. 
\item[(ii)] $\varepsilon\gamma_w^2-\beta_w^2 \neq 0$, by (i). 
\item[(iii)] $\alpha_w = \varepsilon\beta_w^2\gamma_w/(\varepsilon\gamma_w^2-\beta_w^2)$, by (ii). 
\item[(iv)] $w$ is exterior if and only if $\beta_w^2-\varepsilon\gamma_w^2$ is a non-zero square in $\mathbb{F}_q$ (by Lemma~\ref{lem:extcriteria}); in particular, all inflexion points of $\mathcal{C}$ are exterior. 
\end{itemize}
If $w$ is not an inflexion point then the tangent line to $\mathcal{C}$ at $w$ contains exactly one other point of rank $2$, say $v$ with coordinates $(\alpha_v,\beta_v,\gamma_v)$, where $\alpha_v$ is given by (iii). 
We claim that 
\begin{itemize}
\item[(vi)] $v$ is exterior, regardless of whether $w$ is exterior or interior. 
\end{itemize}

{\em Proof of} (vi). In the `free' coordinates $(\alpha,\beta,\gamma)$, the tangent line to $\mathcal{C}$ at $w$ is 
\begin{equation} \label{eq:sig13tan}
\alpha \cdot (\varepsilon\gamma_w^2-\beta_w^2)^2 
- \beta \cdot 2\varepsilon^2\gamma_w^3\beta_w 
+ \gamma \cdot \varepsilon\beta_w^2 (\varepsilon\gamma_w^2+\beta_w^2) = 0.
\end{equation}
Since the rank-$2$ point $v$ lies on this line, we have
\begin{equation} \label{sigma13tangent}
\beta_v\cdot 2\beta_w\gamma_w + \gamma_v\cdot(\varepsilon\gamma_w^2+\beta_w^2) = 0,
\end{equation}
so $\beta_v^2-\varepsilon\gamma_v^2 = \beta_v^2$ or $\frac{(\beta_w^2-\varepsilon\gamma_w^2)^2}{(2\beta_w\gamma_w)^2} \cdot \gamma_v^2$ according to whether $\gamma_v=0$ or not, and hence (iv) implies~(vi). 

If $w$ is not an inflexion point then, by (vi), the tangent line to $\mathcal{C}$ at $w$ has type $o_{13,1}$ or $o_{13,2}$ according to whether $w$ is exterior or interior. 
By (iv), all inflexion points are exterior, so there are $\tfrac{q+1}{2}$ lines of type $o_{13,2}$ in $\pi$ for all $q$. 
The number of lines of type $o_{13,1}$ in $\pi$ is $\tfrac{q+1}{2}$ minus the number of inflexion points, namely $\tfrac{q+1}{2}-3 = \tfrac{q-5}{2}$ or $\tfrac{q+1}{2}-1 = \tfrac{q-1}{2}$ according to whether $q \equiv -1 \pmod 3$ or not. 
The tangent line through a point of inflexion has point-orbit distribution $[0,1,0,q-1]$ and hence type $o_{15,1}$ or $o_{16}$. 
Claim~$\Sigma_{13}$ below shows that it has type $o_{16}$, giving three lines of type $o_{16}$ in $\pi$ when $q \equiv -1 \pmod 3$ and one such line otherwise.

We now count the remaining lines containing a point of rank $2$. 
Suppose first that $w_i$ is an interior rank-$2$ point. 
We have already counted the line $\langle w_i,x \rangle$ (of type $o_{8,2}$) and the tangent line to $\mathcal{C}$ at $w_i$ (of type $o_{13,2}$). 
Let $\ell$ be one of the remaining $q-1$ lines through $w_i$. 
By (vi), $\ell$ is not the tangent to any other point on $\mathcal{C}$, so if it meets $\mathcal{C}$ in a second point then it meets $\mathcal{C}$ in exactly three points. 
Since $w_i$ is interior, this implies that $\ell$ has type $o_{14,2}$, and in particular that it contains exactly one other interior rank-$2$ point, $w_i'$, say. 
Each of the $\tfrac{q+1}{2}$ choices of $w_i$ therefore lies on $\tfrac{q+1}{2}-1=\tfrac{q-1}{2}$ lines of type $o_{14,2}$ (that is, the number of choices of $w_i'$), and $\pi$ contains $\tfrac{1}{2} \cdot \tfrac{q+1}{2} \cdot \tfrac{q-1}{2}$ such lines in total. 
The remaining $(q-1)-\tfrac{q-1}{2} = \tfrac{q-1}{2}$ lines through $w_i$ have type $o_{15,2}$, so $\pi$ contains a total of $\tfrac{(q+1)(q-1)}{4}$ such lines. 

To count the remaining lines through an {\em exterior} rank-$2$ point, we first note the following additional fact, which suggests that separate arguments will be required to treat the cases $q \equiv \pm 1 \pmod 4$. 
Recall here that all inflexion points of $\mathcal{C}$ are exterior. 
\begin{itemize}
\item[(vii)] An exterior rank-$2$ point $v$ on $\mathcal{C}$ lies on exactly two or three tangent lines to $\mathcal{C}$ (one of which is the tangent line at $v$), according to whether $v$ is an inflexion point or not. 
If $v$ is an inflexion point and lies on the tangent line to $\mathcal{C}$ at the point $w$, say, then $w$ is exterior if and only if $q \equiv -1 \pmod 4$. 
If $v$ is not an inflexion point and lies on the tangent lines at the points $w$ and $w'$, then $w$ and $w'$ have the same type (both exterior or both interior) if and only if $q \equiv -1 \pmod 4$. 
\end{itemize}

{\em Proof of} (vii). Denote the (second and third) coordinates of $v$ by $(\beta_v,\gamma_v)$ and view \eqref{sigma13tangent} as a quadratic equation in the coordinates $(\beta_w,\gamma_w)$ of a rank-$2$ point $w$ whose tangent contains $v$:
\begin{equation} \label{Sigma13-vii}
\gamma_w^2 \cdot \varepsilon\gamma_v + \gamma_w\beta_w \cdot 2\beta_v + \beta_w^2 \cdot \gamma_v = 0.
\end{equation}
Suppose first that $v$ is an inflexion point. 
If $v=y=(0,1,0)$ then $w=(0,0,1)$, which is exterior if and only if $-\varepsilon$ is a square. 
Since $\varepsilon$ is a non-square, this is the case if and only if $-1$ is a non-square, namely if and only if $q \equiv -1 \pmod 4$. 
If $q \equiv -1 \pmod 3$ and $v=y_\pm=(-\tfrac{3\varepsilon}{4}, \pm \sqrt{-3\varepsilon}, 3)$, then $w=(-\tfrac{3\varepsilon}{4},\mp\sqrt{-3\varepsilon},1)$, which is exterior if and only if only if $(\mp\sqrt{-3\varepsilon})^2-1^2=-4\varepsilon$ is a square, namely if and only if $-1$ is a non-square. 
Now suppose that $v$ is not an inflexion point. 
Then in particular $v \neq y$, so $\gamma_v \neq 0$, and so any solution $(\beta_w,\gamma_w)$ of \eqref{Sigma13-vii} has $\gamma_w \neq 0$, because $\gamma_w=0$ would imply $\beta_w=0$, contradicting (i). 
Hence, we may view \eqref{Sigma13-vii} as a quadratic in $\frac{\beta_w}{\gamma_w}$. 
The discriminant of this quadratic is $4(\beta_v^2-\varepsilon\gamma_v^2)$, which is a non-zero square by (iv), so there are two solutions $(\beta_w,\gamma_w)$ and $(\beta_w',\gamma_w')$, given by
\[
\frac{\beta_w}{\gamma_w} = \frac{-\beta_v+\sqrt{\beta_v^2-\varepsilon\gamma_v^2}}{\gamma_v} 
\quad \text{and} \quad 
\frac{\beta_w'}{\gamma_w'} = \frac{-\beta_v-\sqrt{\beta_v^2-\varepsilon\gamma_v^2}}{\gamma_v}.
\]
A further calculation shows that
\begin{equation} \label{sig13(vii)generic}
(\beta_w^2-\varepsilon\gamma_w^2)((\beta_w')^2-\varepsilon(\gamma_w')^2) = -4\varepsilon \cdot \frac{\gamma_w^2(\gamma_w')^2}{\gamma_v^2} \cdot (\beta_v^2-\varepsilon\gamma_v^2),
\end{equation}
so (iv) implies that $w$ and $w'$ are of the same type (exterior or interior) if and only if $-1$ is a non-square.
This completes the proof of (vii).

Suppose now that $q \equiv 1 \pmod 4$. 
Consider an exterior rank-$2$ point $w_e$, and recall that we have already counted the line of type $o_{8,1}$ through $w_e$ and the unique point of rank $1$. 
Suppose first that $w_e$ is an inflexion point. 
By (vii), $w_e$ lies on the tangent line to $\mathcal{C}$ at a unique interior rank-$2$ point $z_i$. 
This tangent line has type $o_{13,2}$ and has already been counted. 
The tangent line to $w_e$ at $\mathcal{C}$, which has type $o_{16}$, has also already been counted, so there are $(q+1)-3=q-2$ lines through $w_e$ left to consider. 
Any other line through $w_e$ containing an interior rank-$2$ point has type $o_{14,2}$ and has already been counted. 
There are $\tfrac{q+1}{2}-1=\tfrac{q-1}{2}$ interior rank-$2$ points other than $z_i$, and hence half this many lines of type $o_{14,2}$ through $w_e$. 
This leaves $(q-2)-\tfrac{q-1}{4}=\tfrac{3q-7}{4}$ lines through $w_e$ to consider. 
Any other line through $w_e$ containing a second exterior rank-$2$ point, $w_e'$ say, has type $o_{14,1}$ and contains also a third exterior rank-$2$ point. 
There are $\tfrac{q+1}{2} - 1$ choices for $w_e'$, and therefore $\tfrac{q-1}{4}$ lines of type $o_{14,1}$ through $w_e$. 
The remaining $\tfrac{3q-7}{4}-\tfrac{q-1}{4}=\tfrac{q-3}{2}$ lines through $w_e$ have point-orbit distribution $[0,1,0,q]$ and hence type $o_{15,1}$ or $o_{16}$. 
Claim~$\Sigma_{13}$ below shows that they are all of type $o_{15,1}$. 
We count the total number of lines of types $o_{14,1}$ and $o_{15,1}$ in $\pi$ below, after considering also the non-inflexion points. 

Assuming still that $q \equiv 1 \pmod 4$, suppose now that $w_e$ is not an inflexion point. 
There are $\tfrac{q+1}{2}-3=\tfrac{q-5}{2}$ or $\tfrac{q+1}{2}-1=\tfrac{q-1}{2}$ choices of $w_e$ depending on whether $q \equiv -1 \pmod 3$ or not. 
In either case, by (vii), $w_e$ lies on the tangent lines to $\mathcal{C}$ at one exterior rank-$2$ point $z_e$ and one interior rank-$2$ point $z_i$. 
(Note here that $z_e$ is not an inflexion point, by (vii).) 
These tangent lines have type $o_{13,1}$ and $o_{13,2}$, respectively, and have already been counted. 
By (vi), the tangent line to $w_e$ at $\mathcal{C}$ has type $o_{13,1}$ and has also been counted; denote the second exterior rank-$2$ point on this line by $w_e^T$. 
Any other line through $w_e$ containing an interior rank-$2$ point has type $o_{14,2}$ and has already been counted. 
There are $\tfrac{q-1}{2}$ interior rank-$2$ points other than $z_i$, and half this many lines of type $o_{14,2}$ through $w_e$, leaving $(q+1)-4-\tfrac{q-1}{4}=\tfrac{3q-11}{4}$ lines through $w_e$ to consider. 
Any other line through $w_e$ containing a second exterior rank-$2$ point has type $o_{14,1}$ and contains also a third exterior rank-$2$ point. 
These lines meet a total of $\tfrac{q+1}{2} - 3 = \tfrac{q-5}{2}$ exterior rank-$2$ points, namely those other than $w_e$, $z_e$ and $w_e^T$, so there are $\tfrac{q-5}{4}$ of them through $w_e$. 
The remaining $\tfrac{3q-11}{4} - \tfrac{q-5}{4} = \tfrac{q-3}{2}$ lines through $w_e$ have point-orbit distribution $[0,1,0,q]$ and hence type $o_{15,1}$ or $o_{16}$. 
Claim~$\Sigma_{13}$ below implies that they all have type $o_{15,1}$. 

We can now determine the total number of lines of types $o_{14,1}$ and $o_{15,1}$ in $\pi$ in the case $q \equiv 1 \pmod 4$. 
As explained above, every exterior rank-$2$ point, whether an inflexion point or not, lies on $\tfrac{q-3}{2}$ lines of type $o_{15,1}$, so $\pi$ contains $\tfrac{(q+1)(q-3)}{4}$ such lines in total. 
Now let $N$ denote the total number of lines of type $o_{14,1}$ in $\pi$. 
To calculate $N$, we count in two different ways the number of incident point--line pairs $(w,\ell)$ with $w$ an exterior rank-$2$ point and $\ell$ a line of type $o_{14,1}$. 
On the one hand, there are $3N$ such pairs. 
On the other hand, the number of pairs $(w,\ell)$ is equal to $N_i+N_i'$, where $N_i$ (respectively $N_i'$) is the number of such pairs with $w_e$ an inflexion (respectively non-inflexion) point. 
Each inflexion point lies on $\tfrac{q-1}{4}$ lines of type $o_{14,1}$, so we have $N_i = \tfrac{3(q-1)}{4}$ or $\tfrac{q-1}{4}$ according to whether $q \equiv -1 \pmod 3$ or not. 
Each non-inflexion point lies on $\tfrac{q-5}{4}$ such lines, so $N_i' = \tfrac{q-5}{2} \cdot \tfrac{q-5}{4}$ or $\tfrac{q-1}{2} \cdot \tfrac{q-5}{4}$ in these respective cases. 
Calculating $N=\tfrac{1}{3}(N_i+N_i')$ in each case gives the total number of lines of type $o_{15,1}$ in $\pi$. 

To complete the proof in the case $q \equiv 1 \pmod 4$, it remains to note that all lines not counted thus far have type $o_{17}$ (and that $\pi$ contains $q^2+q+1$ lines in total).

Now suppose that $q \equiv -1 \pmod 4$. 
Consider an exterior rank-$2$ point $w_e$, and recall that we have already counted the line (of type $o_{8,1}$) through $w_e$ and the unique point of rank $1$. 
Suppose first that $w_e$ is an inflexion point. 
By (vii), $w_e$ lies on the tangent line to $\mathcal{C}$ at a unique exterior rank-$2$ point $z_e$. 
This tangent line has already been counted, as has the tangent line to $w_e$ at $\mathcal{C}$, so there are $(q+1)-3=q-2$ lines through $w_e$ left to consider. 
Any line through $w_e$ containing an interior rank-$2$ point has type $o_{14,2}$ and has already been counted. 
There are $\tfrac{q+1}{4}$ such lines through $w_e$, leaving $(q-2)-\tfrac{q+1}{4}=\tfrac{3q-9}{4}$ lines through $w_e$ to consider. 
Apart from the line $\langle w_e,z_e \rangle$, any other line through $w_e$ containing a second exterior rank-$2$ point, $w_e'$ say, has type $o_{14,1}$ and contains also a third exterior rank-$2$ point. 
There are $\tfrac{q+1}{2} - 2=\tfrac{q-3}{2}$ choices for $w_e'$, and half this many lines of type $o_{14,1}$ through $w_e$. 
The remaining $\tfrac{3q-9}{4}-\tfrac{q-3}{4}=\tfrac{q-3}{2}$ lines through $w_e$ have point-orbit distribution $[0,1,0,q]$ and hence type $o_{15,1}$ or $o_{16}$. 
By Claim~$\Sigma_{13}$ below, they are all of type $o_{15,1}$. 
As in the $q \equiv 1 \pmod 4$ case, we delay counting the total number of lines of types $o_{14,1}$ and $o_{15,1}$ until we have also considered the non-inflexion points. 

Now suppose that $w_e$ is not an inflexion point. 
By (vii), there are two possibilities to consider, namely whether $w_e$ lies on the tangent line to $\mathcal{C}$ at some interior rank-$2$ point or not. 
We say that $w_e$ is of {\em class~I} or {\em class~E} in these respective cases. 
There are $\tfrac{q+1}{4}$ points of class~I, each lying on the tangent lines to $\mathcal{C}$ at two of the $\tfrac{q+1}{2}$ interior rank-$2$ points. 
The remaining non-inflexion points $w_e$ are of class~E; there are $\tfrac{q+1}{4}-3=\tfrac{q-11}{4}$ or $\tfrac{q+1}{4}-1=\tfrac{q-3}{4}$ of them according to whether $q \equiv -1 \pmod 3$ or not. 
Regardless of the class of $w_e$, recall yet again that we have already counted the line through $w_e$ and the point of rank $1$. 

First consider a point $w_e$ of class~E.  
By (vii), $w_e$ lies on the tangent lines to $\mathcal{C}$ at two exterior rank-$2$ points, say $z_1$ and $z_2$. 
These tangent lines have already been counted. 
The tangent line to $\mathcal{C}$ at $w_e$ has also already been counted; it has type $o_{13,1}$ and contains an exterior rank-$2$ point $w_e^T \not \in \{w_e,z_1,z_2\}$. 
There are $q-3$ lines through $w_e$ left to consider. 
Let $\ell$ be a line through $w_e$ containing an interior rank-$2$ point. 
Since $\ell$ is not the tangent line to $\mathcal{C}$ at $w_e$, it has type $o_{14,2}$. 
There are $\tfrac{q+1}{4}$ such lines through $w_e$, leaving $(q-3) - \tfrac{q+1}{4} = \tfrac{3q-13}{4}$ lines through $w_e$ to consider. 
Now let $\ell$ be a line through $w_e$ containing an exterior rank-$2$ point $w_e' \not \in \{w_e^T,z_1,z_2\}$. 
Any such line has type $o_{14,1}$. 
There are $\tfrac{q+1}{2}-4=\tfrac{q-7}{2}$ exterior rank-$2$ points other than $w_e$, $w_e^T$, $z_1$ and $z_2$, and therefore $\tfrac{q-7}{4}$ lines of type $o_{14,1}$ through $w_e$. 
The remaining $\tfrac{q-3}{2}$ lines through $w_e$ have point-orbit distribution $[0,1,0,q]$, and Claim~$\Sigma_{13}$ shows that they are all of type $o_{15,1}$. 
Again, we delay counting the total number of lines of types $o_{14,1}$ and $o_{15,1}$. 

Finally, consider a point $w_e$ of class~I.  
By (vii), $w_e$ lies on the tangent lines to $\mathcal{C}$ at two interior rank-$2$ points, say $z_1$ and $z_2$. 
These tangent lines have already been counted. 
The tangent line to $\mathcal{C}$ at $w_e$ has also already been counted; it has type $o_{13,1}$ and contains an exterior rank-$2$ point $w_e^T \neq w_e$. 
There are $q-3$ lines through $w_e$ left to consider. 
Let $\ell$ be a line through $w_e$ containing an interior rank-$2$ point. 
Since $\ell$ is not the tangent line to $\mathcal{C}$ at $w_e$, it has type $o_{14,2}$. 
There are $\tfrac{q+1}{2}-2=\tfrac{q-3}{2}$ interior rank-$2$ points other than $z_1$ and $z_2$, and hence half this many lines of type $o_{14,2}$ through $w_e$. 
This leaves $(q-3) - \tfrac{q-3}{4} = \tfrac{3q-9}{4}$ lines through $w_e$ to consider. 
Now let $\ell$ be a line through $w_e$ containing an exterior rank-$2$ point $w_e' \neq w_e^T$. 
Any such line has type $o_{14,1}$. 
There are $\tfrac{q+1}{2}-2=\tfrac{q-3}{2}$ exterior rank-$2$ points other than $w_e$ and $w_e^T$, and therefore $\tfrac{q-3}{4}$ lines of type $o_{14,1}$ through $w_e$. 
The remaining $\tfrac{q-3}{2}$ lines through $w_e$ have point-orbit distribution $[0,1,0,q]$, and Claim~$\Sigma_{13}$ shows that they are all of type $o_{15,1}$. 

Let us now finally count the total number of lines of types $o_{14,1}$ and $o_{15,1}$ in $\pi$ in the case $q \equiv -1 \pmod 4$. 
As explained in the preceding arguments, every exterior rank-$2$ point $w_e$ lies on $\tfrac{q-3}{2}$ lines of type $o_{15,1}$, so $\pi$ contains a total of $\tfrac{(q+1)(q-3)}{4}$ such lines. 
Now let $N$ denote the number of lines of type $o_{14,1}$ in $\pi$. 
To calculate $N$, we count in two different ways the incident point--line pairs $(w,\ell)$ with $w$ an exterior rank-$2$ point and $\ell$ a line of type $o_{14,1}$. 
This number equals $3N$, so we have $N = \tfrac{1}{3} (N_i + N_i' + N_i'')$, where $N_i$ is the number of such pairs with $w$ an inflexion point, $N_i'$ is the number of such pairs with $w$ a non-inflexion point of class~E, and $N_i''$ is the number of such pairs with $w$ a non-inflexion point of class~I. 
Each inflexion point lies on $\tfrac{q-1}{4}$ lines of type $o_{14,1}$, so $N_i = \tfrac{3(q-3)}{4}$ or $\tfrac{q-3}{4}$ according to whether $q \equiv -1 \pmod 3$ or not. 
There are $\tfrac{q-11}{4}$ or $\tfrac{q-3}{4}$ points of class~E in these respective cases, and each such point lies on $\tfrac{q-7}{4}$ lines of type $o_{14,1}$, so $N_i' = \tfrac{q-11}{4} \cdot \tfrac{q-7}{4}$ or $\tfrac{q-3}{4} \cdot \tfrac{q-7}{4}$. 
Finally, there are $\tfrac{q+1}{4}$ points of class~I, each lying on $\tfrac{q-3}{4}$ lines of type $o_{14,1}$, so $N_i'' = \tfrac{q+1}{4} \cdot \tfrac{q-3}{4}$ regardless of whether $q \equiv -1 \pmod 3$ or not. 
The asserted values of $N$ are now readily calculated. 

To complete the proof in the case $q \equiv 1 \pmod 4$, observe that all lines not counted thus far have type $o_{17}$.
It remains to prove the following:

{\em Claim~$\Sigma_{13}$.} Let $\ell$ be a line in $\pi$ with point-orbit distribution $[0,1,0,q]$, and let $w$ denote the unique point of rank $2$ on $\ell$. 
Then $\ell$ has type $o_{15,1}$ unless $w$ is an inflexion point and $\ell$ is the tangent line to $\mathcal{C}$ at $w$, in which case $\ell$ has type $o_{16}$.

{\em Proof of Claim~$\Sigma_{13}$.} 
Let $(\alpha_w,\beta_w,\gamma_w)$ be the coordinates of $w$. 
Recall that $\varepsilon\gamma_w^2-\beta_w^2 \neq 0$ and $\alpha_w = \beta_w^2\gamma_w/(\varepsilon\gamma_w^2-\beta_w^2)$, and that the tangent line to $\mathcal{C}$ at $w$ is given by \eqref{eq:sig13tan}.
If $\beta_w \neq 0$ then the matrix representative $M_w$ of $w$ obtained from \eqref{eq:sig13rep} can be transformed as follows:
\[
Y_wM_wY_w^T = \left[ \begin{matrix}
-\frac{\varepsilon\beta_w^2\gamma_w}{\beta_w^2-\varepsilon\gamma_w^2} & \beta_w & \cdot \\
\beta_w & \gamma_w & \cdot \\
\cdot & \cdot & \cdot
\end{matrix} \right], 
\quad \text{where} \quad 
Y_w = \left[ \begin{matrix}
1 & \cdot & \cdot \\
\cdot & 1 & \cdot \\
\beta_w^2-\varepsilon\gamma_w^2 & \varepsilon\beta_w\gamma_w & -\beta_w^2
\end{matrix} \right] \in \GL(3,q).
\]
The solid $W$ of Lemma~\ref{lemma_o15_o16} is then represented by the matrices whose third rows and third columns are zero. 
If we now consider the image of a rank $3$ point on $\ell$, with coordinates $(\alpha,\beta,\gamma)$, under the element of $K$ corresponding to $Y_w$, we find that the image of $U := \langle W,\ell \rangle$ (in the notation of Lemma~\ref{lemma_o15_o16}) is represented by
\[
\left[ \begin{matrix}
d_{11} & d_{12} & b\cdot (\alpha(\beta_w^2-\varepsilon\gamma_w^2)+\beta\varepsilon\beta_w\gamma_w)) \\
d_{21}& d_{22} & b\cdot \varepsilon(\gamma\beta_w\gamma_w-\beta\gamma_w^2) \\
b\cdot(\alpha(\beta_w^2-\varepsilon\gamma_w^2)+\beta\varepsilon\beta_w\gamma_w)) & b\cdot \varepsilon(\gamma\beta_w\gamma_w-\beta\gamma_w^2) & b\cdot f_w(\alpha,\beta,\gamma)
\end{matrix} \right],
\]
where $b$ and the $d_{ij}$ range over $\mathbb{F}_q$, and $f_w(\alpha,\beta,\gamma)$ is the left-hand side of \eqref{eq:sig13tan}. 
Lemmas~\ref{lemma_o15_o16} and~\ref{lemma_o15_o16_application} therefore imply that $\ell$ has type $o_{15,1}$ unless $f_w(\alpha,\beta,\gamma)=0$, that is, unless $y$ lies on the tangent line to $\mathcal{C}$ at $w$. 
However, $\ell$ cannot be the tangent line to $\mathcal{C}$ at $w$ unless $w$ is an inflexion point, because in all other cases this line has point-orbit distribution $[0,2,0,q-1]$, in contradiction with the assumption of the claim. 
When $\beta_w=0$, we replace $Y_w$ by the matrix obtained from the identity matrix by swapping the second and third rows. 
The $(3,3)$-entry of the matrix representing $U$ then becomes $b\alpha$, and the result follows because the tangent line to $\mathcal{C}$ at $w$ is the line $\alpha=0$. 
This completes the proof of the claim, and of the lemma.
\end{Prf}

\begin{La} \label{dist12}
The lines in a plane of type $\Sigma_{12}$ are as indicated in Table~\ref{sigma12table}.
\begin{table}[!h]
\begin{tabular}{lcc}
\toprule
Type & $q \equiv 1 \pmod 3$ & $q \not \equiv 1 \pmod 3$ \\
\midrule
$o_{8,1}$ & $\tfrac{q-1}{2}$ & $\tfrac{q-1}{2}$ \\
$o_{8,2}$ & $\tfrac{q-1}{2}$ & $\tfrac{q-1}{2}$ \\
$o_9$ & $2$ & $2$ \\
$o_{13,1}$ & $\tfrac{q-7}{2}$ & $\tfrac{q-3}{2}$ \\
$o_{13,2}$ & $\tfrac{q-1}{2}$ & $\tfrac{q-1}{2}$ \\
$o_{14,1}$ & $\tfrac{(q-1)(q-7)}{24}+1$ & $\tfrac{(q-3)(q-5)}{24}$ \\
$o_{14,2}$ & $\tfrac{(q-1)(q-3)}{8}$ & $\tfrac{(q-1)(q-3)}{8}$ \\
$o_{15,1}$ & $\tfrac{(q-1)^2}{4}$ & $\tfrac{(q-1)^2}{4}$ \\
$o_{15,2}$ & $\tfrac{(q+1)(q-1)}{4}$ & $\tfrac{(q+1)(q-1)}{4}$ \\
$o_{16}$ & $3$ & $1$ \\
$o_{17}$ & $\tfrac{(q-1)(q+2)}{3}$ & $\tfrac{q(q+1)}{3}$ \\
\bottomrule
\end{tabular}
\caption{Data for Lemma~\ref{dist12}.}
\label{sigma12table}
\end{table}
\end{La}

\begin{Prf}
Let $\pi$ denote the representative of $\Sigma_{12}$ given in Table~\ref{table:main}, namely
\begin{equation} \label{eq:sig12rep}
\left[ \begin{matrix}
\alpha & \beta & \cdot \\
\beta & \gamma & \beta \\
\cdot & \beta & \gamma
\end{matrix} \right]. 
\end{equation}
The proof is similar to that of Lemma~\ref{dist13}, and we can avoid repeating certain calculations by noting that the above matrix can be obtained from the one in \eqref{eq:sig12rep} by setting $\varepsilon=1$. 

The point-orbit distribution of $\pi$ is $[1, \frac{q-1}{2}, \frac{q-1}{2}, q^2+1]$. 
The points of rank at most $2$ in $\pi$ lie on the cubic $\mathcal{C} : \alpha (\gamma^2-\beta^2) -\beta^2\gamma = 0$.  
Hence, through the unique point $x : \beta=\gamma=0$ of rank~$1$, there are $\tfrac{q-1}{2}$ lines of each of the types $o_{8,1}$ and $o_{8,2}$, and two lines of type $o_9$. 

The inflexion points of $\mathcal{C}$ are precisely the rank-$2$ points satisfying $\gamma(3\beta^2+\gamma^2)=0$. 
The point $y = (0,1,0)$ is an inflexion point for all $q$. 
It is the only inflexion point unless $q \equiv 1 \pmod 3$, in which case $-3$ is a square and there are two further inflexion points, namely $y_\pm = (-\tfrac{3}{4}, \pm \sqrt{-3}, 3)$.

If $w$ is a point of rank $2$, with coordinates $(\alpha_w,\beta_w,\gamma_w)$, then
\begin{itemize}
\item[(i)] $(\beta_w,\gamma_w) \neq (0,0)$, because $w\neq x$;
\item[(ii)] $\gamma_w^2-\beta_w^2 \neq 0$, by (i);
\item[(iii)] $\alpha_w = \beta_w^2\gamma_w/(\gamma_w^2-\beta_w^2)$, by (ii);
\item[(iv)] $w$ is exterior if and only if $\beta_w^2-\gamma_w^2$ is a non-zero square in $\mathbb{F}_q$, by Lemma~\ref{lem:extcriteria}, and in particular all inflexion points of $\mathcal{C}$ are exterior. 
\end{itemize}
If $w$ is not an inflexion point then the tangent line to $\mathcal{C}$ at $w$ contains exactly one other point of rank $2$, say $v$ with coordinates $(\alpha_v,\beta_v,\gamma_v)$. 
By setting $\varepsilon=1$ in the corresponding argument in the proof of Lemma~\ref{dist12} and applying (iv), we deduce that
\begin{itemize}
\item[(vi)] $v$ is exterior, regardless of whether $w$ is exterior or interior. 
\end{itemize}
By (vi), if $w$ is not an inflexion point then the tangent line to $\mathcal{C}$ at $w$ has type $o_{13,1}$ or $o_{13,2}$ according to whether $w$ is exterior or interior. 
By (iv), all inflexion points are exterior, so there are $\tfrac{q-1}{2}$ lines of type $o_{13,2}$ in $\pi$ for all $q$. 
The number of lines of type $o_{13,1}$ in $\pi$ is $\tfrac{q-1}{2}$ minus the number of inflexion points, namely $\tfrac{q-1}{2}-3 = \tfrac{q-7}{2}$ or $\tfrac{q-1}{2}-1 = \tfrac{q-3}{2}$ according to whether $q \equiv 1 \pmod 3$ or not. 
The tangent line through a point of inflexion contains no other points of $\mathcal{C}$, so has point-orbit distribution $[0,1,0,q-1]$ and hence type $o_{15,1}$ or $o_{16}$. 
By Claim~$\Sigma_{12}$ at the end of the proof, it has type $o_{16}$, yielding a total of three lines of type $o_{16}$ when $q \equiv 1 \pmod 3$ and a unique line of type $o_{16}$ otherwise.

We now count the remaining lines containing a point of rank $2$. 
Suppose first that $w_i$ is an interior rank-$2$ point. 
We have already counted the line $\langle w_i,x \rangle$ (of type $o_{8,2}$) and the tangent line to $\mathcal{C}$ at $w_i$ (of type $o_{13,2}$). 
Let $\ell$ be one of the remaining $q-1$ lines through $w_i$. 
By (vi), $\ell$ is not the tangent to any other point on $\mathcal{C}$, so if it meets $\mathcal{C}$ in a second point then it has type $o_{14,2}$. 
There are $\tfrac{q-1}{2}$ choices of $w_i$ and $\tfrac{q-1}{2}-1$ lines of type $o_{14,2}$ through each $w_i$, so $\pi$ contains $\tfrac{1}{2} \cdot \tfrac{q-1}{2} \cdot \tfrac{q-3}{2}$ such lines in total. 
The remaining $(q-1)-\tfrac{q-3}{2} = \tfrac{q+1}{2}$ lines through $w_i$ have type $o_{15,2}$, so $\pi$ contains a total of $\tfrac{(q+1)(q-1)}{4}$ such lines. 

To count the remaining lines through an exterior point of rank $2$, we first note:
\begin{itemize}
\item[(vii)] An exterior rank-$2$ point $v$ on $\mathcal{C}$ lies on exactly two or three tangent lines to $\mathcal{C}$ (one of which is the tangent line at $v$), according to whether $v$ is an inflexion point or not. 
If $v$ is an inflexion point and lies on the tangent line to $\mathcal{C}$ at the point $w$, say, then $w$ is exterior if and only if $q \equiv 1 \pmod 4$. 
If $v$ is not an inflexion point and lies on the tangent lines at the points $w$ and $w'$, then $w$ and $w'$ have the same type (both exterior or both interior) if and only if $q \equiv 1 \pmod 4$. 
\end{itemize}

{\em Proof of} (vii). 
The asserted fact is analogous to fact (vii) in the proof of Lemma~\ref{dist13}, except that the condition $q \equiv -1 \pmod 4$ there has changed to $q \equiv 1 \pmod 4$. 
In the case where $v$ is not an inflexion point, this change of sign occurs because when setting $\varepsilon=1$ in \eqref{sig13(vii)generic}, we see that the two factors on the left-hand side are both squares if and only if $-1$ is a square (rather than a non-square). 
The case where $v$ is an inflexion can be checked similarly.

Suppose that $q \equiv -1 \pmod 4$. 
Consider an exterior rank-$2$ point $w_e$, and recall that we have already counted the line through $w_e$ and the unique point of rank $1$. 
Suppose first that $w_e$ is an inflexion point. 
By (vii), $w_e$ lies on the tangent line to $\mathcal{C}$ at a unique interior rank-$2$ point $z_i$. 
This tangent line has already been counted. 
The tangent line to $w_e$ at $\mathcal{C}$ (which contains no other rank-$2$ points) has also been counted, so there are $(q+1)-3=q-2$ lines through $w_e$ left to consider. 
Any other line through $w_e$ containing an interior rank-$2$ point has type $o_{14,2}$ and has already been counted. 
There are $\tfrac{q-1}{2}-1=\tfrac{q-3}{2}$ interior rank-$2$ points other than $z_i$, and hence half this many lines of type $o_{14,2}$ through $w_e$, leaving $(q-2)-\tfrac{q-3}{4}=\tfrac{3q-5}{4}$ lines through $w_e$ to consider. 
Any other line through $w_e$ containing a second exterior rank-$2$ point, $w_e'$ say, has type $o_{14,1}$ and contains also a third exterior rank-$2$ point. 
There are $\tfrac{q-1}{2} - 1$ choices for $w_e'$, and half this many, namely $\tfrac{q-3}{4}$, lines of type $o_{14,1}$ through $w_e$. 
The remaining $\tfrac{3q-5}{4}-\tfrac{q-3}{4}=\tfrac{q-1}{2}$ lines through $w_e$ have point-orbit distribution $[0,1,0,q]$ and hence type $o_{15,1}$ or $o_{16}$. 
Claim~$\Sigma_{12}$ below shows that they all have type $o_{15,1}$. 
We count the total number of lines of types $o_{14,1}$ and $o_{15,1}$ after treating also the non-inflexion points. 

Assuming still that $q \equiv -1 \pmod 4$, suppose now that $w_e$ is not an inflexion point. 
By (vii), $w_e$ lies on the tangent lines to $\mathcal{C}$ at one exterior rank-$2$ point $z_e$ and one interior rank-$2$ point $z_i$ (and $z_e$ is never an inflexion point). 
These tangent lines have already been counted. 
By (vi), the tangent line to $w_e$ at $\mathcal{C}$ has type $o_{13,1}$ and has also been counted; denote the second exterior rank-$2$ point on this line by $w_e^T$. 
Any other line through $w_e$ containing an interior rank-$2$ point has type $o_{14,2}$ and has already been counted. 
There are $\tfrac{q-3}{2}$ interior rank-$2$ points other than $z_i$, and half this many lines of type $o_{14,2}$ through $w_e$, leaving $(q+1)-4-\tfrac{q-3}{4}=\tfrac{3q-9}{4}$ lines through $w_e$ to consider. 
Any other line through $w_e$ containing a second exterior rank-$2$ point has type $o_{14,1}$.  
These lines meet a total of $\tfrac{q-1}{2} - 3 = \tfrac{q-7}{2}$ exterior rank-$2$ points, namely those other than $w_e$, $z_e$ and $w_e^T$, so there are $\tfrac{q-7}{4}$ of them through $w_e$. 
The remaining $\tfrac{3q-9}{4} - \tfrac{q-7}{4} = \tfrac{q-1}{2}$ lines through each non-inflexion point $w_e$ have point-orbit distribution $[0,1,0,q]$ and hence type $o_{15,1}$ or $o_{16}$. 
Claim~$\Sigma_{12}$ below implies that they all have type $o_{15,1}$. 

Each exterior rank-$2$ point, whether an inflexion point or not, therefore lies on $\tfrac{q-1}{2}$ lines of type $o_{15,1}$, so there are in total $\tfrac{(q-1)^2}{4}$ lines of type $o_{15,1}$ in $\pi$.
To calculate the total number, call it $N$, of lines of type $o_{14,1}$ in $\pi$, we count in two different ways the number of incident point--line pairs $(w,\ell)$ with $w$ an exterior rank-$2$ point and $\ell$ a line of type $o_{14,1}$. 
On the one hand, there are $3N$ such pairs. 
On the other hand, the number of pairs $(w,\ell)$ is equal to $N_i+N_i'$, where $N_i$ (respectively $N_i'$) is the number of such pairs with $w_e$ an inflexion (respectively non-inflexion) point. 
By the above arguments, $N_i = \tfrac{3(q-3)}{4}$ or $\tfrac{q-3}{4}$ according to whether $q \equiv 1 \pmod 3$ or not, and $N_i' = \tfrac{q-7}{2} \cdot \tfrac{q-7}{4}$ or $\tfrac{q-3}{2} \cdot \tfrac{q-7}{4}$ in these respective cases. 
In each case, $N=\tfrac{1}{3}(N_i+N_i')$.

To complete the proof for $q \equiv -1 \pmod 4$, it remains to note that all lines not counted thus far have type $o_{17}$ (and that $\pi$ contains $q^2+q+1$ lines in total).

Now suppose that $q \equiv 1 \pmod 4$. 
Consider an exterior rank-$2$ point $w_e$, and recall that we have already counted the line through $w_e$ and the unique point of rank $1$. 
Suppose first that $w_e$ is an inflexion point. 
By (vii), $w_e$ lies on the tangent line to $\mathcal{C}$ at a unique exterior rank-$2$ point $z_e$. 
This tangent line has already been counted. 
The tangent line to $w_e$ at $\mathcal{C}$, which contains no other rank-$2$ points, has also been counted, so there are $(q+1)-3=q-2$ lines through $w_e$ left to consider. 
Any line through $w_e$ containing an interior rank-$2$ point has type $o_{14,2}$ and has already been counted. 
There are $\tfrac{q-1}{2}$ interior rank-$2$ points, and hence half this many lines of type $o_{14,2}$ through $w_e$, leaving $(q-2)-\tfrac{q-1}{4}=\tfrac{3q-7}{4}$ lines through $w_e$ to consider. 
Apart from the line $\langle w_e,z_e \rangle$, any other line through $w_e$ containing a second exterior rank-$2$ point, $w_e'$ say, has type $o_{14,1}$ and contains also a third exterior rank-$2$ point. 
There are $\tfrac{q-1}{2} - 2=\tfrac{q-5}{2}$ choices for $w_e'$, and half this many lines of type $o_{14,1}$ through $w_e$. 
The remaining $\tfrac{3q-7}{4}-\tfrac{q-5}{4}=\tfrac{q-1}{2}$ lines through $w_e$ have point-orbit distribution $[0,1,0,q]$ and hence type $o_{15,1}$ or $o_{16}$. 
Claim~$\Sigma_{12}$ implies that they are all of type $o_{15,1}$. 
Again, we delay counting the total number of lines of types $o_{14,1}$ and $o_{15,1}$ until we have treated the case in which $w_e$ is not an inflexion point. 

Now suppose that $w_e$ is not an inflexion point. 
By (vii), there are two possibilities: either $w_e$ lies on the tangent line to $\mathcal{C}$ at some interior rank-$2$ and we say that it is of {\em class~I}, or it does not and we say that it is of {\em class~E}. 
There are $\tfrac{q-1}{4}$ points of class~I, each lying on the tangent lines to $\mathcal{C}$ at two of the $\tfrac{q-1}{2}$ interior rank-$2$ points. 
The remaining non-inflexion points $w_e$ are of class~E; there are $\tfrac{q-1}{4}-3=\tfrac{q-13}{4}$ or $\tfrac{q-1}{4}-1=\tfrac{q-5}{4}$ of them according to whether $q \equiv 1 \pmod 3$ or not. 
In either case, we have already counted the line through $w_e$ and the point of rank $1$. 

First suppose that $w_e$ is of class~E.  
By (vii), $w_e$ lies on the tangent lines to $\mathcal{C}$ at two exterior rank-$2$ points, say $z_1$ and $z_2$. 
These tangent lines have already been counted. 
The tangent line to $\mathcal{C}$ at $w_e$ has also already been counted; it has type $o_{13,1}$ and contains an exterior rank-$2$ point $w_e^T \not \in \{w_e,z_1,z_2\}$. 
There are $q-3$ lines through $w_e$ left to consider. 
Let $\ell$ be a line through $w_e$ containing an interior rank-$2$ point. 
Since $\ell$ is not the tangent line to $\mathcal{C}$ at $w_e$, it also contains a second interior rank-$2$ point, and has type $o_{14,2}$. 
There are $\tfrac{q-1}{2}$ interior rank-$2$ points, and hence half this many lines of type $o_{14,2}$ through $w_e$, leaving $(q-3) - \tfrac{q-1}{4} = \tfrac{3q-11}{4}$ lines through $w_e$ to consider. 
Now let $\ell$ be a line through $w_e$ containing an exterior rank-$2$ point $w_e' \not \in \{w_e^T,z_1,z_2\}$. 
Any such line has type $o_{14,1}$. 
There are $\tfrac{q-1}{2}-4=\tfrac{q-9}{2}$ exterior rank-$2$ points other than $w_e$, $w_e^T$, $z_1$ and $z_2$, and therefore $\tfrac{q-9}{4}$ lines of type $o_{14,1}$ through $w_e$. 
The remaining $\tfrac{q-1}{2}$ lines through $w_e$ have point-orbit distribution $[0,1,0,q]$ and, by Claim~$\Sigma_{12}$, type $o_{15,1}$. 

Finally, suppose that $w_e$ is of class~I.  
By (vii), $w_e$ lies on the tangent lines to $\mathcal{C}$ at two interior rank-$2$ points, say $z_1$ and $z_2$. 
These tangent lines have already been counted. 
The tangent line to $\mathcal{C}$ at $w_e$ has also already been counted; it has type $o_{13,1}$ and contains an exterior rank-$2$ point $w_e^T \neq w_e$. 
There are $q-3$ lines through $w_e$ left to consider. 
Let $\ell$ be a line through $w_e$ containing an interior rank-$2$ point. 
Since $\ell$ is not the tangent line to $\mathcal{C}$ at $w_e$, it also contains a second interior rank-$2$ point, and has type $o_{14,2}$. 
There are $\tfrac{q-1}{2}-2=\tfrac{q-5}{2}$ interior rank-$2$ points other than $z_1$ and $z_2$, and hence half this many lines of type $o_{14,2}$ through $w_e$, leaving $(q-3) - \tfrac{q-5}{4} = \tfrac{3q-7}{4}$ lines through $w_e$ to consider. 
Now let $\ell$ be a line through $w_e$ containing an exterior rank-$2$ point $w_e' \neq w_e^T$. 
Any such line has type $o_{14,1}$. 
There are $\tfrac{q-1}{2}-2=\tfrac{q-5}{2}$ exterior rank-$2$ points other than $w_e$ and $w_e^T$, and hence $\tfrac{q-5}{4}$ lines of type $o_{14,1}$ through $w_e$. 
By Claim~$\Sigma_{12}$, the remaining $\tfrac{q-1}{2}$ lines through $w_e$ have type $o_{15,1}$.

Let us finally count the total number of lines of types $o_{14,1}$ and $o_{15,1}$ in $\pi$ in the case $q \equiv 1 \pmod 4$. 
We have argued, in particular, that every exterior rank-$2$ point lies on $\tfrac{q-1}{2}$ lines of type $o_{15,1}$. 
Hence, $\pi$ contains $\tfrac{(q-1)^2}{4}$ such lines in total. 
Let $N$ denote the number of lines of type $o_{14,1}$ in $\pi$. 
We count in two different ways the incident point--line pairs $(w,\ell)$ with $w$ an exterior rank-$2$ point and $\ell$ a line of type $o_{14,1}$. 
This number equals $3N$, so we have $N = \tfrac{1}{3} (N_i + N_i' + N_i'')$, where $N_i$ is the number of such pairs with $w$ an inflexion point, $N_i'$ is the number of such pairs with $w$ a non-inflexion point of class~E, and $N_i''$ is the number of such pairs with $w$ a non-inflexion point of class~I. 
Each inflexion point lies on $\tfrac{q-5}{4}$ lines of type $o_{14,1}$, so $N_i = \tfrac{3(q-5)}{4}$ or $\tfrac{q-5}{4}$ according to whether $q \equiv 1 \pmod 3$ or not. 
There are $\tfrac{q-13}{4}$ or $\tfrac{q-5}{4}$ points of class~E in these respective cases, and each such point lies on $\tfrac{q-9}{4}$ lines of type $o_{14,1}$, so $N_i' = \tfrac{q-13}{4} \cdot \tfrac{q-9}{4}$ or $\tfrac{q-5}{4} \cdot \tfrac{q-9}{4}$. 
Finally, there are $\tfrac{q-1}{4}$ points of class~I, each lying on $\tfrac{q-5}{4}$ lines of type $o_{14,1}$, so $N_i'' = \tfrac{q-1}{4} \cdot \tfrac{q-5}{4}$ regardless of whether $q \equiv 1 \pmod 3$ or not. 
The asserted values of $N$ are now readily calculated. 

To complete the proof in the case $q \equiv 1 \pmod 4$, observe that all lines not counted thus far have type $o_{17}$.
It remains to prove the following:

{\em Claim~$\Sigma_{12}$.} Let $\ell$ be a line in $\pi$ with point-orbit distribution $[0,1,0,q]$, and let $w$ denote the unique point of rank $2$ on $\ell$. 
Then $\ell$ has type $o_{15,1}$ unless $w$ is an inflexion point and $\ell$ is the tangent line to $\mathcal{C}$ at $w$, in which case $\ell$ has type $o_{16}$.

{\em Proof of Claim~$\Sigma_{12}$.} 
Let $(\alpha_w,\beta_w,\gamma_w)$ be the coordinates of $w$, and recall that $\gamma_w^2-\beta_w^2 \neq 0$ and $\alpha_w = \beta_w^2\gamma_w/(\gamma_w^2-\beta_w^2)$. 
We argue as in the proof of Claim~$\Sigma_{13}$ in Lemma~\ref{dist13}. 
If $\beta_w \neq 0$, we transform the matrix representative $M_w$ of $w$ obtained from \eqref{eq:sig12rep} by setting $\varepsilon=1$ in the matrix $Y_w$ in the proof of Claim~$\Sigma_{13}$. 
The gives the representative of $U:=\langle W,\ell \rangle$ shown there, with $\varepsilon=1$, and with the factor $f_w(\alpha,\beta,\gamma)$ in the $(3,3)$-entry now equal to the left-hand side of \eqref{eq:sig13tan} with $\varepsilon=1$. 
The result follows because the tangent line to $\mathcal{C}$ at $w$ is given by \eqref{eq:sig13tan} with $\varepsilon=1$. 
The argument for $\beta_w=0$ is the same as in Claim~$\Sigma_{13}$ (with $\varepsilon=1$). 
\end{Prf}

\section{Plane stabilisers} \label{s:stabilisers}

We now calculate the stabilisers in $K$ of representatives of each of the plane orbits, and thereby determine the size of each orbit (using the orbit--stabiliser theorem and the fact that $|K| = |\textnormal{PGL}(3,q)| = q^3(q^3-1)(q^2-1)$). 
Although there are various ways to go about this, we take an elementary approach, calculating the stabilisers directly from the orbit representatives. 
The line-orbit distributions determined in the previous section (along with the point-orbit distributions) are useful in guiding these calculations. 
For example, for all orbits $\Sigma_i$ with $i \in \{3,4,5,6,8,9,10,11,15\}$ (and for $\Sigma_{14}'$), a plane $\pi \in \Sigma_i$ contains a unique line $\ell$ of a certain type, so we immediately know that $K_\pi$ must be a subgroup of $K_\ell$. 
(As noted earlier, the line stabilisers have been determined previously, in \cite{LaPo2020}.)

Before proceeding, let us establish some group-theoretic notation (which has already been used in Tables~\ref{TableStabilisers} and \ref{table:line-rank-dist}):
\begin{itemize}
\item $E_q$ denotes an elementary abelian group of order $q$. 
\item $C_k$ denotes a cyclic group of order $k$. 
\item $S_n$ denotes the symmetric group on $n$ letters. 
\item $A \times B$ is the direct product of groups $A$ and $B$; and $A:B$ is a semidirect product with normal subgroup $A$ and subgroup $B$. 
\item $\GO^\pm(2,q)$ denotes the similarity group of a non-degenerate quadratic form of $\pm$ type on a $2$-dimensional vector space over $\mathbb{F}_q$, which has order $2(q-1)(q \mp 1)$. (These groups appear only in the stabilisers of planes of type $\Sigma_6$ and lines of type $o_{8,1}$, $o_{8,2}$ and $o_{10}$.)
\item $E_q^{1+2}$ denotes a group with centre $Z \cong E_q$ and $E_q^{1+2}/Z \cong E_q^2$. 
(This group appears only in the stabilisers of planes of type $\Sigma_{15}$ and lines of type $o_6$.)
\end{itemize}

\begin{La} \label{stab1}
The stabiliser of a plane of type $\Sigma_1$ is isomorphic to $E_q^2:\textnormal{GL}(2,q)$. 
The number of planes of type $\Sigma_1$ is $q^2+q+1$.
\end{La}

\begin{Prf}
Let $\pi$ denote the representative of $\Sigma_1$ given in Table~\ref{table:main}. 
That is, suppose that $\pi$ is represented by the matrix
\[
M = \left[ \begin{matrix}
\alpha & \gamma & \cdot \\
\gamma & \beta & \cdot \\
\cdot & \cdot & \cdot
\end{matrix} \right], 
\]
where, as usual, $(\alpha,\beta,\gamma)$ ranges over $\mathbb{F}_q^3 \setminus \{ (0,0,0) \}$. 
A matrix $g=(g_{ij}) \in \GL(3,q)$ fixes $M$ with respect to the action $M \mapsto gMg^T$ if and only if $g_{31}=g_{32}=0$. 
Upon factoring out scalars, we therefore obtain $K_\pi \cong E_q^2 : \GL(2,q)$. 
Explicitly, if $G$ denotes the group of all $g \le \GL(3,q)$ with $g_{31}=g_{32}=0$, then $K_\pi = G/Z$, where $Z \le \GL(3,q)$ is the group of all scalar matrices. 
If we define $H$ to be the subgroup of all $g \in G$ with $g_{33}=1$ then $G=HZ$, and so the second isomorphism theorem yields $K_\pi = G/Z \cong H/(H\cap Z) = H \cong E_q^2 : \GL(2,q)$. 
\end{Prf}

\begin{La} \label{stab2}
The stabiliser of a plane of type $\Sigma_2$ is isomorphic to $C_{q-1}^2:S_3$. 
The number of planes of type $\Sigma_2$ is $\tfrac{1}{6}q^3(q^2+q+1)(q+1)$.
\end{La}

\begin{Prf}
Let $\pi$ denote the representative of $\Sigma_2$ given in Table~\ref{table:main}. 
That is, suppose that $\pi$ is represented by the matrix
\[
M = \left[ \begin{matrix}
\alpha & \cdot & \cdot \\
\cdot & \beta & \cdot \\
\cdot & \cdot & \gamma
\end{matrix} \right]. 
\]
The matrices $g \in \GL(3,q)$ fixing $M$ with respect to the action $M \mapsto gMg^T$ comprise the subgroup of $\GL(3,q)$ generated by the diagonal matrices and the matrices obtained from the identity by permuting the columns. 
Upon factoring out scalars, this gives $K_\pi \cong C_{q-1}^2 : S_3$. 
(Explicitly, one may argue as in the proof of Lemma~\ref{stab1}: we have $K_\pi = G/Z$ where $G \le \GL(3,q)$ is the subgroup of all matrices $g$ described above, and $G=HZ$ where $H$ is the subgroup of all $g \in G$ with $g_{33}=1$, so $K_\pi = G/Z \cong H/(H\cap Z) = H \cong C_{q-1}^2:S_3$.)
\end{Prf}

\begin{La} \label{stab3}
The stabiliser of a plane of type $\Sigma_3$ is isomorphic to $E_q:C_{q-1}^2$. 
The number of planes of type $\Sigma_3$ is $q^2(q^2+q+1)(q+1)$.
\end{La}

\begin{Prf}
Let $\pi$ denote the representative of $\Sigma_3$ given in Table~\ref{table:main}. 
That is, suppose that $\pi$ is represented by the matrix
\[
M = \left[ \begin{matrix}
\alpha & \cdot & \gamma \\
\cdot & \beta & \cdot \\
\gamma & \cdot & \cdot
\end{matrix} \right]. 
\]
By Lemma~\ref{dist3}, $\pi$ contains a unique line of type $o_5$ and a unique line of type $o_6$, namely the lines $\ell : \gamma=0$ and $\ell' : \beta=0$, respectively. 
Therefore, $K_\pi = K_\ell \cap K_{\ell'}$. 
A matrix $g = (g_{ij}) \in \GL(3,q)$ fixes the matrix representative of $\ell'$ with respect to the action $M \mapsto gMg^T$ if and only if 
\[
g = \left[ \begin{matrix}
g_{11} & g_{12} & g_{13} \\
\cdot & g_{22} & \cdot \\
\cdot & g_{32} & g_{33}
\end{matrix} \right],
\]
and a matrix of this form fixes the matrix representative of $\ell$ if and only if $g_{12}=g_{32}=0$. 
The subgroup $G \le \GL(3,q)$ of all such matrices has shape $E_q : C_{q-1}^3$, and factoring out scalars yields $K_\pi \cong E_q : C_{q-1}^2$. 
\end{Prf}

\begin{La} \label{stab4}
The stabiliser of a plane of type $\Sigma_4$ is isomorphic to $(E_q:C_{q-1}):C_2$. 
The number of planes of type $\Sigma_4$ is $\tfrac{1}{2}q^2(q^3-1)(q+1)$.
\end{La}

\begin{Prf}
Let $\pi$ denote the representative of $\Sigma_4$ given in Table~\ref{table:main}. 
That is, suppose that $\pi$ is represented by the matrix
\[
M = \left[ \begin{matrix}
\alpha & \cdot & \gamma \\
\cdot & \beta & \gamma \\
\gamma & \gamma & \cdot
\end{matrix} \right]. 
\]
By Lemma~\ref{dist4}, $\pi$ contains a unique line of type $o_5$ and a unique line of type $o_{12}$, namely the lines $\ell : \gamma=0$ and $\ell' : \alpha+\beta=0$, respectively. 
Therefore, $K_\pi = K_\ell \cap K_{\ell'}$. 
Since $\ell$ contains precisely two points of rank $1$, namely the points $x : \beta=\gamma=0$ and $y : \alpha=\gamma=0$, $K_\ell$ is equal to the stabiliser in $K$ of the set $\{x,y\}$, and so $K_\pi$ is equal to the stabiliser of $\{x,y\}$ and $\ell'$. 
An element of $K$ represented by $g = (g_{ij}) \in \GL(3,q)$ fixes $\{x,y\}$ if and only if it has the form 
\begin{equation} \label{stab4matrices}
\left[ \begin{matrix}
g_{11} & \cdot & g_{13} \\
\cdot & g_{22} & g_{23} \\
\cdot & \cdot & g_{33}
\end{matrix} \right] 
\quad \text{or} \quad 
\left[ \begin{matrix}
\cdot & g_{11} & g_{13} \\
g_{22} & \cdot & g_{23} \\
\cdot & \cdot & g_{33}
\end{matrix} \right], 
\end{equation}
where the second matrix is obtained from the first by composing with the matrix obtained by swapping the first two columns of the identity. 
Such an element also fixes $\ell'$ if and only if $g_{11}=g_{22}$ and $g_{13} = -g_{23}$. 
Factoring out scalars yields $K_\pi \cong (E_q : C_{q-1}) : C_2$. 
\end{Prf}

\begin{La} \label{stab5}
The stabiliser of a plane of type $\Sigma_5$ is isomorphic to $C_{q-1}:C_2$. 
The number of planes of type $\Sigma_5$ is $\tfrac{1}{2}q^3(q^3-1)(q+1)$.
\end{La}

\begin{Prf}
Let $\pi$ denote the representative of $\Sigma_5$ given in Table~\ref{table:main}. 
That is, suppose that $\pi$ is represented by the matrix
\[
M = \left[ \begin{matrix}
\alpha & \cdot & \gamma \\
\cdot & \beta & \gamma \\
\gamma & \gamma & \gamma
\end{matrix} \right]. 
\]
By Lemma~\ref{dist5}, $\pi$ contains exactly two lines $\ell$ and $\ell'$ of type $o_9$, so $K_\pi$ is equal to the stabiliser in $K$ of the set $\{ \ell,\ell' \}$. 
Writing $x : \beta=\gamma=0$ and $y : \alpha=\gamma=0$ for the rank-$1$ points in $\pi$, we have  $\ell = \langle x,w \rangle$ and $\ell' = \langle y,w \rangle$ where $w$ is the rank-$3$ point with coordinates $(\alpha,\beta,\gamma)=(1,1,1)$. 
Now, $h \in K$ fixes $\{ \ell,\ell' \}$ if and only if (i) it fixes $\{ x,y \}$, and (ii) it fixes $w$. 
Condition (i) is equivalent to $h$ fixing the line of type $o_5$ spanned by $x$ and $y$, and so the same  calculation as in Lemma~\ref{stab4} shows that a matrix representative $g = (g_{ij}) \in \GL(3,q)$ of $h$ must be of the form \eqref{stab4matrices}. 
Such an element additionally satisfies condition (ii) if and only if $g_{11} = g_{33}^2 g_{22}^{-1}$, $g_{13} = g_{33}-g_{11}$ and $g_{23} = g_{33}-g_{22}$. 
Factoring out scalars yields $K_\pi \cong C_{q-1} : C_2$.
\end{Prf}

\begin{La} \label{stab6}
The stabiliser of a plane of type $\Sigma_6$ is isomorphic to $\GO^-(2,q)$, namely, the group of similarities of a non-degenerate quadratic form of $-$ type on a $2$-dimensional vector space over $\mathbb{F}_q$. 
The number of planes of type $\Sigma_6$ is $\tfrac{1}{2}q^3(q^3-1)$.
\end{La}

\begin{Prf}
Let $\pi$ denote the representative of $\Sigma_6$ given in Table~\ref{table:main}. 
That is, suppose that $\pi$ is represented by the matrix
\[
M = \left[ \begin{matrix} \alpha & \beta & \cdot \\ 
\beta & \varepsilon\alpha & \cdot \\ 
\cdot & \cdot & \gamma \end{matrix} \right], 
\quad \text{where } \varepsilon \text{ is a non-square in } \mathbb{F}_q.
\]
By Lemma~\ref{dist6}, $\pi$ contains a unique line of type $o_{10}$, namely the line $\ell : \gamma=0$. 
Since the unique point $x : \alpha=\beta=0$ of rank $1$ in $\pi$ does not lie on $\ell$, we have $K_\pi = K_\ell \cap K_x$. 
The matrix representatives of elements of $K_x$ are precisely the matrices $g = (g_{ij}) \in \GL(3,q)$ with $g_{13}=g_{23}=0$. 
Such a matrix $g$ also fixes the matrix representative of the line $\gamma=0$ if and only if $g_{31}=g_{32}=0$ and the entries $g_{ij}$ with $i,j \le 2$ satisfy
\[
\varepsilon g_{11}x_{12} = g_{21}x_{22} 
\quad \text{and} \quad
\varepsilon (g_{11}^2 - g_{22}^2) = g_{21}^2 - \varepsilon^2 g_{12}^2.
\]
This is the case if and only if
\[
g = \left[ \begin{matrix}
g_{11} & g_{12} & \cdot \\
\pm \varepsilon g_{12} & \pm g_{11} & \cdot \\
\cdot & \cdot & g_{33}
\end{matrix} \right]. 
\]
The upper-left $2 \times 2$ sub-matrices of $g$ comprise an orthogonal group $\GO^-(2,q)$, specifically the group of similarities of the quadratic form associated with the matrix
\[
\left[ \begin{matrix}
1 & \cdot \\
\cdot & -\varepsilon 
\end{matrix} \right], 
\]
which spans the complement of the subspace of the symmetric $2 \times 2$ matrices defined by the upper-left $2 \times 2$ block of $M$. 
Hence, upon factoring out scalars, we obtain $K_\pi \cong \GO^-(2,q)$. 
To calculate the orbit size, recall that $|\GO^-(2,q)| = 2(q^2-1)$. 
\end{Prf}

\begin{La} \label{stab7}
The stabiliser of a plane of type $\Sigma_7$ is isomorphic to $E_q^2:\textnormal{GL}(2,q)$. 
The number of planes of type $\Sigma_7$ is $q^2+q+1$.
\end{La}

\begin{Prf}
As noted in the proof of \cite[Lemma~7.3]{LaPoSh2020}, a plane $\pi \in \Sigma_7$ is the tangent plane to $\mathcal{V}(\mathbb{F}_q)$ at its unique point of $x$ rank $1$. 
Hence, the stabiliser of $\pi$ in $K$ is equal to the stabiliser of $x$ in $K$, which is has shape $E_q^2 : \GL(2,q)$. 
Explicitly, if $\pi$ is the representative of $\Sigma_7$ from Table~\ref{table:main} then $x$ is represented by a matrix whose only non-zero entry is the $(1,1)$-entry, so the elements of $K_\pi = K_x$ correspond to matrices $g = (g_{ij}) \in \GL(3,q)$ with $g_{21}=g_{31}=0$.
\end{Prf}

\begin{La} \label{stab8}
The stabiliser of a plane of type $\Sigma_8$ is isomorphic to $E_q:C_{q-1}^2$. 
The number of planes of type $\Sigma_8$ is $q^2(q^2+q+1)(q+1)$.
\end{La}

\begin{Prf}
Let $\pi$ denote the representative of $\Sigma_8$ given in Table~\ref{table:main}. 
That is, suppose that $\pi$ is represented by the matrix
\[
M = \left[ \begin{matrix} 
\alpha & \beta & \cdot \\ 
\beta & \cdot & \gamma \\ 
\cdot & \gamma & \cdot 
\end{matrix} \right].
\]
By Lemma~\ref{dist6}, $\pi$ contains a unique line of type $o_{12}$, namely the line $\ell : \alpha=0$. 
Since the unique point $x : \beta=\gamma=0$ of rank $1$ in $\pi$ does not lie on $\ell$, we have $K_\pi = K_\ell \cap K_x$. 
As noted in \cite[Section~4]{LaPo2020}, the matrix representatives of elements of $K_\ell$ are precisely the matrices $g = (g_{ij}) \in \GL(3,q)$ with $g_{12}=g_{21}=g_{23}=g_{32}=0$. 
As noted in the proof of Lemma~\ref{stab7}, the matrix representatives of elements of $K_x$ are precisely the matrices $g = (g_{ij}) \in \GL(3,q)$ with $g_{21}=g_{31}=0$. 
Hence, the matrix representatives of $K_\pi = K_\ell \cap K_x$ have the form
\[
\left[ \begin{matrix} 
g_{11} & \cdot & g_{13} \\ 
\cdot & g_{22} & \cdot \\ 
\cdot & \cdot & g_{33} 
\end{matrix} \right].
\]
Factoring out scalars yields $K_\pi \cong E_q : C_{q-1}^2$. 
(Indeed, note that $K_\pi$ is equal to the stabiliser in $K$ of the plane of type $\Sigma_3$ considered in the proof of Lemma~\ref{stab3}.)
\end{Prf}

\begin{La} \label{stab9}
The stabiliser of a plane of type $\Sigma_9$ is isomorphic to $(E_q:C_{q-1}):C_2$.
The number of planes of type $\Sigma_9$ is $\tfrac{1}{2}q^2(q^3-1)(q+1)$.
\end{La}

\begin{Prf}
Let $\pi$ denote the representative of $\Sigma_9$ given in Table~\ref{table:main}. 
That is, suppose that $\pi$ is represented by the matrix
\begin{equation} \label{sigma9stabrep}
M = \left[ \begin{matrix} 
\alpha & \beta & \cdot \\ 
\beta & \gamma & \cdot \\ 
\cdot & \cdot & -\gamma 
\end{matrix} \right].
\end{equation}
As noted in the proof of Lemma~\ref{dist9}, the points of rank at most $2$ in $\pi$ comprise the conic $\mathcal{C} : \beta^2-\alpha\gamma = 0$, which intersects the unique line $\ell : \gamma=0$ of type $o_6$ in $\pi$ in the unique point $x : \beta=\gamma=0$ of rank $1$ in $\pi$. 
Hence, we have $K_\pi = K_\ell \cap K_{\mathcal{C}'}$, where $\mathcal{C}' = \mathcal{C} \setminus \{x\}$ (and $K_{\mathcal{C}'}$ is its setwise stabiliser). 
The matrix representatives of elements of $K_\ell$ are precisely the matrices $g = (g_{ij}) \in \GL(3,q)$ with $g_{21}=g_{31}=g_{32}=0$. 
A matrix $g$ of this form also fixes the matrix representative of $\mathcal{C}'$ (obtained by setting $\alpha = \frac{\beta^2}{\gamma}$ in $M$) if and only if $g_{12}=g_{32}=0$ and $g_{33}^2 = g_{22}^2$. 
That is, the matrix representatives of $K_\pi$ are precisely the matrices of the form
\[
\left[ \begin{matrix} 
g_{11} & g_{12} & \cdot \\ 
\cdot & g_{22} & \cdot \\ 
\cdot & \cdot & \pm g_{22} 
\end{matrix} \right].
\]
Factoring out scalars yields $K_\pi \cong (E_q : C_{q-1}) : C_2$. 
\end{Prf}

\begin{La} \label{stab10}
The stabiliser of a plane of type $\Sigma_{10}$ is isomorphic to $(E_q:C_{q-1}):C_2$.
The number of planes of type $\Sigma_{10}$ is $\tfrac{1}{2}q^2(q^3-1)(q+1)$.
\end{La}

\begin{Prf}
The proof is essentially identical to that of Lemma~\ref{stab9}. 
To use the representative $\pi \in \Sigma_{10}$ from Table~\ref{table:main}, we just replace $\gamma$ in \eqref{sigma9stabrep} by $-\varepsilon\gamma$ with $\varepsilon \in \mathbb{F}_q$ a non-square. 
By Lemma~\ref{dist10} (and its proof), the points of rank at most $2$ in $\pi$ again comprise the conic $\mathcal{C} : \beta^2-\alpha\gamma = 0$, which again intersects the unique line of type $o_6$ in the unique point of rank $1$. 
The subsequent calculations from the proof of Lemma~\ref{stab9} go through verbatim. 
\end{Prf}

\begin{La} \label{stab11}
The stabiliser of a plane of type $\Sigma_{11}$ is isomorphic to $C_{q-1}$ if $q \not \equiv 0 \pmod 3$ and to $E_q$ if $q \equiv 0 \pmod 3$. 
The number of planes of type $\Sigma_{11}$ is $q^3(q^3-1)(q+1)$ and $q^2(q^3-1)(q^2-1)$ in these respective cases.
\end{La}

\begin{Prf}
Let $\pi$ denote the representative of $\Sigma_{11}$ used in the proof of Lemma~\ref{dist11} (as opposed to the one given in Table~\ref{table:main}). 
That is, suppose that $\pi$ is represented by the matrix
\[
M = \left[ \begin{matrix} 
\alpha & \beta & \cdot \\ 
\beta & \cdot & \gamma-\beta \\ 
\cdot & \gamma-\beta & \gamma 
\end{matrix} \right].
\]
By Lemma~\ref{dist11}, $\pi$ contains a unique line of type $o_9$, specifically the line $\ell : \beta=\gamma$. 
Hence, we have $K_\pi \le K_\ell$. 
The matrix representatives of elements of $K_\ell$ are precisely the matrices $g = (g_{ij}) \in \GL(3,q)$ of the form 
\begin{equation} \label{sigma11o9stab}
\left[ \begin{matrix} 
g_{11} & g_{12} & g_{13} \\ 
\cdot & g_{33}^2g_{11}^{-1} & \cdot \\ 
\cdot & -g_{13}g_{33}g_{11}^{-1} & g_{33} 
\end{matrix} \right].
\end{equation}
We must now consider separately the cases where $q \equiv 0 \pmod 3$ or not. 
Recall from the proof of Lemma~\ref{dist11} that the points of rank at most $2$ in $\pi$ comprise the cubic $\mathcal{C} : \beta^2\gamma + \alpha(\beta-\gamma)^2 = 0$. 

First suppose that $q \not \equiv 0 \pmod 3$. 
Then $\mathcal{C}$ contains a unique point of inflexion, namely the point $z$ with coordinates $(\alpha,\beta,\gamma)=(\frac{2}{9},1,-2)$, and $\pi$ contains a unique line of type $o_{16}$, namely the tangent to $\mathcal{C}$ through $z$, which is the line $\ell' : \gamma = 27\alpha+8\beta$. 
Hence, we have $K_\pi = K_\ell \cap K_{\ell'}$. 
An element of $K_\ell$, which is necessarily represented by a matrix $g \in \GL(3,q)$ as in \eqref{sigma11o9stab}, fixes $\ell'$ if and only if it fixes $z$, which is the case precisely when
\[
g_{13} = 3^{-1}(g_{11}-g_{33}) 
\quad \text{and} \quad 
g_{12} = \tfrac{1}{9}g_{11}^{-1}(g_{33}-g_{33}g_{11}-2g_{11}^2).
\]
The matrices of the form \eqref{sigma11o9stab} with these additional constraints comprise a subgroup of $\GL(3,q)$ of shape $C_{q-1}^2$, and by factoring out scalars (as in previous cases) we therefore obtain $K_\pi \cong C_{q-1}$. 

Now suppose that $q \equiv 0 \pmod 3$. 
Then the cubic $\mathcal{C}$ meets the line $\ell$ in the unique point $x : \beta=\gamma=0$ of rank $1$ in $\pi$, so we have $K_\pi = K_\ell \cap K_{\mathcal{C}'}$, where $\mathcal{C}' = \mathcal{C} \setminus \{x\}$. 
A matrix $g \in \GL(3,q)$ of the form \eqref{sigma11o9stab} fixes the matrix representative of $\mathcal{C}'$ if and only if 
\[
g_{11} = g_{33}
\quad \text{and} \quad 
g_{12} = g_{33}^{-1}g_{13}(g_{13}-g_{33}).
\]
Upon factoring out scalars we obtain $K_\pi \cong E_q$. 
\end{Prf}

\begin{La} \label{stab12+13}
The stabiliser of a plane of type $\Sigma_{12}$ is isomorphic to $S_3$ or $C_2$ according to whether $q \equiv 1 \pmod 3$ or not. 
The number of planes of type $\Sigma_{12}$ is $\tfrac{1}{6} q^3(q^3-1)(q^2-1)$ and $\tfrac{1}{2} q^3(q^3-1)(q^2-1)$ in these respective cases. 
The stabiliser of a plane of type $\Sigma_{13}$ is isomorphic to $S_3$ or $C_2$ according to whether $q \equiv -1 \pmod 3$ or not. 
The number of planes of type $\Sigma_{13}$ is $\tfrac{1}{6} q^3(q^3-1)(q^2-1)$ and $\tfrac{1}{2} q^3(q^3-1)(q^2-1)$ in these respective cases.
\end{La}

\begin{Prf}
Let $\pi$ denote the plane represented by the matrix
\[
M = \left[ \begin{matrix}
\alpha & \beta & \cdot \\
\beta & \gamma & \beta \\
\cdot & \beta & \varepsilon\gamma
\end{matrix} \right], 
\quad \text{for some } \varepsilon \in \mathbb{F}_q \setminus \{0\},
\]
so that $\pi$ is the representative of $\Sigma_{12}$, respectively $\Sigma_{13}$, given in Table~\ref{table:main} if $\varepsilon$ is equal to $1$, respectively a non-square. 
Recall from the proofs of Lemmas~\ref{dist13} and \ref{dist12} that the points of rank at most $2$ in $\pi$ lie on the cubic $\mathcal{C} : \alpha (\varepsilon\gamma^2-\beta^2) - \varepsilon\beta^2\gamma = 0$. 
There is a unique point of rank~$1$, namely the point $x : \beta=\gamma=0$. 
The points of inflexion of $\mathcal{C}$ are given by $\gamma(3\beta^2 + \varepsilon\gamma^2)=0$; in particular, the point $y :\alpha=\gamma=0$ is an inflexion point for all $\varepsilon$. 

Write $-3\varepsilon = \delta^2$, so that $\delta$ lies in either $\mathbb{F}_q$ or $\mathbb{F}_{q^2}$, and consider the plane $\overline{\pi}$ obtained by extending the scalars of $M$ to $\mathbb{F}_q(\delta)$. 
Let $\overline{K} \cong \PGL(3,\mathbb{F}_q(\delta))$ denote the group obtained by extending the scalars of $K$ to $\mathbb{F}_q(\delta)$.
The cubic of points of rank at most $2$ in $\overline{\pi}$ contains three points of inflexion over $\mathbb{F}_q(\delta)$, which lie on the line $\overline{\ell} : 4\alpha+\varepsilon\gamma=0$. 
The inflexion points are all exterior rank-$2$ points (see assertion (iv) in the proofs of Lemmas~\ref{dist13} and \ref{dist12}), so $\overline{\ell}$ is of type $o_{14,1}$.
We have $K_\pi \leqslant \overline{K}_{\overline{\ell}} \cap K_x$, with equality when $\delta \in \mathbb{F}_q$. 
If $\delta \not \in \mathbb{F}_q$, we have $K_\pi \leqslant \overline{K}_{\overline{\ell}} \cap K_x \cap K_y$ because $y$ is the unique inflexion point of $\mathcal{C}$. 
(In this case, $\overline{\ell} \cap \pi$ is a line of type $o_{15,1}$.)

For convenience, let us now apply the transformation $M \mapsto XMX^T$, where 
\[
X = \left[ \begin{matrix}
4 & \delta & -1 \\
-4 & \cdot & 4 \\
4 & -\delta & -1
\end{matrix} \right] \in \GL(3,\mathbb{F}_q(\delta)). 
\]
Under this transformation,  $\overline{\ell}$ is mapped to the line $\overline{\ell}'$ represented by the matrix
\[
\left[ \begin{matrix}
a & \cdot & \cdot \\
\cdot & -(a+b) & \cdot \\
\cdot & \cdot & b
\end{matrix} \right], 
\quad \text{where } (a,b) \text{ ranges over all non-zero values in } \mathbb{F}_q(\delta)^2,
\]
and $y,x$ are mapped to the points $y',x'$ represented by the matrices $\text{diag}(1,0,-1)$ and $v \otimes v$ with $v = (1,-1,1)$, respectively. 
We denote the image of $\overline{\pi}$ under this transformation by $\overline{\pi}'$. 

By \cite[Section~4]{LaPo2020} (cf. Table~\ref{table:line-rank-dist}), $\overline{K}_{\overline{\ell}'} \leqslant \PGL(3,\mathbb{F}_q(\delta))$ is a group of shape $C_2^2 : S_3$, generated by elements corresponding to (i) diagonal matrices with all non-zero entries equal to $\pm 1$, which fix pointwise the rank-$2$ points given by $(a,b)=(1,-1)$, $(1,0)$ and $(0,1)$; and (ii) permutation matrices, which permute these points. 
The intersection $\overline{K}_{\overline{\ell}'} \cap K_{x'}$ is isomorphic to $S_3$, generated by elements corresponding to the matrices
\[
A = \left[ \begin{matrix}
\cdot & \cdot & 1 \\
\cdot & 1 & \cdot \\
1 & \cdot & \cdot
\end{matrix} \right] 
\quad \text{and} \quad 
B = \left[ \begin{matrix}
\cdot & \cdot & 1 \\
-1 & \cdot & \cdot \\
\cdot & -1 & \cdot
\end{matrix} \right], 
\]
and $\overline{K}_{\overline{\ell}'} \cap K_{x'} \cap K_{y'}$ is isomorphic to $C_2$, generated by the element corresponding to $A$. 

We can now explicitly calculate the elements of $K_\pi$ by noting that $K_{\overline{\pi}} = X^{-1} K_{\overline{\pi}'} X$. 
If $\delta \in \mathbb{F}_q$ then $K_\pi$ is the image of $X^{-1} \langle A,B \rangle X$ in $\PGL(3,q)$, so $K_\pi \cong S_3$. 
If $\delta \not \in \mathbb{F}_q$ then $K_\pi$ is contained in the image of $X^{-1} \langle A \rangle X$ in $\PGL(3,q)$, which is isomorphic to $C_2$. 
Since the element of $\PGL(3,q)$ corresponding to $X^{-1}AX = \text{diag}(1,-1,1)$ 
fixes $\pi$, we have $K_\pi \cong C_2$. 
The result now follows by taking $\varepsilon=1$ for $\Sigma_{12}$, in which case $\delta \in \mathbb{F}_q$ if and only if $q \equiv 1 \pmod 3$, and $\varepsilon$ a non-square in $\mathbb{F}_q$ for $\Sigma_{13}$, in which case $\delta \in \mathbb{F}_q$ if and only if $q \equiv -1 \pmod 3$.
\end{Prf}

\begin{La} \label{stab14}
Suppose that $q \not \equiv 0 \pmod 3$. 
The stabiliser of a plane of type $\Sigma_{14}$ is isomorphic to $C_3$. 
The number of planes of type $\Sigma_{14}$ is $\tfrac{1}{3} q^3(q^3-1)(q^2-1)$.
\end{La}

\begin{Prf}
Let $\pi$ denote the representative of $\Sigma_{14}$ given in Table~\ref{table:main}, namely
\begin{equation} \label{Sigma14matrix-stab}
\left[ \begin{matrix}
\alpha & \beta & \cdot \\
\beta & c\gamma & \beta-\gamma \\
\cdot & \beta-\gamma & \gamma
\end{matrix} \right],
\end{equation}
where $c$ satisfies the condition ($\dagger$) given in Table~\ref{table:main}. 
That is, $c$ is some fixed element of $\bF_q \setminus \{0,1\}$ such that (i) $-3c$ is a square in $\mathbb{F}_q$ and (ii) $f(c):=\tfrac{\sqrt{c}+1}{\sqrt{c}-1}$ is a non-cube in $\mathbb{F}_q(\sqrt{-3})$. 

We first claim that $c \neq -1$. 
To verify this, begin by noting that $a=\sqrt{-1}$ is a solution of the equation $f(-1)=a^3$ in every field in which $-1$ is a square. 
Now suppose for a contradiction that $c = -1$. 
Then (i) says that $3$ is a square in $\mathbb{F}_q$. 
If $q \equiv 1 \pmod 3$ then $-3$ is a square in $\mathbb{F}_q$, so $-1$ must also be a square $\mathbb{F}_q$, and so (ii) does not hold because $f(-1)$ has the cube root $\sqrt{-1}$ in $\mathbb{F}_q(\sqrt{-3})=\mathbb{F}_q$. 
If $q \equiv -1 \pmod 3$ then $-3$ is a non-square in $\mathbb{F}_q$, so $-1$ must be a non-square in $\mathbb{F}_q$, but in this case $-1$ is a square in $\mathbb{F}_q(\sqrt{-3})$, so again (ii) does not hold. 

Recall that point-orbit distribution of $\pi$ is $[1, \frac{q\mp 1}{2}, \frac{q\mp 1}{2}, q^2\pm 1]$ according to whether $q \equiv \pm 1 \pmod 3$. 
The points of rank at most $2$ in $\pi$ lie on the cubic $\mathcal{C} : \alpha f_c(\beta,\gamma)-\beta^2\gamma = 0$, where $f_c(\beta,\gamma) := (c-1)\gamma^2+2\beta\gamma-\beta^2$. 
The unique point of rank $1$ is $x : \beta=\gamma=0$.

Consider an element $h$ of $K$, represented by a matrix $g=(g_{ij}) \in \GL(3,q)$. 
If $h$ fixes $\pi$ then, in particular, it must fix $x$, and so $g$ must have the form
\begin{equation} \label{sig14stab-x}
g = \left[ \begin{matrix}
g_{11} & g_{12} & g_{13} \\
\cdot & g_{22} & g_{23} \\
\cdot & g_{32} & g_{33}
\end{matrix} \right]. 
\end{equation}
Now consider the rank-$2$ points $w$ and $w'$ with $(\alpha,\beta,\gamma)=(0,1,0)$ and $(c^{-1},1,1)$, respectively, and let $M_w$ and $M_w'$ denote the corresponding matrix representatives obtained from \eqref{Sigma14matrix-stab}.
Note that $\pi$ is spanned by $x$, $w$ and $w'$. 
Hence, given that $g$ has the form \eqref{sig14stab-x}, $h$ now fixes $\pi$ if and only if it maps both $w$ and $w'$ to points in $\pi$, that is, if and only if $gM_wg^T$ and $gM_w'g^T$ have the form \eqref{Sigma14matrix-stab} for some $\alpha$, $\beta$ and $\gamma$ (and are of rank $2$). 

We first claim that $g_{33} \neq 0$. 
For a contradiction, suppose the contrary, and note that the supposedly non-zero determinant of $g$ is then $-g_{11}g_{23}g_{32}$. 
Consider the image of $w$ under $h$. 
The $(2,2)$-entry of $gM_wg^T$ must equal $c$ times the $(3,3)$-entry; the $(3,3)$-entry is zero, and the $(2,2)$-entry is $2g_{22}g_{23}$, so we must have $g_{22}=0$ because $g_{23} \neq 0$. 
The $(1,3)$-entry of $gM_wg^T$ is $(g_{11}+g_{13})g_{32}$ and must vanish, so $g_{13}=-g_{11}$. 
Finally, the $(2,3)$-entry of $gM_wg^T$, which is $g_{23}g_{32}$, must equal the sum of the $(1,2)$-entry, which is now $g_{12}g_{23}$, and the $(3,3)$-entry, which is zero. 
This implies that $g_{12}=g_{32}$. 
At this point, the image of $w$ under $h$ lies in $\pi$ as required. 
Now consider the image of $w'$ under $h$. 
The $(1,3)$-entry of $gM_w'g^T$, which is $g_{32}(g_{11}+cg_{32})$, must vanish, so $g_{11}=-cg_{32}$ because $g_{32} \neq 0$. 
The $(1,2)$- and $(3,3)$-entries of $gM_w'g^T$ sum to $cg_{32}(g_{23}+g_{32})$ and this sum must equal the $(2,3)$-entry, which is zero, so $g_{23} = -g_{32}$. 
At this point, the $(2,2)$-entry of $gM_w'g^T$ is $g_{32}^2$. 
It must equal $c$ times the $(3,3)$-entry. 
The $(3,3)$-entry is $cg_{32}^2$, so we must have $c^2=1$ because $g_{32} \neq 0$. 
However, $c \neq 1$ by condition ($\dagger$), and $c \neq -1$ as shown above, so we have reached a contradiction. 
Therefore, $g_{33} \neq 0$ as claimed. 

Knowing now that $g_{33} \neq 0$, consider again the image of $w$ under $h$. 
The $(2,2)$-entry of $gM_wg^T$, which is $2g_{22}g_{23}$, must equal $c$ times the $(3,3)$-entry. 
The $(3,3)$-entry is $2g_{32}g_{33}$, so we must have $g_{32} = c^{-1}g_{33}^{-1}g_{22}g_{23}$. 
Setting the $(1,3)$-entry of $gM_gg^T$ to zero then implies that $g_{12} = -c^{-1}g_{33}^{-2}(g_{11}+g_{13})g_{22}g_{23}$. 
Note that at this point $\text{det}(g)=g_{11}g_{22}g_{33}^{-1}(cg_{33}^2-g_{23}^2)$, so in particular $cg_{33}^2-g_{23}^2 \neq 0$. 
Now consider the image of $w'$ under $h$. 
The requirement that the $(2,2)$-entry be equal to $c$ times the $(3,3)$-entry implies that $(g_{22}^2-g_{33}^2)(cg_{33}^2-g_{22}^2)=0$, and so we must have $g_{22} = \pm g_{33}$. 
Setting the $(1,3)$-entry of $gM_w'g^T$ then yields $g_{13} = g_{11}g_{23}(g_{23} \mp g_{33})(cg_{33}^2-g_{22}^2)^{-1}$ according to whether $g_{22} = \pm g_{33}$. 
In summary, we now have $g=g_\pm$ where
\[
g_\pm := \left[ \begin{matrix}
g_{11} & \mp c^{-1}g_{11}g_{23}(cg_{33} \mp g_{23})(cg_{33}^2-g_{23}^2)^{-1} & g_{11}g_{23}(g_{23} \mp g_{33})(cg_{33}^2-g_{23}^2)^{-1} \\
\cdot & \pm g_{33} & g_{23} \\
\cdot & \pm c^{-1}g_{23} & g_{33}
\end{matrix} \right]. 
\]
It remains to satisfy the condition that the sum of the $(1,2)$- and $(3,3)$-entries should be equal to the $(2,3)$-entry in both of $gM_wg^T$ and $gM_w'g^T$, for each of $g=g_\pm$. 
Consider first $g=g_+$. 
We deduce from considering the images of $w$ and $w'$ respectively that
\begin{align*}
cg_{33}(g_{11}-g_{33}) - g_{23}(g_{11}+2g_{33}) - g_{23}^2 &= 0, \\
cg_{33}(g_{11}-g_{33}) - cg_{23}(g_{11}+2g_{33}) - g_{23}^2 &= 0.
\end{align*}
Summing these equations gives $(c-1)g_{23}(g_{11}+2g_{33})=0$, and we recall that $c-1 \neq 0$ by condition ($\dagger$). 
If $g_{23}=0$ then $g=g_+$ is a scalar multiple of the identity matrix, so to obtain a non-trivial element of $K_\pi$ we must have $g_{11} = -2g_{33}$. 
Putting this into either of the above equations then implies that $g_{32}^2 = -3c g_{33}^2$. 
Since $-3c$ is a square, by condition ($\dagger$), we therefore have $g_{32} = \pm \sqrt{-3c} g_{33}$. 
In summary, under the assumption that $g$ is of the form $g_+$, we have concluded that $g$ is either a scalar matrix or a matrix of the form
\[
A_\pm(\xi) := \left[ \begin{matrix}
-2\xi & \pm 2c^{-2} (c\xi \mp \sqrt{-3c} \;\xi) \sqrt{-3c} & \pm 2c^{-1} (\xi \mp \sqrt{-3c} \;\xi) \sqrt{-3c} \\
\cdot & \xi & \pm \sqrt{-3c} \;\xi \\
\cdot & \pm c^{-1} \sqrt{-3c} \;\xi & \xi
\end{matrix} \right] 
\]
for some non-zero $\xi \in \mathbb{F}_q$. 
The matrices $A_\pm(\xi)$ have order $3$ modulo scalars, so we obtain a subgroup of $K_\pi$ of order $3$. 
It remains to show that this subgroup is the whole of $K_\pi$, that is, that no further non-trivial elements arise under the assumption that $g$ instead has the form $g_-$. 

Suppose therefore that $g=g_-$. 
Consider $gM_w'g^T$ first. 
Setting the sum of the $(1,2)$- and $(3,3)$-entries equal to the $(2,3)$-entry gives us the constraint
\[
g_{11}(g_{23}+g_{33}) + c^{-1}(cg_{33}^2+g_{23}^2) + 2g_{23}g_{33} = 0.
\]
Notice that we cannot have $g_{23}+g_{33}=0$, because then the above equation would imply that $c=1$, in contradiction with condition ($\dagger$). 
Therefore, $g_{11} = -c^{-1}(g_{23}+g_{33})^{-1}(cg_{33}^2+g_{23}^2+2cg_{23}g_{33})$. 
Considering now the same condition in $gM_wg^T$, we find that $g_{23}$ and $g_{33}$ must satisfy
\begin{equation} \label{sig14finalCubic}
(c+1)g_{23}^3 + 6cg_{23}^2g_{33} + 3c(c+1)g_{23}g_{33}^2 + 2c^2g_{33}^3 = 0.
\end{equation}
Notice that neither of $g_{23}$ or $g_{33}$ can be zero, because if either one is zero then \eqref{sig14finalCubic} implies that both are zero because $c \neq \{0,-1\}$, and this contradicts $\text{det}(g) \neq 0$. 
Hence, we may view the left-hand side of \eqref{sig14finalCubic} as a cubic in $\theta:=g_{23}g_{33}^{-1}$. 
Upon changing the parameter from $\theta$ to $\theta-\tfrac{2c}{c+1}$ and dividing by $c+1$, this cubic becomes 
\[
\theta^3 + \frac{3c(c-1)^2}{(c+1)^2} \theta - \frac{4c^2(c-1)^2}{(c+1)^3}.
\]
The discriminant of the above cubic is $-108c^3(c-1)^4 = 4(-3c)^3(c-1)^4$, which is a square in $\mathbb{F}_q$ because $-3c$ is a square, by condition ($\dagger$). 
By \cite[Theorem~3]{DicksonCubics}, the cubic therefore has a solution in $\mathbb{F}_q$ if and only if 
\[
\frac{1}{2} \left( \frac{4c^2(c-1)^2}{(c+1)^3} + \sqrt{-3} \cdot \frac{2c(c-1)^2\sqrt{-3c}}{3(c+1)^2} \right) 
= 
\frac{(-\sqrt{c})^3 (c-1)^3}{(c+1)^3} \cdot \frac{\sqrt{c}-1}{\sqrt{c}+1}
\]
is a cube in $\mathbb{F}_q(\sqrt{-3})$. 
By condition ($\dagger$), this is not the case. 
Therefore, there are no non-trivial elements of $K_\pi$ represented by matrices of the form $g=g_-$, and so $K_\pi \cong C_3$ as claimed. 
\end{Prf}

\begin{La} \label{stab14dash}
Suppose that $q \equiv 0 \pmod 3$. 
The stabiliser of a plane of type $\Sigma_{14}'$ is isomorphic to $E_q : C_{q-1}$. 
The number of planes of type $\Sigma_{14}'$ is therefore $q^2(q^3-1)(q+1)$.
\end{La}

\begin{Prf}
Let $\pi$ denote the representative of $\Sigma_{14}'$ used in the proof of Lemma~\ref{dist14dash} (as opposed to the one given in Table~\ref{table:main}). 
That is, suppose that $\pi$ is represented by the matrix
\[
M = \left[ \begin{matrix} 
-\beta & \cdot & \gamma-\beta \\ 
\cdot & \alpha+\beta & \beta \\ 
\gamma-\beta & \beta & \cdot 
\end{matrix} \right].
\]
By Lemma~\ref{dist11}, $\pi$ contains a unique line of type $o_9$, specifically the line $\ell : \beta=\gamma$. 
Hence, we have $K_\pi \le K_\ell$. 
The matrix representatives of elements of $K_\ell$ are precisely the matrices $g = (g_{ij}) \in \GL(3,q)$ of the form 
\[
\left[ \begin{matrix} 
g_{11} & \cdot & g_{21}g_{33}g_{11}^{-1} \\ 
g_{21} & g_{11}^2g_{33}^{-1} & g_{23} \\ 
\cdot & \cdot & g_{33} 
\end{matrix} \right].
\]
The points of rank at most $2$ in $\pi$ comprise the cubic $\mathcal{C} : \beta\gamma(\beta+\gamma) + \alpha(\beta-\gamma)^2 = 0$, which meets the line $\ell$ in the unique point $x : \beta=\gamma=0$ of rank $1$ in $\pi$. 
Therefore, $K_\pi = K_\ell \cap K_{\mathcal{C}'}$, where $\mathcal{C}' = \mathcal{C} \setminus \{x\}$. 
A matrix $g \in \GL(3,q)$ of the above form fixes the matrix representative of $\mathcal{C}'$ if and only if $g_{23} =  -g_{21}^2g_{33}(g_{11}^{-1})^2$. 
Factoring out scalars gives $K_\pi \cong E_q : C_{q-1}$. 
\end{Prf}

\begin{La} \label{stab15}
The stabiliser of a plane of type $\Sigma_{15}$ is isomorphic to $E_q^{1+2}:C_{q-1}$, where $E_q^{1+2}$ has centre $Z \cong E_q$ and $E_q^{1+2}/Z \cong E_q^2$. 
The number of planes of type $\Sigma_{15}$ is $(q^3-1)(q+1)$.
\end{La}

\begin{Prf}
Let $\pi$ denote the representative of $\Sigma_{15}$ given in Table~\ref{table:main}. 
That is, suppose that $\pi$ is represented by the matrix
\[
M = \left[ \begin{matrix} 
\alpha & \beta & \gamma \\ 
\beta & \gamma & \cdot \\ 
\gamma & \cdot & \cdot 
\end{matrix} \right].
\]
As noted in the proof of Lemma~\ref{dist15}, the points of rank at most $2$ in $\pi$ are all on the line $\ell : \gamma=0$. 
This happens to be the same line of type $o_6$ as in the proof of Lemma~\ref{stab9}, so we know that the matrix representatives of elements of $K_\ell$ are the matrices $g = (g_{ij}) \in \GL(3,q)$ with $g_{21}=g_{31}=g_{32}=0$. 
An element of $K_\ell$ fixes $\pi$ if and only if it fixes the set of rank-$3$ points in $\pi$ (i.e. those with $\gamma \neq 0$ in $M$). 
This imposes the additional constraint $g_{11} = g_{33}^2g_{22}^2$, so that the matrix representatives of elements of $K_\pi$ are precisely the invertible matrices of the form
\[
\left[ \begin{matrix} 
g_{33}^{-1}g_{22} & g_{12} & g_{13} \\ 
\cdot & g_{22} & g_{23} \\ 
\cdot & \cdot & g_{33}
\end{matrix} \right].
\]
Note that setting $g_{22}=g_{33}=1$ yields a subgroup of shape $E_q^{1+2}$, with centre defined by the constraint $g_{12}=g_{23}=0$. 
Hence, upon factoring out scalars we obtain $K_\pi \cong E_q^2 : C_{q-1}$. 
\end{Prf}

\section*{Acknowledgements}

The first author acknowledges the support of {\em The Scientific and Technological Research Council of Turkey} T\"UB\.{I}TAK (project no.~118F159).


\begin{thebibliography}{1}

\bibitem{AlLa2019}
N~Alnajjarine and M.~Lavrauw, ``Determining the rank of tensors in $\mathbb{F}_q^2 \otimes \mathbb{F}_q^3 \otimes \mathbb{F}_q^3$'', in: D.~Slamanig, E.~Tsigaridas and Z.~Zafeirakopoulos (eds.), MACIS 2019: Mathematical Aspects of Computer and Information Sciences, {\em Lecture Notes in Computer Science} {\bf 11989}, Springer, Cham, 2020. 

\bibitem{Campbell1927}
A.~D.~Campbell, ``Pencils of conics in the Galois fields of order $2^n$'', 
\newblock {\em Amer. J. Math.} {\bf 49} (1927) 401--406.

\bibitem{Casas-Alvero2010}
E. Casas-Alvero. {\em Analytic Projective Geometry}, European Mathematical Society, Z\"urich, 2014.

\bibitem{DicksonCubics}
L.~E.~Dickson, ``Criteria for the irreducibility of functions in a finite field'', {\em Bull. Amer. Math. Soc.} {\bf 13} (1906) 1--8.

\bibitem{Dickson1908}
L.~E.~Dickson, ``On families of quadratic forms in a general field'', 
\newblock {\em Quarterly J. Pure Appl. Math.} {\bf 45} (1908) 316--333.

\bibitem{Jordan1906}
C. Jordan, ``R\'eduction d'un r\'eseau de formes quadratiques ou bilin\'eaires: premi\`ere partie",
\newblock {\em Journal de math\'ematiques pures et appliqu\'ees} (1906) 403--438.

\bibitem{Jordan1907}
C. Jordan: ``R\'eduction d'un r\'eseau de formes quadratiques ou bilin\'eaires: deuxi\`eme partie", Gauthier-Villars, 
\newblock {\em Journal de math\'ematiques pures et appliqu\'ees} (1907) 5--51.

\bibitem{LaPo2020}
M.~Lavrauw and T.~Popiel, ``The symmetric representation of lines in $\text{PG}(\mathbb{F}_q^3 \otimes \mathbb{F}_q^3)$, {\em Discrete Math.} {\bf 343} (2020) 111775.

\bibitem{LaPoSh2020}
M.~Lavrauw, T.~Popiel and J.~Sheekey, ``Nets of conics of rank one in $\text{PG}(2,q)$, $q$ odd'', {\em J. Geometry} {\bf 111} (2020) 36.

\bibitem{LaSh2015}
M.~Lavrauw and J.~Sheekey, ``Canonical forms of $2 \times 3 \times 3$ tensors over the real field, algebraically closed fields, and finite fields'', {\em Linear Algebra Appl.} {\bf 476} (2015) 133--147.

\bibitem{Wall1977}
C.~T.~C.~Wall, ``Nets of conics''. {\em Math. Proc. Cambridge Philos. Soc.} {\bf 81} (1977) 351--364.

\bibitem{Wilson1914}
A.~H.~Wilson, ``The canonical Types of Nets of Modular Conics", {\em Amer. J. Math.} {\bf 36} (1914) 187--210.

\end{thebibliography}
\end{document}